\documentclass{gtart_h}

\def\ifplaintex{\expandafter\ifx\csname documentclass\endcsname\relax}

\def\gtp{{\mathsurround=0pt\it $\cal G\mskip-2mu$eometry \&\ 
$\cal T\!\!$opology $\cal P\!$ublications}}  

\def\recd{{\small Received:\qua\receiveddate\ifx\reviseddate\relax
\else\qquad Revised:\qua\reviseddate\fi\par}} 

\def\lognumber#1{\def\thelognumber{#1}}
\def\volumenumber#1{\def\thevolumenumber{#1}}
\def\volumeyear#1{\def\thevolumeyear{#1}}
\def\papernumber#1{\def\thepapernumber{#1}}
\def\pagenumbers#1#2{\def\startpage{#1}\def\finishpage{#2}}
\def\published#1{\def\publishdate{#1}}

\def\received#1{\def\receiveddate{#1}}

\def\accepted#1{\def\accepteddate{#1}}
\def\asciititle#1{\def\theasciititle{#1}}

\def\asciiaddress#1{\def\theasciiaddress{#1}}
\def\asciiemail#1{\def\theasciiemail{#1}}

\long\def\asciiabstract#1{\long\def\theasciiabstract{#1}}
\def\asciikeywords#1{\def\theasciikeywords{#1}}

\let\\\par\let\thelognumber\relax\let\thevolumenumber\relax
\let\thepapernumber\relax\let\thevolumeyear\relax\let\startpage\relax
\let\finishpage\relax\let\publishdate\relax\let\receiveddate\relax
\let\reviseddate\relax\let\accepteddate\relax\let\theasciititle\relax
\let\theasciiauthors\relax\let\theasciiaddress\relax
\let\theasciiabstract\relax\let\theasciikeywords\relax

\let\theasciiemail\relax

\ifplaintex
\font\logobig=cmssbx10 scaled 3836
\font\logomed=cmssbx10 scaled 2557
\else
\font\logobig=cmssbx10 scaled 4200
\font\logomed=cmssbx10 scaled 2800
\fi

\long\def\makeagttitle{   
\count0=\startpage
\agt\hfill      
\hbox to 45truept{\vbox to 0pt{\vglue -13truept{\logomed A\kern -.37em{\logobig 
T}\kern -.38em G}\vss}\hss}
\break
{\small Volume \thevolumenumber\ (\thevolumeyear)
\startpage--\finishpage\nl
Published: \publishdate}

\vglue .25truein

{\parskip=0pt\leftskip 0pt plus
1fil\def\\{\par\smallskip}{\Large\bf\thetitle}\par\medskip} \vglue
0.05truein

{\parskip=0pt\leftskip 0pt plus 1fil\def\\{\par}{\sc\theauthors}
\par\medskip}%
 
\vglue 0.03truein 

{\small\leftskip 25truept\rightskip 25truept{\bf Abstract}\stdspace\theabstract

{\bf AMS Classification}\stdspace\theprimaryclass
\ifx\thesecondaryclass\relax\else; \thesecondaryclass\fi\par
{\bf Keywords}\stdspace \thekeywords\par}\vglue 7truept

}   

\ifplaintex
\font\phead=cmsl9 scaled 950
\font\pnum=cmbx10 scaled 913
\font\pfoot=cmsl9 scaled 950
\headline{\vbox to 0pt{\vskip -4.5mm\line{\small\phead\ifnum
\count0=\startpage ISSN 1472-2739 (on-line) 1472-2747 (printed)
\hfill {\pnum\folio}\else\ifodd\count0\def\\{ }%
\ifx\theshorttitle\relax\thetitle\else\theshorttitle\fi\hfill{\pnum\folio}
\else\def\\{ and }{\pnum\folio}\hfill\ifx\theshortauthors\relax\theauthors
\else\theshortauthors\fi\fi\fi}\vss}}
\footline{\vbox to 0pt{\vglue 0mm\line{\small\pfoot\ifnum\count0=\startpage
\copyright\ \gtp\hfill\else
\agt, Volume \thevolumenumber\ (\thevolumeyear)\hfill\fi}\vss}}
\else
\font=cmsl9 scaled 1050
\font\lnum=cmbx10 
\font=cmsl9 scaled 1050
\makeatletter
\def\@oddhead{{\small\lhead\ifnum\count0=\startpage ISSN 1472-2739 
(on-line) 1472-2747 (printed)\hfill {\lnum\number\count0}\else\ifodd\count0
\def\\{ }\ifx\theshorttitle\relax \thetitle \else\theshorttitle\fi\hfill
{\lnum\number\count0}\else\def\\{ and }{\lnum\number\count0}
\hfill\ifx\theshortauthors\relax 
\theauthors\else\theshortauthors\fi\fi\fi}}\def\@evenhead{\@oddhead}
\def\@oddfoot{\small\lfoot\ifnum\count0=\startpage\copyright\ \gtp\hfill\else
\agt, Volume \thevolumenumber\ (\thevolumeyear)\hfill\fi}
\def\@evenfoot{\@oddfoot}
\makeatother
\fi
\let\maketitlepage\makeagttitle
\let\makeshorttitle\maketitlepage
\let\maketitle\maketitlepage

\newwrite\gtoutfile
\long\gdef\makeheadfile{  
{\def\\{, }\def\s{ }
\immediate\openout\gtoutfile head.xxx
\immediate\write\gtoutfile{Proxy-for: \ifx\theasciiauthors\relax
\theauthors\else\theasciiauthors\fi\s<\ifx\theasciiemail\relax\theemail\else\theasciiemail\fi>}
\immediate\write\gtoutfile{\noexpand\\}
\immediate\write\gtoutfile{Authors: \ifx\theasciiauthors\relax
\theauthors\else\theasciiauthors\fi}
{\def\\{ }\immediate\write\gtoutfile{Title: \ifx\theasciititle\relax
\thetitle\else\theasciititle\fi}}
\immediate\write\gtoutfile{Subj-class: GT or SG, GR etc}
\immediate\write\gtoutfile{MSC-class: \theprimaryclass\ifx\thesecondaryclass\relax\else, \thesecondaryclass\fi}
\immediate\write\gtoutfile{Journal-ref: Algebr. Geom. Topol. \thevolumenumber\s
(\thevolumeyear) \startpage-\finishpage}
\immediate\write\gtoutfile{Comments: Published by Algebraic and
Geometric Topology at}
\immediate\write\gtoutfile{\s\s\s  http://www.maths.warwick.ac.uk/agt/AGTVol\thevolumenumber/agt-\thevolumenumber-\thepapernumber.abs.html}
\immediate\write\gtoutfile{\noexpand\\}
\immediate\write\gtoutfile{}
\ifx\theasciiabstract\relax
\immediate\write\gtoutfile{\theabstract}\else
\immediate\write\gtoutfile{\theasciiabstract}\fi
\immediate\write\gtoutfile{}
\immediate\write\gtoutfile{\noexpand\\}
\immediate\write\gtoutfile{}
\immediate\closeout\gtoutfile}}  

\def\maketitlepage{\makeagttitle\makeheadfile}
\let\makeshorttitle\maketitlepage
\let\maketitle\maketitlepage

\lognumber{5}
\volumenumber{5}
\volumeyear{2005}
\papernumber{5}
\pagenumbers{71}{106}
\received{9 February 2004} 
\accepted{28 December 2004}
\published{3 February 2005}

\usepackage{pstricks, amssymb, xr, graphicx, subfigure, amsmath}
\input epsf
\def\relabelbox{%
  \hbox\bgroup%
}%
\def\endrelabelbox{%
\egroup%
}%
\def\relabel #1#2 {%
  \special{ps:/a {} def}%
  \smash{\rlap{#2}}%
}%
\def\adjustrelabel <#1,#2> #3#4 {%
  \special{ps:/a {} def}%
  \smash{\rlap{\kern #1 \raise #2\hbox{#4}}}%
}%
\def\extralabel <#1,#2> #3 {\smash{\rlap{\kern #1 \raise #2\hbox{#3}}}}%

\def\hfl#1{\smash{\mathop{\hbox to 10mm{\rightarrowfill}}\limits^{\textstyle
#1}}}

\newtheorem{thm}{Theorem}[section]
\newtheorem{lemma}[thm]{Lemma}
\newtheorem{sublemma}[thm]{Sublemma}
\newtheorem{proposition}[thm]{Proposition}
\newtheorem{corollary}[thm]{Corollary}

\theoremstyle{definition}

\newtheorem{remark}[thm]{Remark}
\newtheorem{notation}[thm]{Notation}
\newtheorem{definition}[thm]{Definition}

\newcommand{\IT}{\mbox{\rm Int}}

\newcommand{\CA}{{\cal A}}

\newcommand{\CF}{{\cal F}}
\newcommand{\CG}{{\cal G}}
\newcommand{\CH}{{\cal H}}
\newcommand{\CI}{{\cal I}}

\newcommand{\CL}{{\cal L}}

\newcommand{\ZZ}{\mathbb{Z}}
\newcommand{\RR}{\mathbb{R}}

\newcommand{\NN}{\mathbb{N}}
\newcommand{\bp}{\begin{proof}}
\newcommand{\eop}{\end{proof}}

\begin{document}
\title[Clover calculus via basic algebraic topology]{Clover calculus for homology $3$--spheres\\via basic algebraic topology}
\asciititle{Clover calculus for homology 3-spheres via 
basic algebraic topology}
\authors{Emmanuel Auclair\\Christine Lescop}
\address{Institut Fourier (UMR 5582 du CNRS), B.P. 74\\38402 Saint-Martin d'H\`eres cedex, France}

\asciiaddress{Institut Fourier (UMR 5582 du CNRS), B.P. 74\\38402 Saint-Martin d'Heres cedex, France}

\asciiemail{auclaire@ujf-grenoble.fr, lescop@ujf-grenoble.fr}
\gtemail{\mailto{auclaire@ujf-grenoble.fr}{\rm\qua 
and\qua}\mailto{lescop@ujf-grenoble.fr}}
\urladdr{http://www-fourier.ujf-grenoble.fr/~lescop }

\begin{abstract}
We present an alternative definition for the Goussarov--Habiro filtration of the $\ZZ$--module freely generated by oriented integral homology $3$--spheres, by means of Lagrangian-preserving homology handlebody replacements ({\CL}P--surgeries).
Garoufalidis, Goussarov and Polyak proved that the graded space $(\CG_n)_n$ associated to this filtration is generated by Jacobi diagrams. Here, we express elements associated to 
{\CL}P--surgeries as explicit combinations of these Jacobi diagrams in $(\CG_n)_n$.
The obtained coefficient in front of a Jacobi diagram is computed like its weight system
with respect to a Lie algebra equipped with a non-degenerate invariant bilinear form, where cup products in $3$--manifolds play the role of the Lie bracket and the linking number replaces the invariant form.
In particular, this article provides an algebraic version of the graphical clover calculus developed by Garoufalidis, Goussarov, Habiro and Polyak. This version induces splitting formulae for all finite type invariants of homology $3$--spheres.
\end{abstract}

\asciiabstract{%
We present an alternative definition for the Goussarov--Habiro
filtration of the Z-module freely generated by oriented integral
homology 3-spheres, by means of Lagrangian-preserving homology
handlebody replacements (LP-surgeries).  Garoufalidis, Goussarov and
Polyak proved that the graded space (G_n)_n associated to this
filtration is generated by Jacobi diagrams. Here, we express elements
associated to LP-surgeries as explicit combinations of these Jacobi
diagrams in (G_n)_n.  The obtained coefficient in front of a Jacobi
diagram is computed like its weight system with respect to a Lie
algebra equipped with a non-degenerate invariant bilinear form, where
cup products in $3$--manifolds play the role of the Lie bracket and
the linking number replaces the invariant form.  In particular, this
article provides an algebraic version of the graphical clover calculus
developed by Garoufalidis, Goussarov, Habiro and Polyak. This version
induces splitting formulae for all finite type invariants of homology
3-spheres.}

\asciikeywords{3-manifolds, homology spheres, finite type invariants, 
Jacobi diagrams, Borromeo surgery, clover calculus, clasper calculus, Goussarov-Habiro filtration}

\maketitle

\section{Introduction}

In 1995, in \cite{oht}, Tomotada Ohtsuki introduced a notion of finite
type invariants for homology $3$--spheres (that are compact oriented
$3$--manifolds with the same homology with integral coefficients as
the standard $3$--sphere $S^3$), following the model of the theory of
Vassiliev invariants for knots in the ambient space $\RR^3$.  He
defined a filtration of the real vector space freely generated by
homology $3$--spheres and began the study of the associated graded
space. In \cite{le}, Thang Le finished identifying this graded space
to a space of Jacobi diagrams called $\CA_{\RR}(\emptyset)$. The
Jacobi diagrams, precisely defined in \mbox{Subsection \ref{S21}}, are
represented by trivalent finite graphs with additional orientation
information.

Similar filtrations of the $\ZZ$--module freely generated by homology $3$--spheres and 
their relationships have been studied by Garoufalidis, Goussarov, Polyak and others. See \cite{ggp} and references therein.
Over $\ZZ[1/2]$, all of them are equivalent to the Ohtsuki filtration \cite{ggp}.

Among these filtrations, the most convenient one is the Goussarov--Habiro one 
where the Matveev Borromeo surgeries \cite{mat} (defined in \mbox{Subsection \ref{S23}}) 
play the role of the crossing changes in the knot case. It allowed Garoufalidis, 
Goussarov and Polyak to define a set of generators $\Psi_n(\Gamma)$ for the degree $n$ 
part $\CG_n$ of the associated
Goussarov--Habiro graded $\ZZ$--module, for Jacobi diagrams $\Gamma$ with at most $n$ 
vertices \cite{ggp}. See Subsection \ref{SII3}. Garoufalidis, Goussarov and Polyak also 
gave some graphical rules that allow one to reduce an element to a combination of 
their generators. This set of rules is the so-called {\em clover calculus.\/}
Here, these rules are enclosed in two propositions \ref{Review} and \ref{ihx}. 

Our main theorem \ref{TT} expresses elements of $\CG_n$
associated to the {\em {\CL}P--surgeries \/}  defined in Subsection \ref{S22}, 
as explicit combinations of the $\Psi_n(\Gamma)$, in terms of intersection forms
(or cup products) and 
linking numbers. Therefore, this article presents a completely algebraic version of the 
Garoufalidis--Goussarov--Habiro--Polyak clover calculus. 
Furthermore, it tightens the links between Jacobi diagrams and topology by relating 
the vertices of the Jacobi diagrams to cup products in $3$--manifolds and the 
diagram edges with linking numbers.

We also give an alternative definition of the Goussarov--Habiro filtration of 
the $\ZZ$--module of integral homology $3$-spheres, by means of {\CL}P--surgeries. 
See Corollary \ref{nouvelledef}.

Let us now give a slightly more specific description of our main
theorem \ref{TT}.

A {\em homology genus $g$ handlebody\/}
is an oriented compact $3$--manifold with the same integral homology as the standard  genus $g$
handlebody $H_g$. The boundary $\partial A$ of such a manifold $A$ is then homeomorphic 
to the genus $g$ surface $\partial H_g$. The {\em Lagrangian\/} $\CL_A$ of $A$
is the kernel of the map induced by the inclusion from $H_1(\partial A;\ZZ)$ to 
$H_1(A;\ZZ)$. A {\em Lagrangian-preserving surgery\/} or {\em {\CL}P--surgery\/} on a 
homology sphere $M$
consists in removing the interior of such a homology handlebody $(A \subset M)$ and 
replacing it by another such $B$ whose boundary $\partial B$ is identified to 
$\partial A$ so that $\CL_A=\CL_B$. 

In our definition of the Goussarov--Habiro filtration $(\CF_n)_{n \in \NN}$ of the $\ZZ$--module $\CF=\CF_0$ freely generated by the oriented homology spheres up to orientation-preserving homeomorphisms, the $n^{th}$ module $\CF_n$ is 
generated by {\em brackets\/} $[D]$ of so-called {\em $n$--component {\CL}P--surgeries\/}
$D$ that are made of $n$
disjoint {\CL}P--surgeries $(A_i,B_i)$ in $M$. (The $A_i$ are disjoint in $M$.) 
$$[D]= \sum_{J \subset \{1,\dots,n\}}(-1)^{\sharp J} 
\left( \left(M \setminus \sqcup_{j \in J}\mbox{Int}(A_j)\right) \cup_{\partial} \left(\sqcup_{j \in J}B_j\right)\right).$$
Our main result expresses the bracket $[D]$ of an $n$--component {\CL}P--surgery $D$ as 
$$[D]=\sum_{\Gamma}\ell(D;\Gamma)\Psi_n(\Gamma)$$
in $\CG_n=\CF_n/\CF_{n+1}$, where the coefficient $\ell(D;\Gamma)$ of $\Psi_n(\Gamma)$
is an explicit function of the cup products in the manifolds $(A_i \cup -B_i)$, 
of the
linking pairings on $H_1(A_i) \otimes H_1(A_j)$, $i \neq j$, and of variations of the Rohlin invariant when replacing $A_i$ by $B_i$.

Let us roughly define $\ell(D;\Gamma)$ when $n$ is the number of vertices of $\Gamma$
and when $\Gamma$ admits no non-trivial automorphism. The general definition of
$\ell(D;\Gamma)$ is given in Subsection \ref{subsdeflinkLP}.
When a bijection $\sigma$ from the set of vertices of $\Gamma$ to $\{1, \dots, n\}$ 
is given, the algebraic intersection of surfaces (or the cup product) of each $(A_i \cup -B_i)$ is placed
at the vertex $\sigma^{-1}(i)$. The cup products are next contracted along the 
edges with respect to the linking pairing to produce a number $\ell(D;\Gamma;\sigma)$, 
and $\ell(D;\Gamma)=\sum_{\sigma}\ell(D;\Gamma;\sigma)$. 
This construction is similar to the construction of weight systems associated to 
Lie algebras.

The proof of Theorem \ref{TT} goes as follows. We first prove that the standard 
Goussarov, Garoufalidis and Polyak generators have appropriate coefficients in 
Subsection \ref{SSJD}. Then we use the similarities
between the behaviour of the bracket in $\CG_n$ and the behaviour of our 
coefficients to reduce the proof to this former case.

Though this article is largely inspired by \cite{ggp}, it is written in a self-contained 
way
in an
attempt to replace all the graphical arguments in \cite{ggp} by more intrinsic 
arguments of geometric or algebraic topology.

Theorem \ref{TT} can be used to derive formulae on the behaviour 
under {\CL}P--surgeries of all finite-type invariants of homology spheres in the 
Goussarov--Habiro sense.
For example, it immediately leads to splitting formulae for the restriction to 
homology spheres of the Kontsevich--Kuperberg--Thurston universal finite-type invariant $Z_{KKT}$ of rational homology spheres. 
In \cite{L2}, the second author proved that these
formulae generalise to rational homology spheres and to rational homology 
handlebody replacements that preserve the rational Lagrangians. 
These generalized splitting formulae are fairly 
easy to guess
from the Kontsevich--Kuperberg--Thurston construction
(but much harder to prove in general), they actually led the second author to
the formulae of Theorem \ref{TT}.
These formulae had been previously noticed by
G.~Kuperberg and D.~Thurston in the special case
where rational homology handlebodies are reglued by a homeomorphism that induces the 
identity in homology \cite{KT}. This special case is sufficient to prove 
that $Z_{KKT}$ is universal among finite-type invariants of homology spheres. 

The second author thanks Thang Le, Gregor Masbaum and Dylan Thurston for useful and 
pleasant conversations.

\section{Background}
\setcounter{equation}{0}

\subsection{Jacobi diagrams} \label{S21}

In what follows, a {\em  Jacobi diagram\/} $\Gamma$ is a trivalent graph  
without simple loop like $\begin{pspicture}[.2](0,0)(.6,.4)
\psline{-*}(0.05,.2)(.25,.2)
\pscurve{-}(.25,.2)(.4,.05)(.55,.2)(.4,.35)(.25,.2)
\end{pspicture}$.
Let $V(\Gamma)$ and $E(\Gamma)$ denote the set of vertices of $\Gamma$ and 
the set of edges of $\Gamma$, respectively.
A {\em half-edge\/} $c$ of $\Gamma$ is a pair $c=(v(c);e(c))$ where 
$v(c) \in V(\Gamma)$, $e(c) \in E(\Gamma)$ and $v(c)$ belongs to $e(c)$.
The set of half-edges of $\Gamma$ will be denoted by $H(\Gamma)$ and its 
two natural projections above onto $V(\Gamma)$ and $E(\Gamma)$ will be 
denoted by $v$ and $e$, respectively.
An {\em automorphism\/} of a Jacobi diagram $\Gamma$ 
is a permutation $\phi$ of $H(\Gamma)$ verifying the two following conditions
$$\begin{array}{ccc}
\big( e(c)=e(c') \big) &\Rightarrow& \big( e\big( \phi(c)\big) =
e\big( \phi(c')\big) \big) \\
\big( v(c)=v(c')\big) &\Rightarrow& \big( v\big( \phi(c)\big) =
v\big( \phi(c')\big) \big) 
\end{array}$$ 
for any 
$c,c'\in H(\Gamma)$.
An automorphism $\phi$ of a Jacobi diagram $\Gamma$ 
{\em preserves the vertices of $\Gamma$ \/} if
$$\forall c\in H(\Gamma), v\big( \phi(c)\big) =v(c).$$
Let $\mbox{\rm Aut}(\Gamma)$ be the set of automorphisms of $\Gamma$.
Let  $\mbox{\rm Aut}_V(\Gamma)$ denote the set of automorphims of $\Gamma$
that preserve the vertices of $\Gamma$. Let $\sharp \mbox{\rm Aut}_V(\Gamma)$
denote the number of automorphisms
of $\Gamma$ that preserve the vertices.  
A {\em vertex-orientation\/} of a Jacobi diagram $\Gamma$ is an {\em orientation\/} 
of each vertex of $\Gamma$, that is a cyclic order of the three
half-edges that meet at that vertex. Two vertex-orientations of $\Gamma$ are 
equivalent if 
and only if the cardinality of the set of vertices where they differ is even.
An {\em orientation\/} of $\Gamma$ is an equivalence class of vertex-orientations.
An {\em oriented\/} Jacobi diagram is a Jacobi diagram carrying an 
orientation. 
A Jacobi diagram $\Gamma$ is {\em reversible\/} if there exists an automorphism $\phi$ 
of $\Gamma$ that
reverses an orientation of $\Gamma$.
For any automorphism $\phi$ of $\Gamma$, set
$$\mbox{\rm sign}(\phi )=
\left\{
\begin{array}{cl}
1&{\rm if\ } \phi {\rm\  preserves\ the\ orientation}\\
-1&{\rm if\ } \phi {\rm\  reverses\ the\ orientation.}
\end{array}\right.$$
Note that, for all $\phi\in \mbox{\rm Aut}_V(\Gamma)$, 
$\mbox{\rm sign}(\phi)=1$.\\ 
The {\em degree\/} of a Jacobi diagram is 
half the number of all its vertices. 
Let $\CA_k $ denote the free abelian group generated by the degree $k$ 
oriented Jacobi diagrams, quotiented out by the following relations AS and IHX.
\begin{center}
$\begin{array}{c}
\begin{pspicture}[.2](0,-.2)(.8,1)
\psset{xunit=.7cm,yunit=.7cm}
\psarc[linewidth=.5pt](.5,.5){.2}{-70}{15}
\psarc[linewidth=.5pt](.5,.5){.2}{70}{110}
\psarc[linewidth=.5pt]{->}(.5,.5){.2}{165}{250}
\psline{*-}(.5,.5)(.5,0)
\psline{-}(.1,.9)(.5,.5)
\psline{-}(.9,.9)(.5,.5)
\end{pspicture}
+
\begin{pspicture}[.2](0,-.2)(.8,1)
\psset{xunit=.7cm,yunit=.7cm}
\pscurve{-}(.9,.9)(.3,.7)(.5,.5)
\pscurve[border=2pt]{-}(.1,.9)(.7,.7)(.5,.5)
\psline{*-}(.5,.5)(.5,0)
\end{pspicture}=0\\
\mbox{\small\rm AS-relation}
\end{array}
\hspace{40pt}
\begin{array}{c}
\begin{pspicture}[.2](0,-.2)(.8,1)
\psset{xunit=.7cm,yunit=.7cm}
\psline{-*}(.1,1)(.35,.2)
\psline{*-}(.5,.5)(.5,1)
\psline{-}(.75,0)(.5,.5)
\psline{-}(.25,0)(.5,.5)
\end{pspicture}
+
\begin{pspicture}[.2](0,-.2)(.8,1)
\psset{xunit=.7cm,yunit=.7cm}
\psline{*-}(.5,.6)(.5,1)
\psline{-}(.8,0)(.5,.6)
\psline{-}(.2,0)(.5,.6)
\pscurve[border=2pt]{-*}(.1,1)(.3,.3)(.7,.2)
\end{pspicture}
+
\begin{pspicture}[.2](0,-.2)(.8,1)
\psset{xunit=.7cm,yunit=.7cm}
\psline{*-}(.5,.35)(.5,1)
\psline{-}(.75,0)(.5,.35)
\psline{-}(.25,0)(.5,.35)
\pscurve[border=2pt]{-*}(.1,1)(.2,.75)(.7,.75)(.5,.85)
\end{pspicture}
=0 \\
\mbox{\small\rm IHX-relation}
\end{array}$
\end{center}
Each of these relations relate diagrams which can be represented by 
immersions that are identical outside the part of them represented in the 
pictures. In the pictures, the cyclic order of the half-edges is represented
by the counterclockwise order. For example, 
AS identifies the sum of two diagrams which only differ by the orientation
at one vertex to zero. 
The space $\CA_0$ is equal to $\ZZ$ generated by the empty diagram.
In what follows, if $\Gamma$ is an oriented Jacobi diagram, then $-\Gamma$ denotes
the same Jacobi diagram with the opposite orientation. If $\Gamma$ is reversible, 
then $\Gamma=-\Gamma$.

\subsection{$Y$\!--graphs and the Goussarov--Habiro filtration} \label{S23}
 
Here, we briefly review the {\em $Y$\!--surgery\/}, or the 
{\em surgery along $Y$\!--links\/}, which is presented in \cite{ggp}. 
The $Y$\!--surgery is equivalent to the {\em Borromeo transformation\/} 
in Matveev's work \cite{mat}.

\begin{figure}[ht!] 
\centerline{
\subfigure[\label{ssa}]{
\relabelbox \small 
\epsfysize 1.2truein\epsfbox{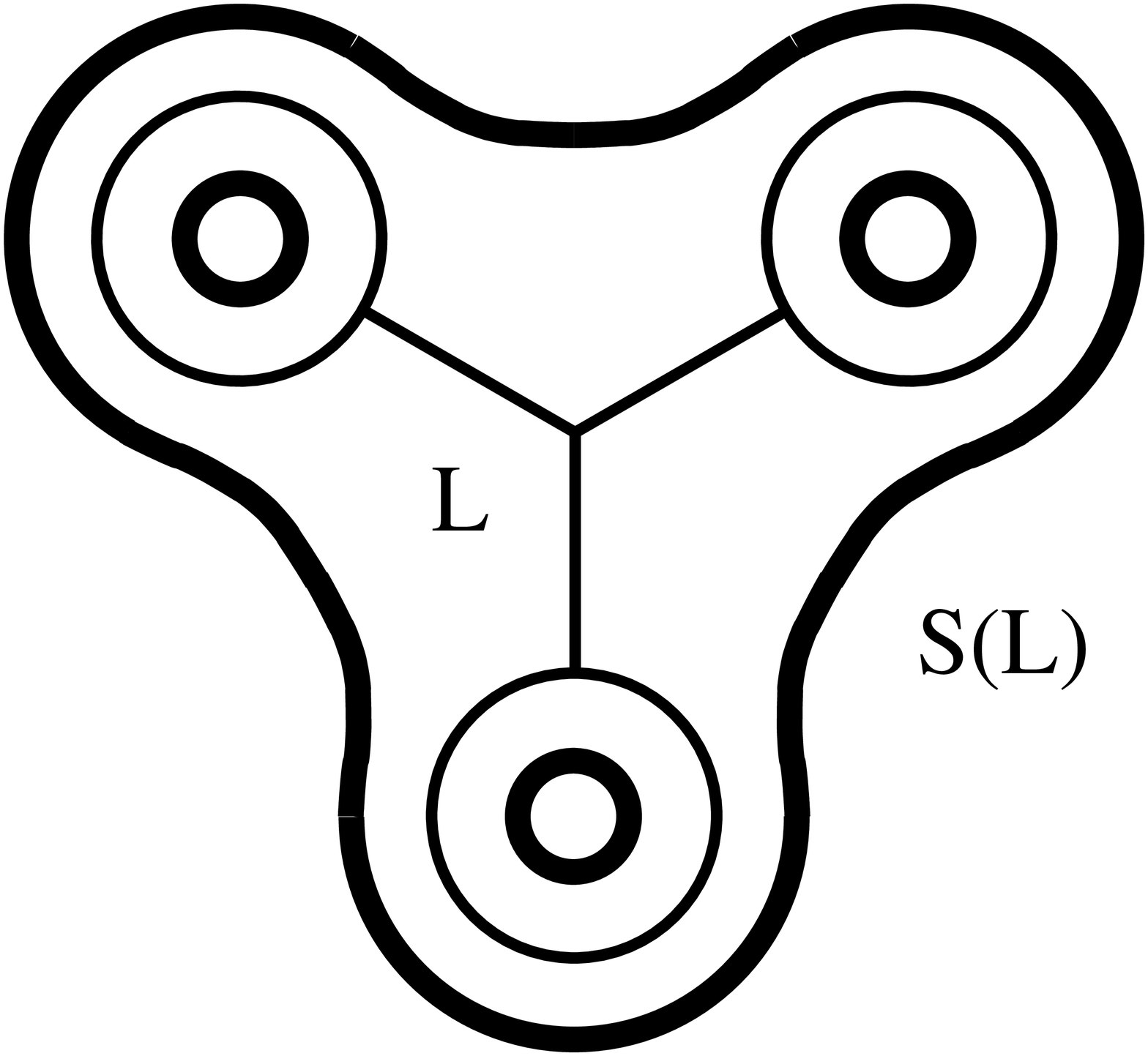}
\relabel {L}{$\Lambda$}
\relabel {S(L)}{$\Sigma (\Lambda )$}
\endrelabelbox}
\hspace{40pt}
\subfigure[\label{ssb}]{
\relabelbox \small
\epsfysize 1.2truein\epsfbox{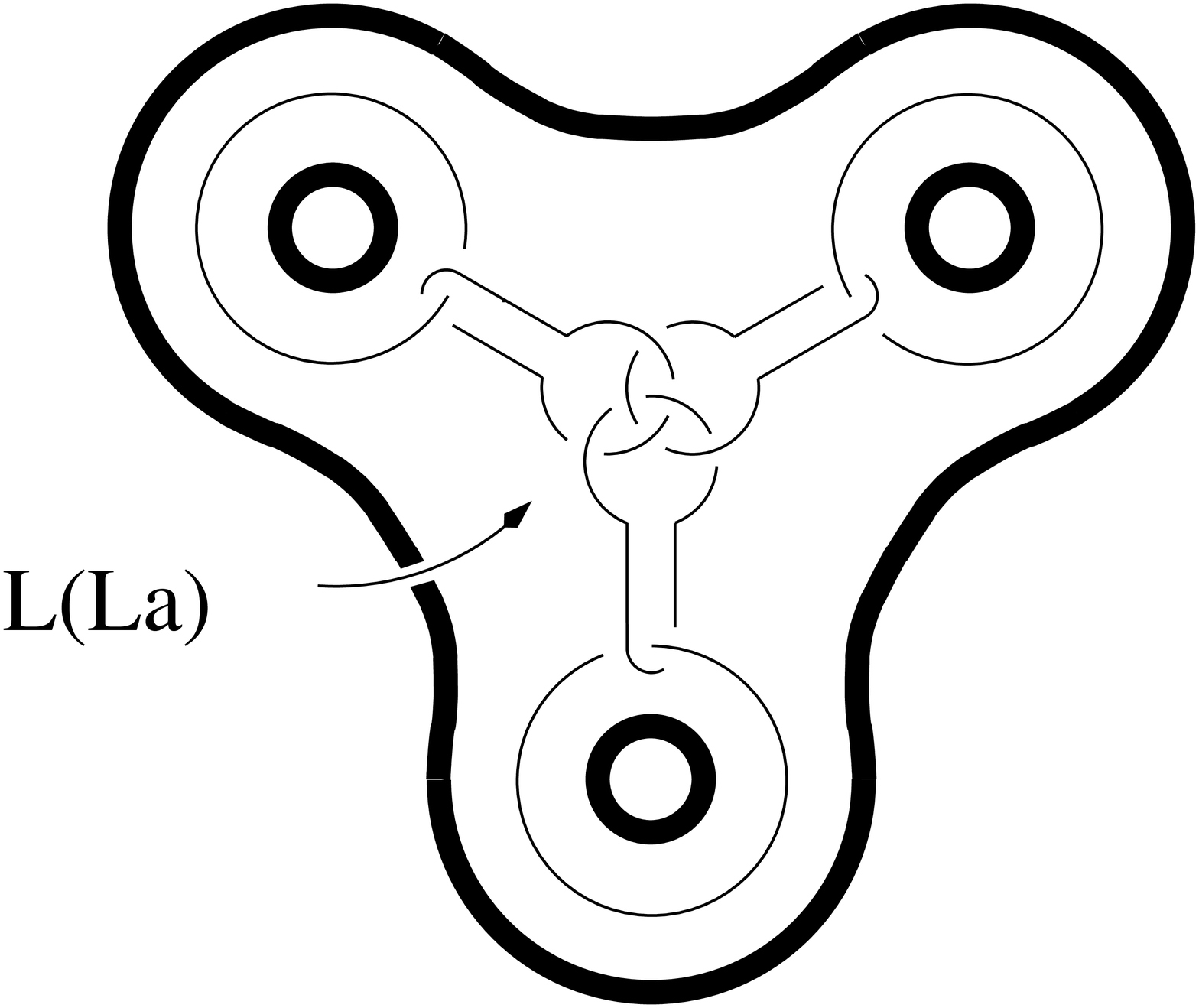}
\relabel {L(La)}{$L(\Lambda )$}
\endrelabelbox}
}\vspace{-2mm}
\caption{$Y$\!--graph and associated link}  \label{fig1}
\end{figure}

Let $\Lambda$ be the graph embedded in the surface $\Sigma(\Lambda)$
shown in \mbox{Figure \ref{ssa}.} 
In the $3$--handlebody $(N=\Sigma(\Lambda) \times [-1,1])$,
the edges of $\Lambda$ are framed by a vector field normal
to $\Sigma(\Lambda)=\Sigma(\Lambda) \times\{0\}$.  $\Sigma(\Lambda)$ is called a {\em framing surface\/} for $\Lambda$.
Let $L(\Lambda)\subset N$ be the link presented in 
\mbox{Figure \ref{ssb}}
with six framed components 
that inherit their framings from $\Sigma(\Lambda)$.
\setcounter{subfigure}{0}

Let $M$ be a $3$--manifold. A {\em $Y$\!--graph\/} in $M$ is an 
embedding $\phi$ of $N$ (or $\Sigma(\Lambda)$) into $M$ up to isotopy. 
Such an isotopy class is determined by the framed image of the framed
unoriented graph $\Lambda $ under $\phi$. 
A {\em leaf} of a $Y$\!--graph $\phi$ 
is the image under $\phi$ of a simple loop of our graph $\Lambda$.
An {\em edge} of $\phi$ 
is an edge of $\phi(\Lambda)$ that is not a leaf. 
The {\em vertex} of $\phi$
is  the unique vertex of $\phi(\Lambda)$ adjacent to the 
three edges.
With this terminology, a $Y$\!--graph  has one vertex, three 
edges and three leaves:

\begin{figure}[ht!]
\centerline{\relabelbox \small \epsfysize 3cm \epsfbox{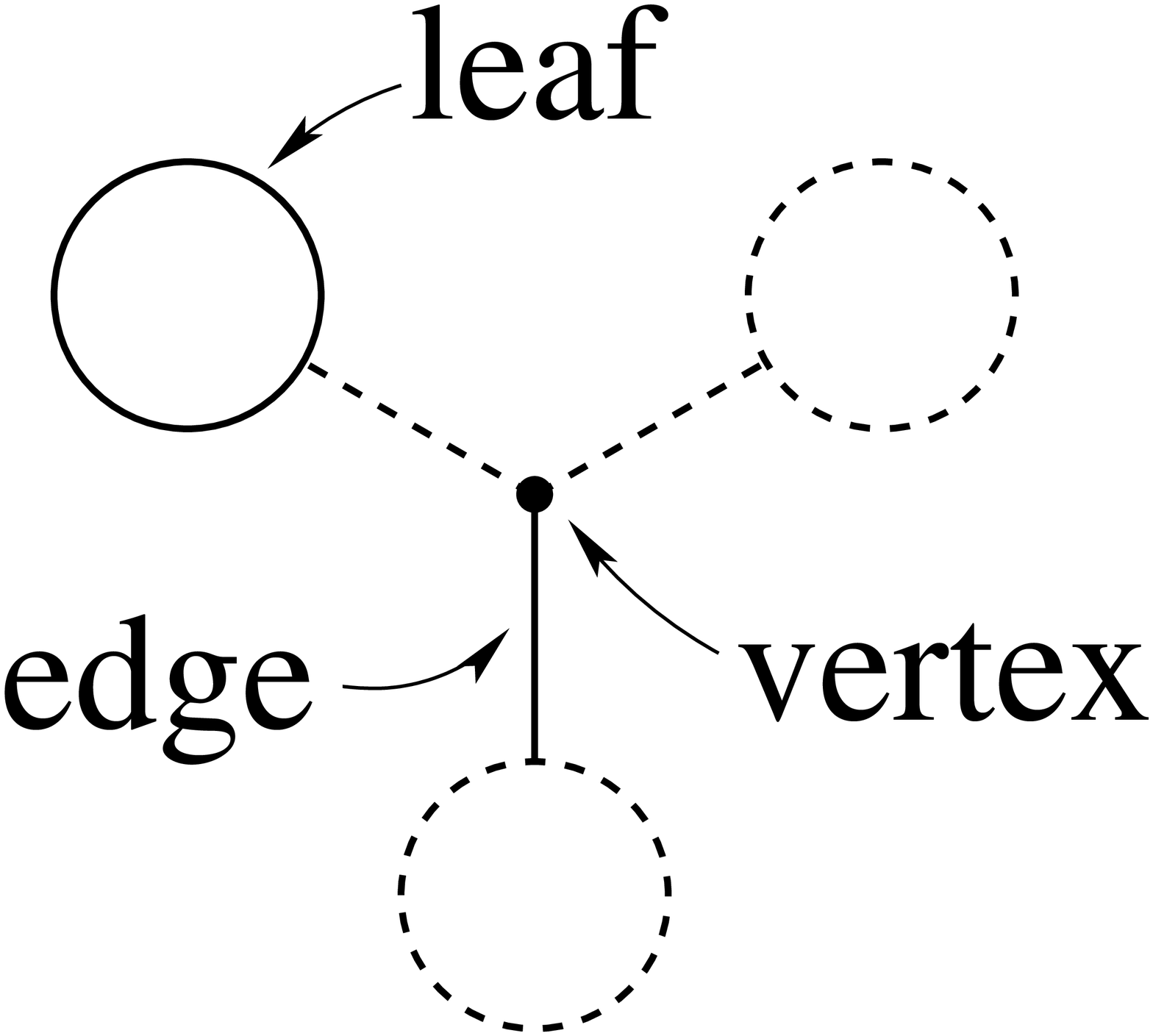}
\relabel {leaf}{leaf}
\relabel {edge}{edge}
\relabel {vertex}{vertex}
\endrelabelbox}
\end{figure}

Let $G\subset M$ be a \mbox{$Y$\!--graph.} A leaf $l$ of a $Y$\!--component of $G$
is {\em trivial} if $l$  bounds an 
embedded disc that induces the framing of $l$,
in $M\setminus G$.

The {\em $Y$\!--surgery along the $Y$\!--graph\/} $\phi(\Lambda)$ is 
the surgery along the framed link $\phi(L(\Lambda))$ (see 
\cite[Chapter 9]{Rol}, 
or \cite[Chapter 11]{Lic}
for details about surgery on framed knots).
The resulting manifold is denoted by $M_{\phi(\Lambda)}$. 
An {\em $n$--component $Y$\!--link\/}
$G\subset M$ is an embedding of the disjoint  union of $n$ copies of $N$ 
into $M$ up to isotopy. The  
$Y$\!--surgery along a $Y$\!--link $G$ is defined as the surgery
along each  $Y$\!--component of $G$. The resulting manifold is 
denoted by $M_G$.

In this article, 
the homology coefficients will always be integers.
A {\em $\ZZ$--sphere\/} is a  compact oriented $3$--manifold $M$ such that 
$H_{\ast}(M)=H_{\ast}(S^3)$. It is also called a 
{\em homology sphere}. 
A {\em homology handlebody}  or {\em $\ZZ$--handlebody\/} is an oriented, 
compact $3$--manifold $A$ with
the same homology (with integral coefficients) as the standard (solid) 
handlebody $H_g$ below.

\begin{figure}[ht!]
\centerline{\relabelbox \small
\epsfysize 1.5cm \epsfbox{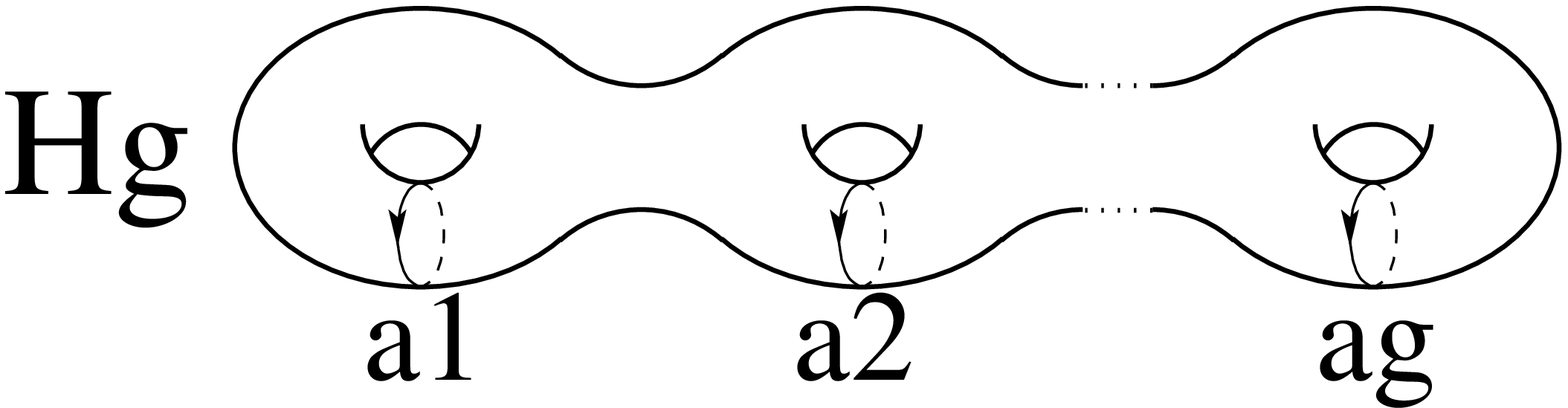}
\relabel {Hg}{$H_g$}
\relabel {a1}{$a_1$}
\relabel {a2}{$a_2$}
\relabel {ag}{$a_g$}
\endrelabelbox}
\end{figure}

Note that the boundary $\partial A$ of such a $\ZZ$--handlebody $A$ is homeomorphic to the
boundary $\Sigma_g$ of $H_g$. For any surface ${\Sigma}$,
let $\langle,\rangle_{\Sigma}$ be the intersection form on $H_1(\Sigma)$.
For a $\ZZ$--handlebody $A$,  
${\cal L}_A$ denotes the kernel of the map from $H_1(\partial A)$ to $H_1(A)$ induced by the inclusion.
It is a Lagrangian of $(H_1(\partial A);\langle,\rangle_{\partial A})$. It is called 
the {\em Lagrangian}\/ of $A$.

If $A$ is a $\ZZ$--handlebody  and if $G$ is a $Y$\!--link in the interior $\IT(A)$ of $A$,
then $A_G$ is still a $\ZZ$--handlebody whose boundary $\partial A$ is canonically
identified with $\partial A_G$, so that $\CL_A=\CL_{A_G}$. 
Similarly, if $G$ is a $Y$\!--link in 
a homology sphere $M$, then $M_G$ is still a homology sphere.

Let $\mathcal{F}$ be the abelian  group freely generated by the 
oriented $\ZZ$--spheres up to orientation-preserving diffeomorphisms.
Let $M$ be a $\ZZ$--sphere and let $G\subset M$ be a $Y$\!--link with $n$ 
components indexed by $\{1,\dots, n\}$. For any subset $J\subset \{1,\dots, n\}$, 
let $G(J)$ be the $Y$\!--sublink 
of $G$ made of the components of $G$ whose indices are in $J$.
Set
$$
[M,G]=\sum_{J\subset \{1,\dots , n\}} (-1)^{\sharp J} M_{G(J)}\in \mathcal{F}.
$$
Let $\mathcal{F}_n$ denote the subgroup of 
$\mathcal{F}$ generated
by all the elements $[M,G]$, where $G$ is an $n$--component $Y$\!--link in a $\ZZ$--sphere $M$.
This defines a filtration 
$$\mathcal{F}_0=\mathcal{F}\supset
 \mathcal{F}_1\supset \dots \mathcal{F}_n\supset \dots$$ 
of $\mathcal{F}$. It is called the Goussarov--Habiro filtration 
(see \cite{ggp} and \cite{hab}). 
Set $$\mathcal{G}_n=\mathcal{F}_n / \mathcal{F}_{n+1}.$$

\subsection{Linking Jacobi diagrams to the Goussavov--Habiro filtration} 
\label{SII3}

\begin{figure}[ht!] 
\centerline{
\subfigure[\label{ssf1}]{\relabelbox \small
\epsfysize 2cm \epsfbox{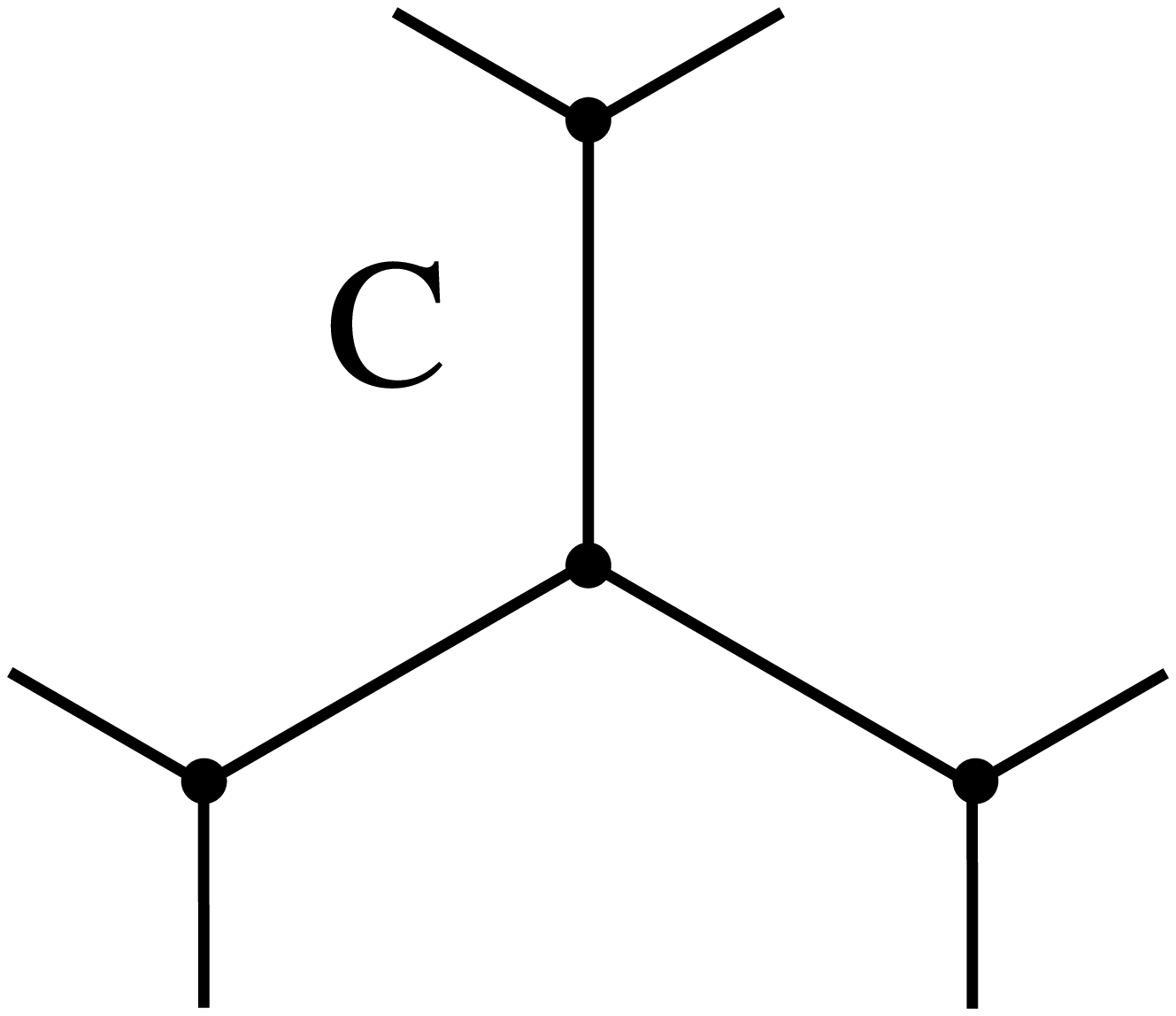}
\relabel {C}{$\Gamma$}
\endrelabelbox
\hspace{15pt} \raisebox{2cm}[0pt][0pt]{\raisebox{-1cm}{$\longrightarrow$}} 
\hspace{15pt}
\includegraphics[scale=0.2, height=2cm]{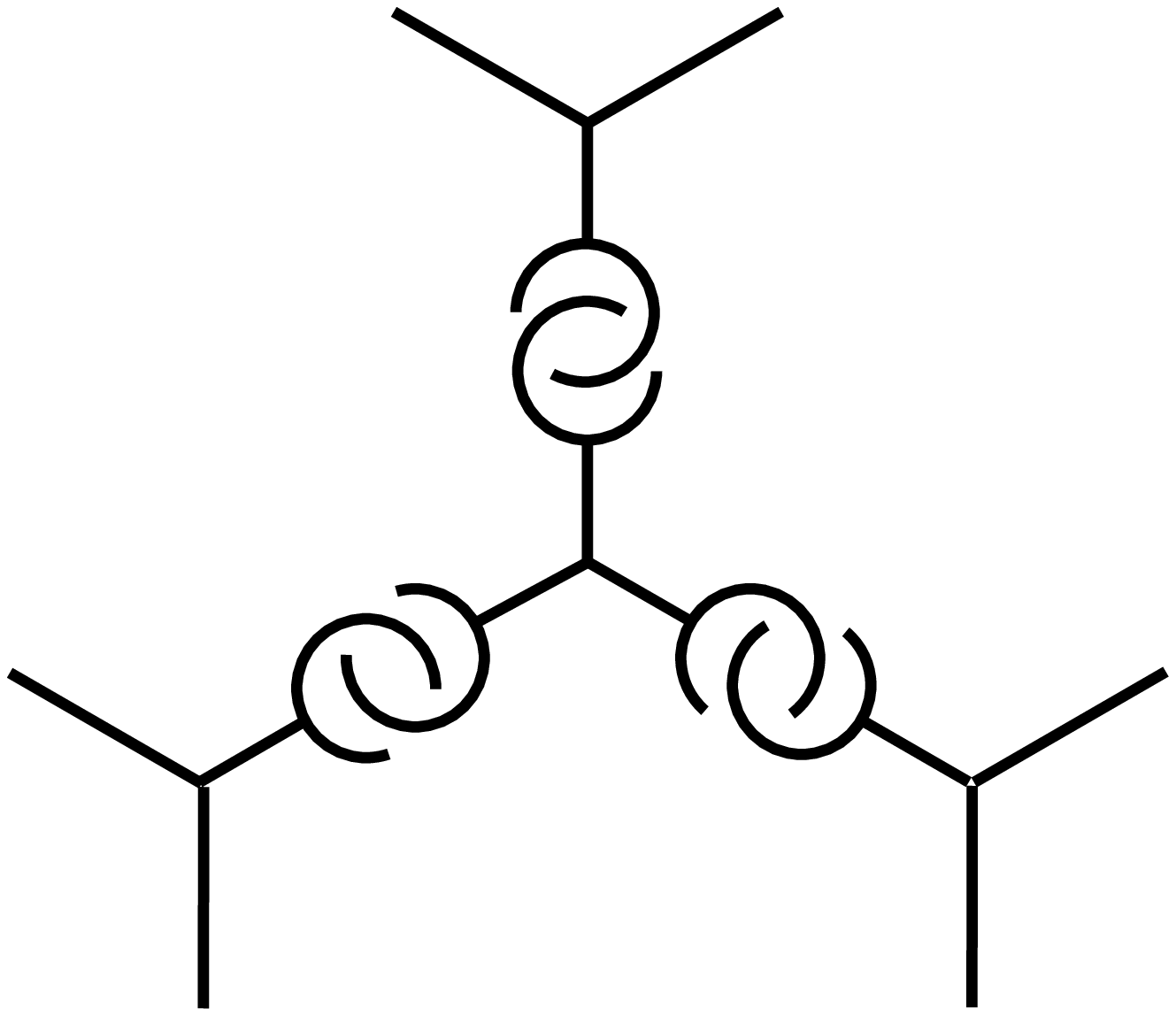}}
\hspace{40pt}
\subfigure[\label{ssf3}]{\relabelbox \small
\epsfysize 2cm \epsfbox{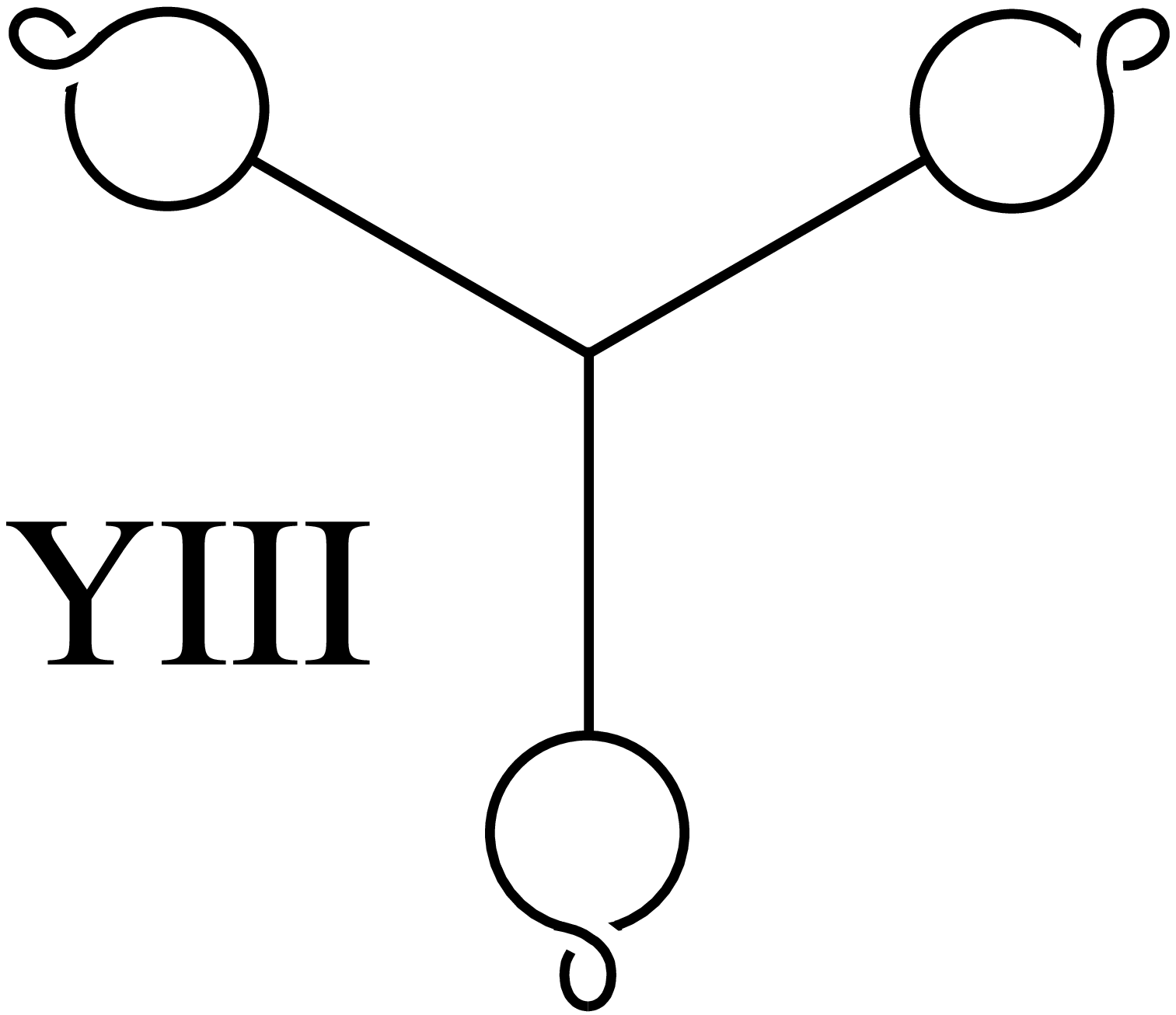}
\relabel {YIII}{$Y_{I\!I\!I}$}
\endrelabelbox
}}
\vspace{-2mm}
\caption{Turning a Jacobi diagram into a $Y$\!--link} \label{abs}
\end{figure}

Below, following \cite{ggp}, we describe a surjective map from
$\oplus_{2k\leq n} \mathcal{A}_k$  to $\mathcal{G}_{n}$, whose tensor product
by $\ZZ[1/2]$ is an isomorphism. 
Let $k$ and $n$ be integers such that $2k\leq n$. Let $\Gamma $ be a 
 degree $k$ oriented
Jacobi diagram. Let $\tilde{\Gamma}$ be an arbitrary framed
embedding of $\Gamma $ in 
$S^3$, where the framing is induced by a regular projection of 
$\tilde{\Gamma}$ in $\RR^2$ that induces the counterclockwise
orientation of the trivalent vertices of $\Gamma$. 
Insert a Hopk link on each edge of $\tilde{\Gamma}$ as illustrated 
in \mbox{Figure \ref{ssf1}.} Let $G(\tilde{\Gamma})$ denote the resulting 
$Y$\!--link in $S^3$. 

Let $Y_{I\!I\!I}$ be the framed $Y$\!--graph embedded in $S^3$ 
shown in \mbox{Figure \ref{ssf3}.}
Let $\phi_n(\Gamma)$ be the disjoint union of $G(\tilde{\Gamma})$ and of  
$n-2k$ 
copies  of $Y_{I\!I\!I}$.

\begin{thm}{\rm\cite[Theorem 4.13]{ggp}}\label{base}\qua
The  linear map 
$$\begin{array}{cccc}
\Psi_n\co&\bigoplus_{2k\leq n}\mathcal{A}_k & \longrightarrow & \mathcal{G}_n\\
        & \Gamma & \longmapsto & \overline{[S^3,\phi_n(\Gamma)]}
\end{array}$$ 
does not depend on the choice of $\phi_n$, is well-defined and is surjective. 
Moreover, in $\mathcal{G}_n$,
$$2\Psi_n(\bigoplus_{2k<n} \mathcal{A}_k)=0.$$   
 \end{thm}

That $\Psi_n$ is independent of the choice of $\phi_n$, factors through $AS$
and satisfies $2\Psi_n(\oplus_{2k<n} \mathcal{A}_k)=0$ is a consequence of Proposition~\ref{Review} proved below.
In this article, the class $\overline{[M,G]} \in \CG_n$ of the bracket $[M,G]$ of any $n$--component $Y$\!--link, will be expressed as an explicit combination of the
$\Psi_n(\Gamma)$ for oriented Jacobi diagrams $\Gamma$ with at most $n$ vertices.
Therefore, the surjectivity of $\Psi_n$ will be reproved.
For the sake of completeness, a proof that $\Psi_n$ factors through IHX is given in the appendix.

\section{Statement of the main result} \label{SIII}
\setcounter{equation}{0}

\subsection{{\CL}P--surgeries} \label{S22}

An {\em $n$--component {\CL}P--surgery\/} is a 3--tuple  
$$D=\left( M;n;(A_i,B_i)_{i=1,\dots,n}\right) $$
where
\begin{itemize}
\item $M$ is a homology sphere, $n \in \NN$,
\item for any $i=1,2, \dots n$, $A_i$ and $B_i$ are $\ZZ$--handlebodies
whose boundaries are identified by implicit diffeomorphisms 
(we shall write $\partial B_i=\partial A_i$), so that $\CL_{B_i}=\CL_{A_i}$,
\item the disjoint union of the $A_i$ is embedded in $M$. We shall write
$$\sqcup_{i=1}^n A_i \subset M.$$
\end{itemize}
For such an {\CL}P--surgery $D$, and for any subset $J\subset \{1, \dots, n\}$,
set $$D(J)=\big( M;\sharp J; (A_i,B_i)_{i\in J}\big).$$
Let $M_{D(J)}$ denote the homology sphere obtained by replacing $A_i$ by $B_i$ 
for any element $i$ of $J$.
$$M_{D(J)}=\big( M \setminus \sqcup_{i \in J} \mbox{\rm Int}(A_i) \big) 
\bigcup_{ \sqcup_{i \in J} \partial A_i} \big( \sqcup_{i \in J} B_i \big). $$
Define 
$$[D]=\sum_{J\subset\{1,\dots,n\}}(-1)^{\sharp J} M_{D(J)} \ \in \mathcal{F}.$$

The following proposition will be proved in \mbox{Subsection \ref{NI}}.

\begin{proposition} \label{T5}
For any $n$--component {\CL}P--surgery $D$,
$$[D]\in \mathcal{F}_n.$$
\end{proposition}

Conversely, any $n$--component $Y$\!--link $G=(G_i)_{i\in \{1,\dots, n\}}$ in a homology sphere $M$, 
induces the $n$--component 
{\CL}P--surgery $$(M;G)=\big( M; n; (A_i,B_i)_{i=1,\dots, n}\big)$$
such that, for any $i$ in $\{1,\dots, n\}$, 
$A_i$ is a regular neighbourhood of
the $Y$\!--component $G_i$ of $G$, and $B_i=(A_i)_{G_i}$. Then
$[(M;G)]=[M,G]$ and $(M;G)$ is called the {\em {\CL}P--surgery induced by $G$\/}.

This allows us to give the following alternative definition for the 
Goussarov--Habiro filtration.

\begin{corollary} \label{nouvelledef}
$\mathcal{F}_n$ is the subspace of $\mathcal{F}$ generated by the 
elements $[D]$, where $D$ runs among the $n$--component {\CL}P--surgeries. 
\end{corollary}

In what follows, for any $n$--component {\CL}P--surgery $D$, $\overline{[D]}$ 
denotes the class of $[D]$ in $\mathcal{G}_n$. It is called the {\em bracket\/} 
of $D$.

\subsection{The linking number of an {\CL}P--surgery with respect to a Jacobi diagram} \label{subsdeflinkLP}

This subsection is devoted to the definition of the linking number $\ell(D;\Gamma)$ of an $n$--component {\CL}P--surgery $D$ with respect to a degree $k$ Jacobi diagram $\Gamma$, with $2k \leq n$.

Let $\Gamma$ be an oriented degree $k$ Jacobi diagram.
Define a map 
$$h\co H(\Gamma)\longrightarrow \{1,2,3\}$$  such that,
for any vertex $w$ of $\Gamma$, the map
$$h_w \co v^{-1}(w)\longrightarrow \{1,2,3\}$$
is a bijection. Set 
$$\mbox{\rm sign}(h)=\prod_{w\in V(\Gamma)}\mbox{\rm sgn}(h_w)$$
where, for any vertex $w$ of\/ $\Gamma$, 
$\mbox{\rm sgn}(h_w)=1$ \/if the orientation of
$w$ is induced by the order of the half-edges given by $h_w$, and\/
$\mbox{\rm sgn}(h_w)=-1$\/ otherwise.
A {\em coloration of $\Gamma$} is a bijection
$\sigma \co  V(\Gamma ) \longrightarrow  \{1,\dots, 2k\}$.
Below, $\sigma$ also denotes the induced map 
$\sigma \circ v \co H(\Gamma)\longrightarrow  \{1,\dots, 2k\}$.
Let 
$D=\big( M;2k;(A_i,B_i)\big)$ be a $2k$--component {\CL}P--surgery.
Let us define the {\em linking number $\ell(D;\Gamma;\sigma)$ of $D$ with respect to
$\Gamma$ and $\sigma$.} 

The boundary of an oriented manifold is always oriented with the outward 
normal first convention.
The Mayer--Vietoris boundary map 
$$\partial_{i,MV} \co H_2(A_i \cup_{\partial A_i} -B_i) \longrightarrow 
\CL_{A_i},$$
that maps the homology class of an oriented surface to the oriented boundary 
of its intersection with $A_i$, is an isomorphism.
This isomorphism carries the triple intersection of surfaces in the closed $3$--manifold 
$(A_i \cup_{\partial A_i} -B_i)$ on 
$\bigotimes^3 H_2(A_i \cup_{\partial A_i} -B_i)$ to a linear form 
$\CI(A_i,B_i)$ on
$\bigotimes_{j=1}^3 \CL^{(j)}_{A_i}$ which is antisymmetric with respect to 
the permutation of two factors, where $\CL^{(j)}_{A_i}$ denotes the 
$j^{\mbox{\small th}}$ copy of $\CL_{A_i}$.
Then the linear form $\CI(A_i,B_i)$ is an element of
$\bigotimes_{j=1}^3 \big(\CL^{(j)}_{A_i}\big)^{\ast}$ where 
$\big(\CL^{(j)}_{A_i}\big)^{\ast}$ denotes the dual 
$\mbox{Hom}(\CL^{(j)}_{A_i};\ZZ)$ of $\CL^{(j)}_{A_i}$.
Let $c\in H(\Gamma)$.
Define 
$$X(c)=\left(\CL^{(h(c))}_{A_{\sigma(c)}}\right)^{\ast}.$$
The linear form $\CI(A_i,B_i)$ belongs to 
$$\bigotimes_{\{c \in H(\Gamma);\  \sigma(c) =i\}}X(c).$$ 
Then define
$$T(D; \Gamma ;\sigma )=\mbox{\rm sign}(h) \bigotimes_{w\in V(\Gamma)} 
\CI (A_{\sigma (w)},B_{\sigma (w)})  
\in \bigotimes_{c \in H(\Gamma)}X(c).$$
Note that $T(D;\Gamma;\sigma)$ is independent of $h$.
 
\begin{notation}\label{notation}
Let $A$ be a $\ZZ$--handlebody.
Then $H_1(A)$ is canonically isomorphic to 
$\frac{H_1(\partial A)}{\CL_{A}}$.
Furthermore, the intersection form $\langle\ ,\,\rangle_{\partial A}$ induces the map
$$\begin{array}{llll} \langle\ \  , \ .\rangle \co & H_1(\partial A) &\longrightarrow & 
\CL_{A}^{\ast}\\
& x & \longmapsto &\langle\ .\ ,\ x\rangle\end{array}$$
that in turn induces an isomorphism  
from $\frac{H_1(\partial A)}{\CL_{A}}$ to 
$\CL_{A}^{\ast}$. Then 
$$\varphi_{A} \co H_1(A)\longrightarrow \CL_{A}^{\ast}$$
will denote the composition of these two isomorphisms.
\end{notation}
For $\{i,j\} \subset  \{1,2,\dots,2k\}$, the linking number in $M$ induces 
a bilinear form on $H_1(A_i) \times H_1(A_j)$ that is viewed as a linear form 
on $\CL ^{\ast}_{A_i}\otimes \CL ^{\ast}_{A_j}$ via 
$\varphi _{A_i}^{-1}\otimes \varphi _{A_j}^{-1}$.
Therefore, for each edge $f \in E(\Gamma)$ made of two half-edges $c$ and $d$ 
(such that $e^{-1}(f)=\{c,d\}$), the linking number
yields a contraction
$$ \ell_f \co  X(c) \otimes X(d) \longrightarrow \ZZ.$$
Applying all these contractions to the tensor $T(D;\Gamma ;\sigma)$  maps 
$T(D;\Gamma ;\sigma)$  to the integral {\em linking number $\ell(D;\Gamma;\sigma)$ 
of $D$ with respect to $\Gamma$ and $\sigma$.\/} 

For any automorphism $\phi$ in $\mbox{\rm Aut}(\Gamma)$, let
$$\phi_v \co V(\Gamma)\longrightarrow V(\Gamma)$$
denote the bijection such that  $v\circ \phi=\phi_v \circ v$.
Let $\mbox{\rm Bij}(\Gamma)$ denote the set of colorations of $\Gamma$. Then
$\mbox{\rm Aut}(\Gamma)$ acts on $\mbox{\rm Bij}(\Gamma)$ by the action
$$\phi \ .\ \sigma=\sigma \circ (\phi_v)^{-1}.$$ 
Let $\mbox{\rm Bij}(\Gamma)/\mbox{\rm Aut}(\Gamma)$ denote the quotient of 
$\mbox{\rm Bij}(\Gamma)$ under this action.
Note that, for any automorphism $\phi$ of $\Gamma$,
$$\ell(D;\Gamma;\sigma)={\rm sign}(\phi).\ell(D;\Gamma; \phi \ .\ \sigma).$$

The following lemma is proved at the end of the next subsection.
\begin{lemma} \label{prior}
There exists an integer $\ell _0 (D;\Gamma;\sigma)$ such that
$$\ell  (D;\Gamma;\sigma)=\sharp \mbox{\rm Aut}_V(\Gamma)\ .\ 
\ell _0 (D;\Gamma;\sigma).$$
\end{lemma}

In what follows, for any $\ZZ$--sphere $M$, $\mu(M)\in \ZZ/2\ZZ$ will 
denote the {\em Rohlin invariant} of $M$ that is the reduction mod $2$ of the Casson 
invariant (see \cite[Proposition 1.3, Definition 1.6]{gm}).

For any $n$--component {\CL}P--surgery $D=\big(M;n;(A_i,B_i)\big)$
and for any subset $J\subset \{1,\dots, n\}$,
set
$$\mathcal{L}(D(\bar{J}))=
\prod_{i\in (\{1,\dots, n\}\setminus J)} 
\big(\mu((M\setminus \IT(A_i))\cup B_i)-\mu(M)\big).$$

Let $\Gamma$ be an oriented degree $k$ Jacobi diagram. Let 
$D=\big(M;n;(A_i,B_i)\big)$ be an $n$--component {\CL}P--surgery with $2k\leq n$. 
Here, we define the {\em linking number $\ell(D;\Gamma)$ of $D$ with respect to $\Gamma$\/}.
\begin{itemize}
\item If $2k=n$ and if $\Gamma$ is not reversible, then set
$$\ell(D;\Gamma)=\sum_{\sigma\in 
\mbox{\rm Bij}(\Gamma)}
\frac{\ell(D;\Gamma;\sigma)}{\sharp \mbox{\rm Aut}(\Gamma)}
\  \in\ZZ.$$
Note that
$$\ell(D;\Gamma)=\sum_{\overline{\sigma}
\in \mbox{\rm Bij}(\Gamma)/\mbox{\rm Aut}(\Gamma)}\ell _0(D;\Gamma;\sigma).$$
\item If $2k=n$ and if $\Gamma$ is reversible, then set
$$\ell(D;\Gamma)=
\sum_{\overline{\sigma}\in \mbox{\rm Bij}(\Gamma)/\mbox{\rm Aut}(\Gamma)}
\overline{\ell_0(D;\Gamma;\sigma)}\  
\in\ZZ/2\ZZ$$
where $\overline{\ell_0(D;\Gamma;\sigma)}\in \ZZ/2\ZZ$ denotes the  
mod $2$ reduction of
$\ell_0(D;\Gamma;\sigma)$.
\item If $2k<n$, then set
$$\ell(D;\Gamma)=\sum_{\{J\subset \{1,\dots, n\}\ ;\  \sharp J =2k\}}
\ell\big( D(J);\Gamma\big) \ .\ \mathcal{L}\big( D(\bar{J}) \big)\  
\in \ZZ/2\ZZ.$$ 
\end{itemize}

\subsection{Expression of brackets of {\CL}P--surgeries 
in terms of Jacobi diagrams} \label{S26}

Let $n \in \NN$. Let $\mathcal{J}_{n}$ be a set 
of oriented Jacobi diagrams of degree at most $n/2$ that contains one Jacobi diagram in each isomorphism class of non-oriented Jacobi diagrams of degree at most $n/2$.
The main goal of this paper is to show the following result.

\begin{thm} \label{TT}
Let $D$ be an $n$--component {\CL}P--surgery. Then 
$$\overline{[D]}=\sum_{\Gamma \in \mathcal{J}_{n}}  
\ell(D;\Gamma ).  
\Psi_{n}(\Gamma ) \in \mathcal{G}_{n}.$$
\end{thm}

\bp[Proof of Lemma \ref{prior}]
For any $i\in \{1,\dots,2k\}$, let $(a_j^i)_{j\in J_i}$ be a basis of $\CL_{A_i}$,
where $J_i=\{1,\dots, g_i\}$ and $g_i$ is the genus of $\partial A_i$.
Let $(z_j^i)_{j\in J_i}$ be the basis of $H_1(A_i)$ such that, for any $k$ and $l$ 
in $J_i$, 
$(\varphi _{A_i}(z_k^i))(a_l^i)=\delta _{kl}$.
Let $c_1$, $c_2$ and $c_3$ be the three half-edges of $\Gamma$ such that
$\sigma(c_k)=i$ and $h(c_k)=k$. Then
$$\CI(A_i,B_i)=
\sum_{(j_1,j_2,j_3) \in J_i^3} 
\CI\left( A_i,B_i\right) (a^i_{j_1}, a^i_{j_2}, a^i_{j_3})\ 
\varphi _{A_i}(z^i_{j_{1}})\otimes \varphi _{A_i}(z^i_{j_{2}}) \otimes \varphi _{A_i}(z^i_{j_{3}})$$
$$=\sum_{\stackrel{(j_1,j_2,j_3) \in J_i^3}{j_1<j_2<j_3}}
\left( 
\CI(A_i,B_i)(a^i_{j_1},a^i_{j_2},a^i_{j_3}) 
\sum_{\tau \in \mathcal{S}_3} \left( \mbox{sgn}(\tau) \bigotimes_{k=1,2,3} 
\varphi _{A_i}(z^i_{j_{\tau (k)}})\right) \right)$$
where $\varphi _{A_i}(z^i_{j_{\tau (k)}}) \in X(c_k)$, $\mathcal{S}_3$ denotes 
the set of the permutations of $\{1,2,3\}$ and 
 $\mbox{\rm sgn}(\tau)$ denotes the signature of the permutation $\tau$.

Let $\CH(\Gamma)$ denote the set of maps $h^{\prime} \co H(\Gamma)\longrightarrow 
\{1,2,3\}$ such that $h^{\prime}(v^{-1}(w))=\{1,2,3\}$ for any  $w\in V(\Gamma)$. Set 
$$J=\left\{ (j_1^1,j_2^1,j_3^1,\dots, j_1^{2k},j_2^{2k},j_3^{2k})
\in \prod_{i=1}^{2k}
(J_i)^3\ ;\forall i\in \{1,\dots,2k\}, \ j_1^i<j_2^i<j_3^i \right\}.$$ 
For any $j\in J$, set
$\mathcal{J}(j)=\prod_{i=1}^{2k} \CI (A_i, B_i)(a^i_{j_1},a^i_{j_2},a^i_{j_3})$. 
Then 
$$T(D;\Gamma;\sigma)=\sum_{j\in J} \mathcal{J}(j)
\bigg(\sum_{h^{\prime}\in \CH(\Gamma)} 
\mbox{sign}(h^{\prime}) \bigotimes_{c\in H(\Gamma)} 
\varphi _{A_{\sigma(c)}}\left(z^{\sigma (c)}_{j^{\sigma(c)}_{h^{\prime}(c)}}\right)\bigg).$$
Then
$\ell (D;\Gamma;\sigma)=\sum _{h^{\prime}\in \CH(\Gamma)} 
\ell(D;\Gamma;\sigma;h^{\prime})$
where 
$$\ell (D;\Gamma;\sigma;h^{\prime})=\mbox{sign}(h^{\prime})\ \sum_{j\in J} \mathcal{J}(j)
\bigg(
\prod_{e=(c_1,c_2)\in E(\Gamma)} \ell k\big(
z^{\sigma (c_1)}_{j^{\sigma(c_1)}_{h^{\prime}(c_1)}}, 
z^{\sigma (c_2)}_{j^{\sigma(c_2)}_{h^{\prime}(c_2)}}\big)\bigg).$$
For any automorphism $\zeta\in \mbox{Aut}_V(\Gamma)$,
$\ell(D;\Gamma;\sigma;h^{\prime}\circ \zeta)=\ell(D;\Gamma;\sigma;h^{\prime})$.
Then $\ell _0(D;\Gamma;\sigma)$ is the sum of the integers 
$\ell(D;\Gamma;\sigma;h^{\prime})$
running over all classes $\overline{h^{\prime}}$ of 
$\CH(\Gamma)/\mbox{Aut}_V(\Gamma)$.
\eop

\section{Proof of the theorem}
\setcounter{equation}{0}

\subsection{Proof of Theorem \ref{TT} 
for {\CL}P--surgeries induced by Jacobi diagrams} \label{SSJD}

Here we prove \mbox{Theorem \ref{TT}} 
when  $D=(S^3;\phi_{n}(\Gamma_Y))$, where the $Y$\!--link $\phi_n(\Gamma_Y)$
is the image of a Jacobi diagram $\Gamma_Y$ under the map $\phi_{n}$ 
of \mbox{Subsection \ref{SII3}}.
It is a direct corollary of the proposition below (and of Theorem \ref{base}).

\begin{proposition} \label{fondJac}
Let $\Gamma$ be an oriented degree $k$ Jacobi diagram.
Let $\Gamma_Y$ be an oriented degree $k^{\prime}$ Jacobi diagram.
Let\/ $n$ be an integer such that $n\geq max(2k,2k^{\prime})$. Then
$$\ell \big( \phi _{n} (\Gamma_Y) ;\Gamma \big)= \left\{
\begin{array}{ll}
1 & {\rm if}\  \Gamma_Y\cong \Gamma \\
-1& {\rm if}\  \Gamma_Y\cong -\Gamma\\
0 & {\rm if}\  \Gamma_Y\ncong \pm \Gamma  
\end{array} \right.$$
where $\Gamma \cong \Gamma^{\prime}$ iff $\Gamma$ and $\Gamma^{\prime}$ are isomorphic as oriented
Jacobi diagrams.
\end{proposition}

\begin{figure}[ht!]
\refstepcounter{figure} \label{fig100}
\begin{center}
\centerline{
\relabelbox \small
\epsfysize 3cm \epsfbox{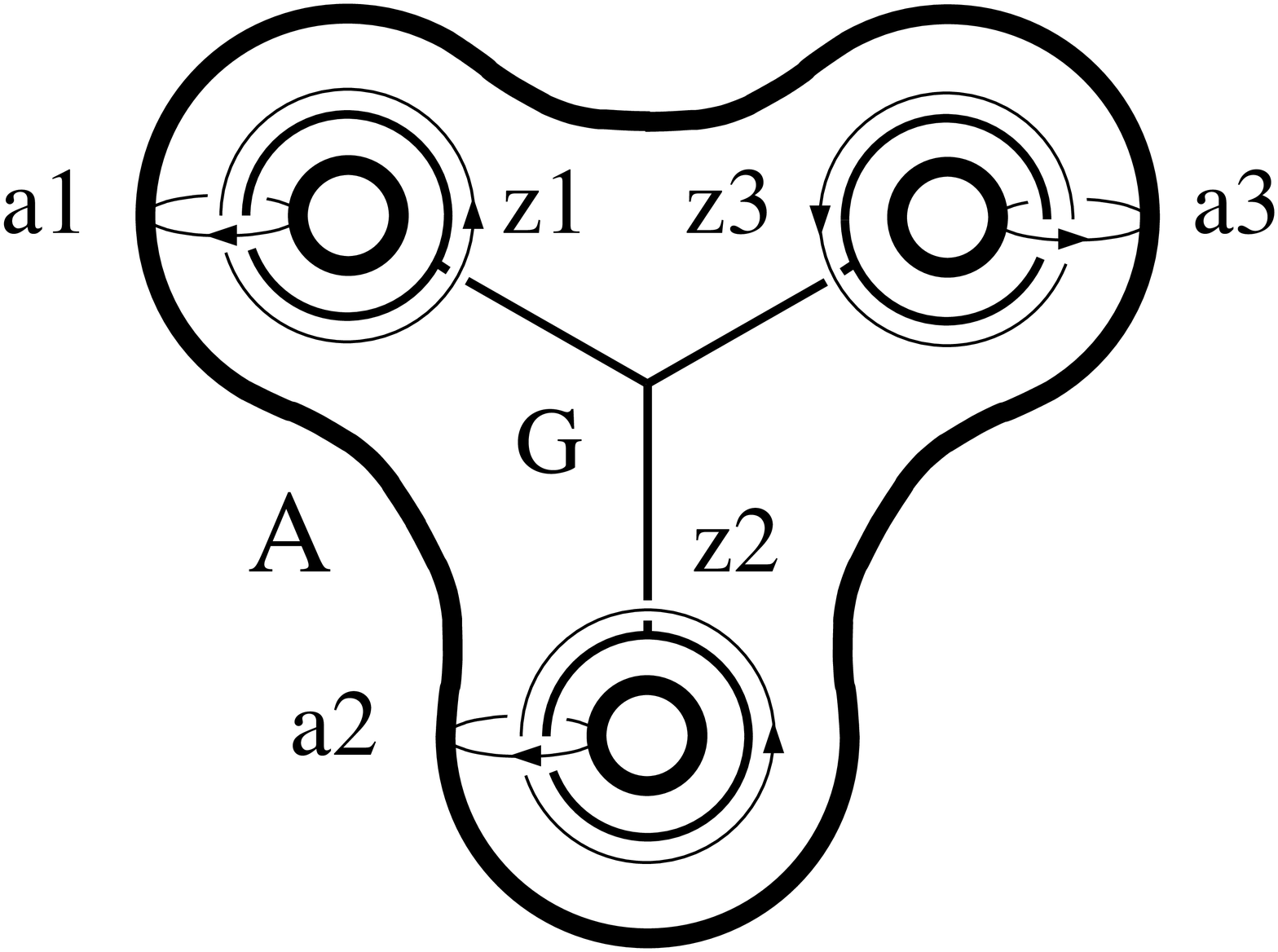}
\relabel {A}{$A$}
\relabel {a1}{$a_1$}
\relabel {a2}{$a_2$}
\relabel {a3}{$a_3$}
\relabel {z1}{$z_1$}
\relabel {z2}{$z_2$}
\relabel {z3}{$z_3$}
\relabel {G}{$G$}
\endrelabelbox}
\mbox{\small Figure \ref{fig100}} 
\end{center}
\end{figure}

\begin{lemma} \label{plop}
Let $G$ be  a framed $Y$\!--graph embedded in the interior of a $3$--handlebody 
$A$ as in \mbox{Figure \ref{fig100}.}
Let $B$ be the $\ZZ$--handlebody obtained by $Y$\!--surgery on $A$ along $G$.
Let  $(a_1,a_2,a_3)\subset \partial A$ be the oriented curves 
represented in \mbox{Figure \ref{fig100}.}
Then $(a_1,a_2,a_3)$ is a basis of $\CL_A=\CL_B$ and
$$|\big( \CI(A,B)\big) (a_1\otimes a_2 \otimes a_3)|=1.$$
\end{lemma}

\proof
This can be computed directly, or we can use that 
$$A\cup_{\partial A} (-B)=S^1\times S^1\times S^1$$ 
($A\cup(-B)$ is the manifold obtained by surgery on the $0$--framed Borromean link
in $S^3$ that is $(S^1)^3$, see \cite[13.1.5]{Thu}).
Let $S_1$, $S_2$ and $S_3$ be the three following surfaces in $(S^1)^3$. 
\begin{center}
$\begin{array}{ccc}
S_1&=&\{\star \} \times S^1\times S^1\\ 
S_2&=&S^1 \times \{\star\} \times S^1 \\
S_3&=&S^1\times S^1 \times \{\star\} .
\end{array}$
\end{center}
Let 
$\mathcal{I} \in \Big( \bigotimes^3H_2\big( (S^1)^3\big) \Big)^{\ast}$ be the
intersection form of   
$A\cup_{\partial A} (-B)=(S^1)^3$.  Since 
$S_1 \cap S_2 \cap S_3=\{\star \} \times \{\star \} \times\{\star \}$
is a single transverse intersection point, then 
$$|\mathcal{I}(S_1\otimes S_2 \otimes S_3)|=1.$$ 
By the isomorphism from $H_2\big( (S^1)^3\big) $ to
$\CL_A$ induced by the Mayer--Vietoris boundary map, $(a_1,a_2,a_3)$
can be seen as a basis of $H_2\big( (S^1)^3\big) $. Therefore,
$\mathcal{I}(S_1\otimes S_2 \otimes S_3)$ is a multiple of 
$\big( \CI(A,B)\big) (a_1\otimes a_2 \otimes a_3)$. Then
$$|\big( \CI(A,B)\big) (a_1\otimes a_2 \otimes a_3)|=1.\eqno{\qed}$$

\begin{lemma}\label{premiermu}$\qquad\mu(S^3_{Y_{I\!I\!I}})=1$
\end{lemma}

This lemma is a direct consequence of \mbox{Corollary \ref{mu}} that
is proved in Subs\-sec\-tion \ref{induced}. It relies on the results of
Subsections \ref{NI} and \ref{clovercalculus} that are logically
independent of the proof below that illustrates our formulae.

\bp[Proof of Proposition~\ref{fondJac}]$\phantom{99}$

$\bullet$\qua First assume that $2k=2k^{\prime}=n.$

Let $\sigma$ be a coloration of $\Gamma$.
Let $D=\big(S^3;\phi_n(\Gamma _Y)\big)=\big( M;n;(A_i,B_i)\big)$
be the {\CL}P--surgery induced by the $Y$\!--link $\phi_n(\Gamma _Y)$.
Each pair $(A_i,B_i)$ is a copy of the pair $(A,B)$ presented in 
\mbox{Lemma \ref{plop}}.
Let $i\in \{1,\dots ,n\}$. 
Let  
$(a^i_1, a^i_2 , a^i_3)$  be the basis of  
$\CL_{A_i}$ that corresponds to the curves $(a_1, a_2, a_3)$ in 
\mbox{Figure \ref{fig100}}.
Let $(z^i_1, z^i_2 , z^i_3)$ be the basis of $H_1(A_i)$ 
that corresponds to the curves $(z_1,z_2,z_3)$ in 
\mbox{Figure \ref{fig100}}.
Under the (implicit from now on) isomorphism
$$\varphi_{A_i} \co H_1(A_i)\longrightarrow \CL_{A_i}^{\ast}$$
presented in \mbox{Notation \ref{notation}},
$(z^i_1, z^i_2 , z^i_3)$ is the dual basis to $(a^i_1, a^i_2 , a^i_3)$, i.e.
$$\big( \varphi_{A_i}(z^i_k) \big) (a^i_l)=\delta _{kl}.$$
Then 
$$\CI(A_i,B_i)= \sum_{\tau} \mbox{\rm sgn}(\tau)\  
\CI\big(A,B\big)(a_1\otimes a_2 \otimes a_3)\ 
z^i_{\tau(1)}\otimes z^i_{\tau(2)} \otimes z^i_{\tau(3)}.
$$
Since $|\big( \CI(A,B)\big) (a_1\otimes a_2 \otimes a_3)|=1$ by 
\mbox{Lemma \ref{plop}},
$$
T(D;\Gamma;\sigma)=\mbox{\rm sign}(h)\ 
\sum_{(\tau _i)\in (\mathcal{S}_3)^{n}}
\bigg( \big(\prod_{i=1}^{n}\mbox{\rm sgn}(\tau _i)\big)\ 
\bigotimes_{c\in H(\Gamma)} z^{\sigma(c)}_{\tau_{\sigma(c)}(h(c))}
\bigg)
$$
where $h$ is as in Subsection~\ref{subsdeflinkLP}.
For any $\tau=(\tau_i)_{i=1,\dots, n}\in (\mathcal{S}_3)^{n}$, let
$\zeta(\sigma; \tau)$ denote the map
$$\begin{array}{rccl}
\zeta(\sigma;\tau) \co & H(\Gamma)&\longrightarrow & 
\{1,\dots, n\}\times \{1,2,3\}\\
& c &\longmapsto & \big( \sigma(c),\tau_{\sigma(c)}(h(c))\big).
\end{array}$$
Let 
$$\begin{array}{rccl}
\xi \co & H(\Gamma_Y) & \longrightarrow & \{1,\dots, n\}\times \{1,2,3\}\\
 & c & \longmapsto & (\xi_1(c),\xi_2(c))
\end{array}
$$
be the bijection such that, for any half-edge $c$ of $\Gamma_Y$, 
$z^{\xi_1(c)}_{\xi_2(c)}$ is the core of the leaf corresponding to $c$.
Set $\phi (\sigma;\tau)=\xi^{-1}\circ \zeta(\sigma; \tau)$. 
Then $\phi (\sigma;\tau)$ is a 
bijection from $H(\Gamma)$ to 
$H(\Gamma _Y)$ such that
\begin{itemize}
\item for any $c,c'$ in $H(\Gamma)$, 
$v(\phi (\sigma;\tau)(c))=v(\phi (\sigma;\tau)(c'))$  
if and only if $v(c)=v(c')$ 
\item for any edge $e=(c,c')$ of $H(\Gamma)$, 
$$
\ell k(z^{\sigma(c)}_{\tau_{\sigma(c)}(h(c))}, 
z^{\sigma(c')}_{\tau_{\sigma(c')}(h(c'))})\ = 
\left\{
\begin{array}{cl}
1& \text{if } 
\phi(\sigma;\tau)(c) \text{ and }
\phi(\sigma;\tau)(c')\\
&\text{belong to the same edge of }\Gamma _Y\\
0& \text{otherwise.}
\end{array}\right.$$
\end{itemize}
Therefore
$$\ell(D;\Gamma;\sigma)=\sum_{\{\tau \in (\mathcal{S}_3)^{n}; \phi(\sigma;\tau) \;\mbox{\tiny is an isomorphism} \} } \mbox{sign}(\phi(\sigma;\tau))$$
where $$\mbox{\rm sign}\left(\phi(\sigma;\tau)\right)=
\mbox{\rm sign}(h)\ (\prod_{i=1}^{n}\mbox{\rm sgn}(\tau _i)).$$
Hence, if $\Gamma \ncong \pm \Gamma _Y$, then for any coloration $\sigma$ of 
$\Gamma$, $\ell(D;\Gamma;\sigma)=0$, and $\ell(D;\Gamma)=0$.

Otherwise, there exists a coloration $\sigma$ of 
$\Gamma$ and a map $\tau \in (\mathcal{S}^3)^{n}$ such that 
$\phi(\sigma;\tau)$ is an orientation-preserving isomorphism from 
$\Gamma$ to $\mbox{\rm sign}\left(\phi(\sigma;\tau)\right) \Gamma _Y$.
For any map $\tau^{\prime}\in (\mathcal{S}^3)^{n}$ such that $\phi(\sigma;\tau^{\prime})$ 
is an isomorphism, 
$(\phi(\sigma; \tau^{\prime}))^{-1}\circ \phi(\sigma;\tau)$ is an automorphism of $\Gamma$
that preserves the vertices.

Then $\mbox{\rm sign}\big(\phi(\sigma;\tau)\big)=
\mbox{\rm sign}\big(\phi(\sigma;\tau^{\prime})\big)$. 

Conversely, any automorphism
of $\mbox{\rm Aut}_V(\Gamma)$ provides such a map $\tau^{\prime}$. Then
$$\ell(D;\Gamma;\sigma)=\mbox{\rm sign}\big(\phi(\sigma;\tau)\big) \ \sharp\mbox{\rm Aut}_V(\Gamma).$$   
For any other pair $(\sigma^{\prime};\tau^{\prime})$ such that $\phi(\sigma^{\prime};\tau^{\prime})$ is an 
isomorphism from $\Gamma$ to $\Gamma _Y$, 
$\sigma^{\prime}$ is obtained from $\sigma$ by composition by an 
automorphism of $\Gamma$. Then 
$\ell(D;\Gamma)=\ell_0(D;\Gamma;\sigma)=\mbox{\rm sign}\big(\phi(\sigma;\tau)\big)$ and Proposition~\ref{fondJac} is proved in this case.

$\bullet$\qua If $2k^{\prime}<2k=n$, then   
$\ell(\phi_{2k}(\Gamma _Y);\Gamma)=0$ because when $A$ is the regular 
neighbourhood of $Y_{I\!I\!I}$, the elements of
$H_1(A)$ do not link any element of the other $H_1(A_i)$'s.

$\bullet$\qua When $2k<n$, let $J\subset \{1,\dots, n\}$ and let
$\overline{J}=\{1,\dots,n\} \setminus J$.  Let $J_Y$ be the set of
indices of the $2k^{\prime}$--component $Y$\!--link
$G(\tilde{\Gamma}_Y)$. See Subsection~\ref{SII3}.  $\sharp
J_Y=2k^{\prime}$. Set $\overline{J_Y}=\{1,\dots,n\} \setminus J_Y$.

If $J_Y\cap \overline{J}\neq 
\emptyset$, then $\CL (D(\overline{J}))=0$
since $S^3_{Y_0}=S^3$
when $Y_0$ is a $Y$\!--graph in $S^3$ with a trivial leaf.

If $\overline{J_Y}\cap J \neq \emptyset$, then
$\ell (D(J);\Gamma)=0$ like in the previous case.
Then 
$$\ell(D;\Gamma)= \left\{ \begin{array}{ll}
0& \text{if } k\neq k^{\prime}\\
\ell(D(J_Y);\Gamma)\ .\ \CL(D(\overline{J_Y})) & \text{if }k=k^{\prime}. 
\end{array}
\right.$$
Then $\CL(D(\overline{J_Y}))=1$ by \mbox{Lemma \ref{premiermu}} and 
$\ell(D(J_Y);\Gamma)=\ell(\phi_{2k}(\Gamma _Y);\Gamma)$. 
Thus the result follows from the first case.
\eop

\subsection{Decomposition of {\CL}P--surgeries into surgeries on $Y$\!--links} 
\label{NI}

In this subsection, we recall known
facts and state useful lemmas about the theory of Borromeo surgeries \cite{mat,ggp}. We shall see
all these facts as consequences of the following single lemma~\ref{llll}. As an  application
of the theory of Borromeo surgeries, we shall prove \mbox{Proposition \ref{T5}}.

\begin{lemma}{\rm\cite[Lemma 2.1]{ggp}}\label{llll}\qua
Let $M$ be an oriented $3$--manifold (with possible boundary). Let $G$ be a $Y$\!--graph in 
$M$ with a trivial leaf 
that bounds a disc $D$ in $M\setminus G$. Then
\begin{itemize}
\item for any framed graph $T_0$ in $M\setminus G$ that does not meet
$D$, the pair $(M_G,T_0)$ is diffeomorphic to the pair $(M,T_0)$.
\item If $T$ is a framed graph 
in $M\setminus G$ that meets $\IT(D)$ at exactly one
point, then the pair $(M_G,T)$ is diffeomorphic to the pair $(M,T_G)$,
where $T_G$ is the framed graph in $M$ presented in  
\mbox{Figure \ref{bord}.} 
\end{itemize}
\refstepcounter{figure} \label{bord}
\begin{center}
\centerline{\relabelbox \small
\epsfysize 2cm \epsfbox{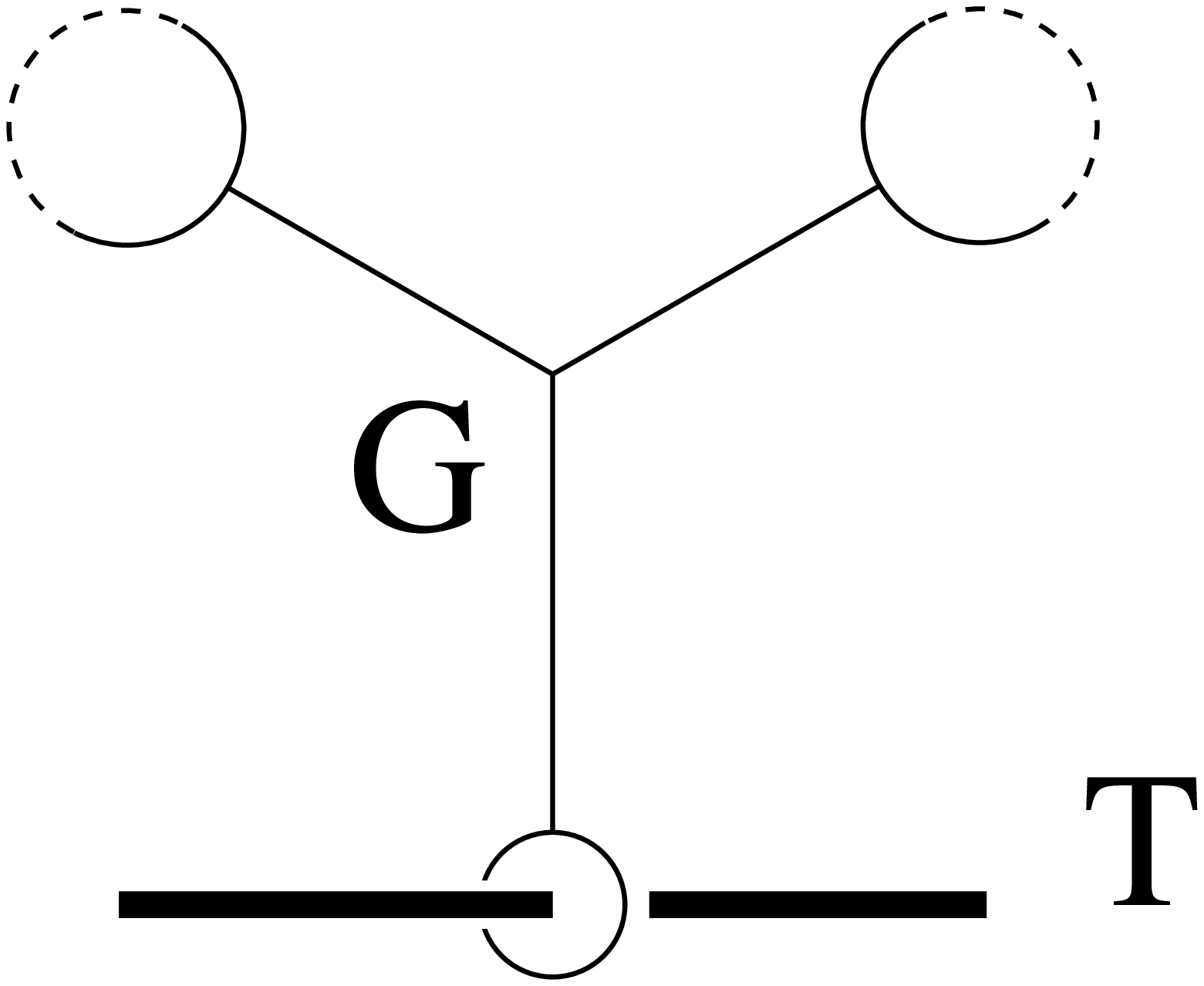}
\relabel {G}{$G$}
\relabel {T}{$T$}
\endrelabelbox
\hspace{20pt} \raisebox{1cm}[0pt][0pt]{$\longrightarrow$} 
\hspace{20pt}
\relabelbox \small
\epsfysize 2cm \epsfbox{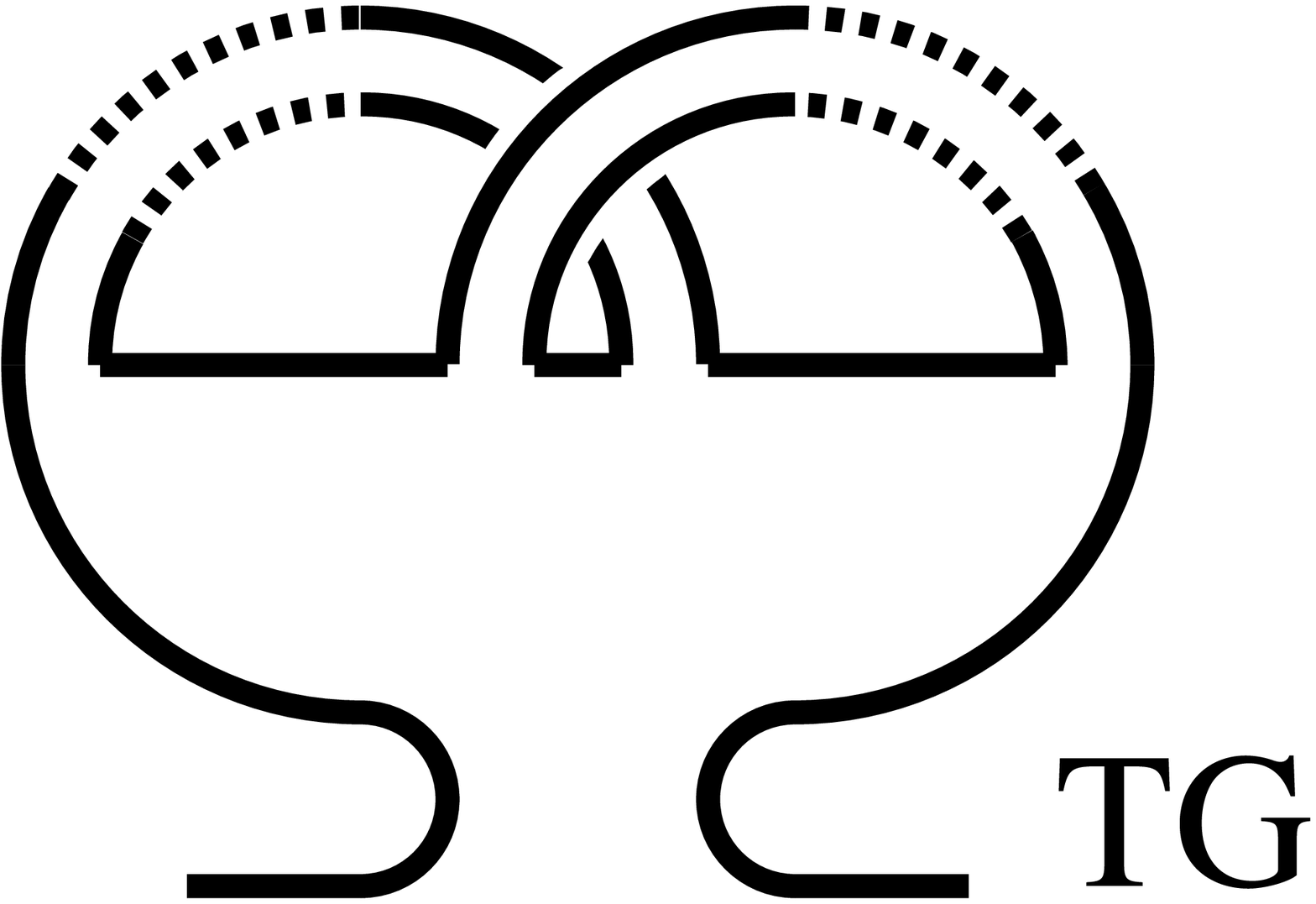}
\relabel{TG}{$T_G$}
\endrelabelbox}
\mbox{\small\rm Figure \ref{bord}}
\end{center}
\end{lemma}

\begin{corollary} \label{cinv}
Let $M$ be an oriented $3$--manifold. Let $\Sigma$ denote a genus 
$1$ surface in $M$. Let $I_1$ and $I_2$ be two intervals 
such that
\begin{itemize}
\item $\partial \Sigma=I_1\cup I_2$
\item $I_1\cap I_2=\partial I_1 =\partial I_2$
\item $I_1$ and $I_2$ are framed by a vector field normal to the surface $\Sigma$.  
\end{itemize}
Let $T$ be a framed graph such that $I_1=T\cap \Sigma$. Then there exists a
$Y$\!--graph $G$ in $M\setminus T$ with a trivial leaf
that is a meridian curve of $I_1$ such that 
the pair $(M_G,T)$ is diffeomorphic to the pair $(M,\left( T\setminus \IT(I_1)\right) \cup I_2)$
\end{corollary}

\begin{lemma}{\rm\cite[Theorem 3.2]{ggp}}\label{T2}\qua
Let $\Lambda$ be the $Y$\!--graph in the $3$--handlebody $(N=\Sigma(\Lambda) \times [-1,1])$ presented in 
\mbox{Figure \ref{ssa}}. Then there exists a $Y$\!--graph $\Lambda^{-1}$
in $N\setminus \Lambda$ such that 
the $Y$\!--surgery along $\Lambda \cup \Lambda^{-1}$ does not change $N$.
In particular,
if $M$ is a $3$--manifold then,
for any $Y$\!--graph $G$ in $M$, there exists a \mbox{$Y$\!--graph}
$G^{-1}$ in a regular neighbourhood of $G$
such that $M_{G\cup G^{-1}}=M$.
\end{lemma}

\begin{figure}[ht!] 
\centerline{\relabelbox \small
\epsfysize 2cm \epsfbox{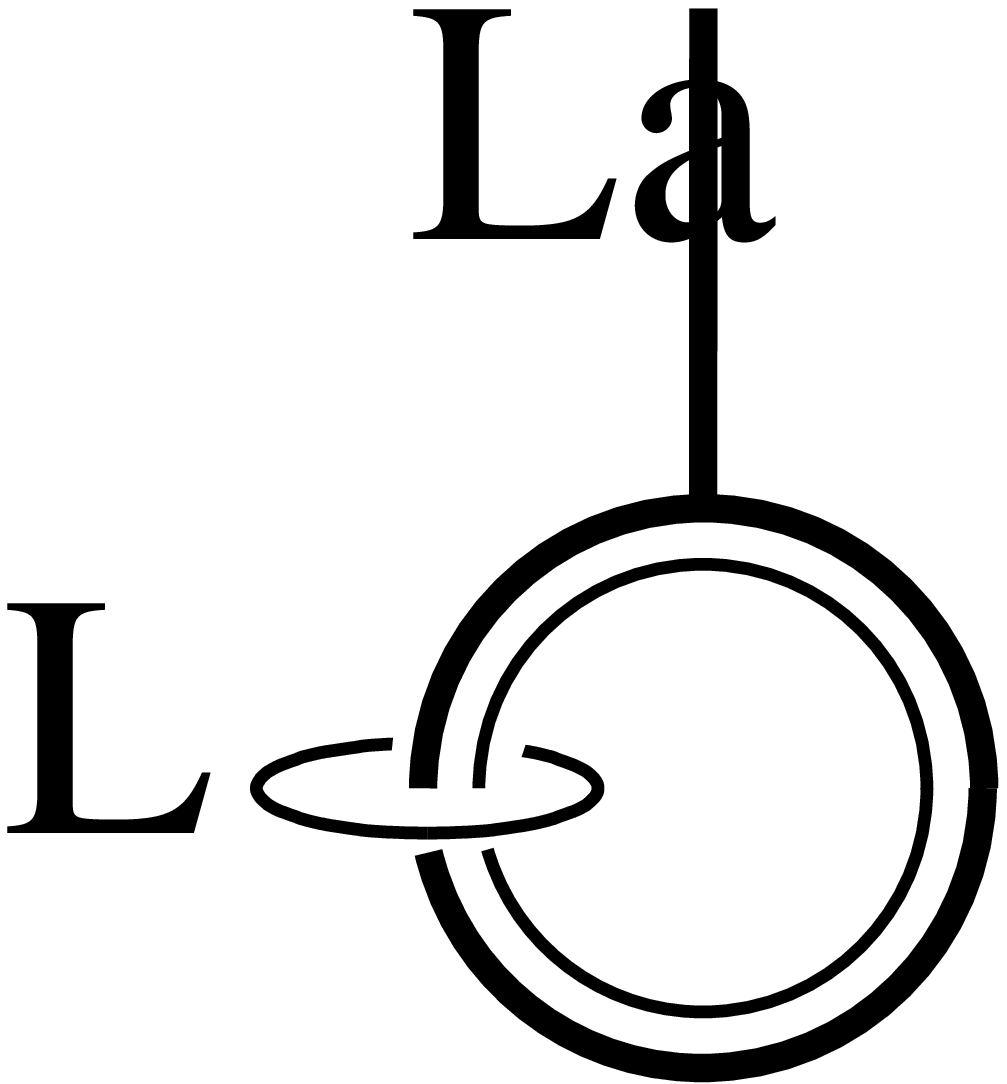}
\relabel {L}{$L$}
\relabel {La}{$\Lambda$}
\endrelabelbox
}
\caption{Trivializing a leaf} \label{gloup}
\end{figure}

\bp
Let $L$ be the framed link in $N\setminus \Lambda$  made of the two framed 
knots presented in 
\mbox{Figure \ref{gloup}},
such that $N_L=N$ and such that $\Lambda$
is isotopic to a $Y$\!--graph $\Lambda _0$ with a trivial leaf in $N_L$.
Let $L^{-1}$ denote a framed link in $N\setminus (\Lambda \cup L)$
such that the surgery along $L\cup L^{-1}$ is trivial in $N\setminus \Lambda$.
Then $L^{-1}$ corresponds to a framed link $L'$ in $N_L$.
Then $$N_{\Lambda}=N_{\Lambda \cup L \cup L^{-1}}=
(N_L)_{\Lambda_0\cup L'}.$$
The $Y$\!--surgery along $\Lambda_0$ is fully determined by 
\mbox{Lemma \ref{llll}}. It takes the tube piercing the trivial leaf and makes it 
describe its complement in the boundary of a genus one surface. 

By \mbox{Corollary \ref{cinv}}, there exists a $Y$\!--graph $\Lambda
^{-1}_0$ in $N_L\setminus (\Lambda_0\cup L')$ that undoes it.
$$N=\big( (N_L)_{\Lambda _0 \cup \Lambda ^{-1}_0} \big) _{L'}.$$
After surgery on $L'$, that does not change $N$ since the surgery on $L$ did not 
change $N$, $\Lambda_0^{-1}$ 
corresponds to a $Y$\!--graph 
$\Lambda^{-1}$ in $N\setminus \Lambda$ such that the $Y$\!--surgery
along $\Lambda \cup \Lambda^{-1}$ is trivial.\eop

\begin{remark}
\label{rksury}
What is used in the above proof and will be used again is the following principle.
Up to surgery along links, one leaf of a $Y$\!--graph can be assumed to bound a disk $D$ (pierced by surgery arcs). Then surgery along that $Y$\!--graph amounts to move the pack $T$ of framed surgery arcs piercing $D$ as indicated in Lemma~\ref{llll}, that therefore fully determines the effect of the surgery along the $Y$\!--graph.
\end{remark}

\begin{figure}[ht!]
\begin{center}
\centerline{\relabelbox \small
\epsfysize 1.5cm \epsfbox{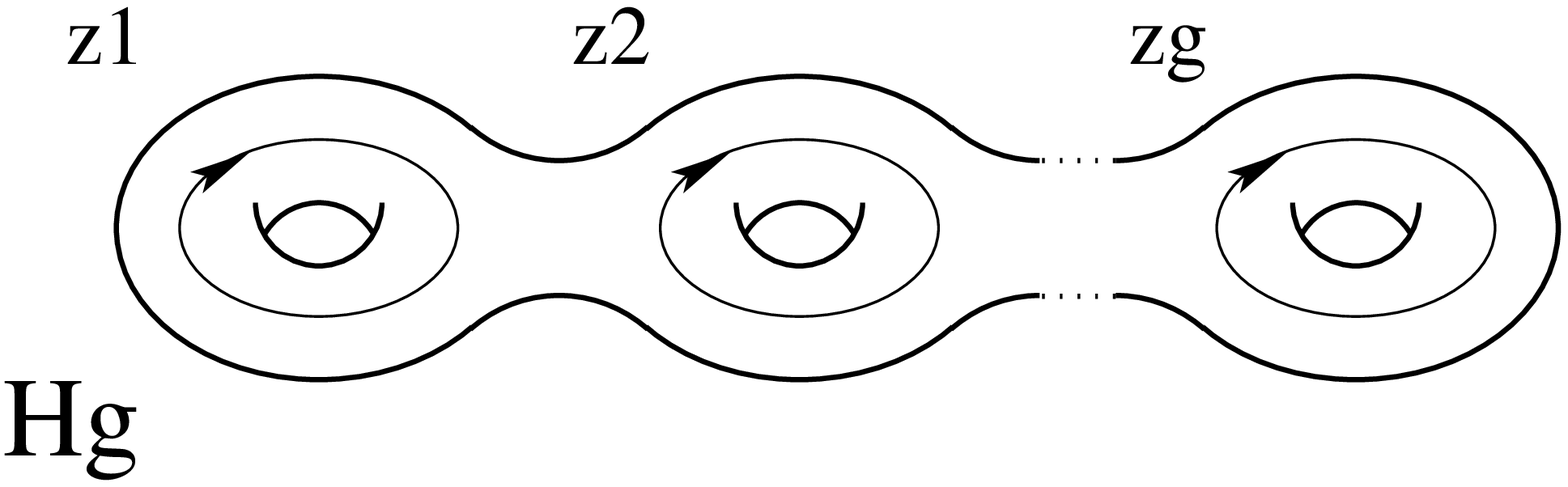}
\relabel {Hg}{$H_g$}
\relabel {z1}{$z_1$}
\relabel {z2}{$z_2$}
\relabel {zg}{$z_g$}
\endrelabelbox }
\refstepcounter{figure} \label{bou}
\mbox{\small Figure \ref{bou}}
\end{center}
\end{figure}

\begin{lemma} \label{lll}
Let $\phi$ be an embedding of the genus $g$ handlebody $H_g$ into
$S^3$. Let $z_1,\dots, z_g$ denote the curves in $\partial H_g$
presented in \mbox{Figure \ref{bou}.}  If each curve $\phi(z_i)$
bounds an embedded surface in $S^3\setminus \IT(\phi(H_g)) $, then
there exists a \mbox{$Y$\!--link} $G$ in $S^3 \setminus \phi(H_g) $
such that $S^3_G=S^3$ and the curves $\phi(z_i)$ bound embedded discs
in $S^3_G\setminus \IT(\phi(H_g))$.
\end{lemma}

\bp  
Thanks to \mbox{Corollary \ref{cinv}} and to the fact that any orientable surface
is a connected sum of genus one surfaces, 
there exists a $Y$\!--link $G_1$ in $S^3\setminus \phi(H_g)$ 
such that $S^3_{G_1}=S^3$
and $\phi(z_1)$ bounds 
an embedded disc $D_1$ in $S^3\setminus \phi(H_g)$ after 
\mbox{$Y$\!--surgery} 
on $G_1$. In particular the lemma is true for $g=1$. 
Assume that the lemma is true for handlebodies of genus \mbox{$g-1$}.
We shall use this induction
hypothesis for a regular neighbourhood $N$ of $\phi(H_g)\cup D_1$ 
equipped with the curves
$\phi(z_2),\dots, \phi(z_g)$ that are
still homologically trivial in $S^3\setminus N$.  
Let $\hat{G_1}\subset S^3_{G_1}$
be the union of the $3$--handlebodies reglued during the $Y$\!--surgery on $G_1$.
By induction hypothesis, there exists
a $Y$\!--link $G_2$ in $S^3_{G_1}\setminus N$ 
such that $(S^3_{G_1})_{G_2}=S^3$
and the curves $\phi(z_2),\dots, \phi(z_g)$ bound  embedded discs in 
$S^3\setminus N$ 
after 
\mbox{$Y$\!--surgery} on $G_2$.  After a possible isotopy in $S^3_{G_1}$, 
$G_2$ avoids $\hat{G_1}$. Then 
$G_2$ corresponds to a $Y$\!--link $G^{\prime}_2$ in $S^3\setminus 
\big( \phi(H_g)\cup G_1\big) $ 
such that $G_1\cup G^{\prime}_2$ satisfies the conclusion of the lemma. 
\eop
  
Let $M$ be a $\ZZ$--sphere.
Let $\ell k$ denote the linking number in $M$.
Let $L\subset M$ be a link, and let $L_1,\dots, L_n$
be the components of $L$. Then $L$ is {\em algebraically split}
if and only if 
$$\big( i\neq j\big) \Rightarrow \big( \ell k(L_i,L_j)=0\big) .$$  
Then \mbox{Lemma \ref{lll}} induces the following corollary.

\begin{corollary}[\mbox{\cite[Lemma 2]{mat}} or \mbox{\cite[Lemma 1.2]{mn}}]  
\label{T3}
Let $L$ be an algebraically split link in $S^3$. 
Then there exists a \mbox{$Y$\!--link} 
$G$ in $S^3\setminus L$ such that $S^3_G=S^3$ and $L$
is trivially embedded in $S^3_G$. 
\end{corollary}

\bp
Embed $H_g$ in $S^3$ so that the curves $z_i$ are the components of
$L$.
\eop

\begin{thm}{\rm\cite[Theorem 2]{mat}}\label{T1}\qua
If $M$ and $M^{\prime}$ are homology spheres, then there exists
a \mbox{$Y$\!--link} $G$ in $M$ such that $M_G=M^{\prime}$. 
\end{thm}

\bp
Since any $\ZZ$--sphere can be obtained by surgery on $S^3$ along an 
algebraically split link framed by $\pm 1$ (see \cite[Lemma 2.1]{gm}), 
and since the surgery on the trivial knot in $S^3$ framed by $\pm 1$
gives $S^3$, \mbox{Theorem \ref{T1}} is an easy corollary
of Lemma~\ref{T2} and \mbox{Corollary \ref{T3}.}
\eop

Then we can prove the following useful lemma (see \cite[Theorem 2.5]{hbg},
 too).
\begin{lemma} \label{V1}
Let $A$ and $B$ be two $\ZZ$--handlebodies with the same genus, whose boundaries
are identified so that $\CL_A=\CL_B$.
Then there exists a \mbox{$Y$\!--link} $G$ embedded in the interior of $A$ 
such that 
$A_G=B$,
where the identification of $\partial A$ with $\partial B$
is induced by the natural identification of $\partial A$ with $\partial A_G$.  
\end{lemma}

\begin{figure}[ht!]
\begin{center} 
\centerline{\relabelbox \small
\epsfysize 1.5cm \epsfbox{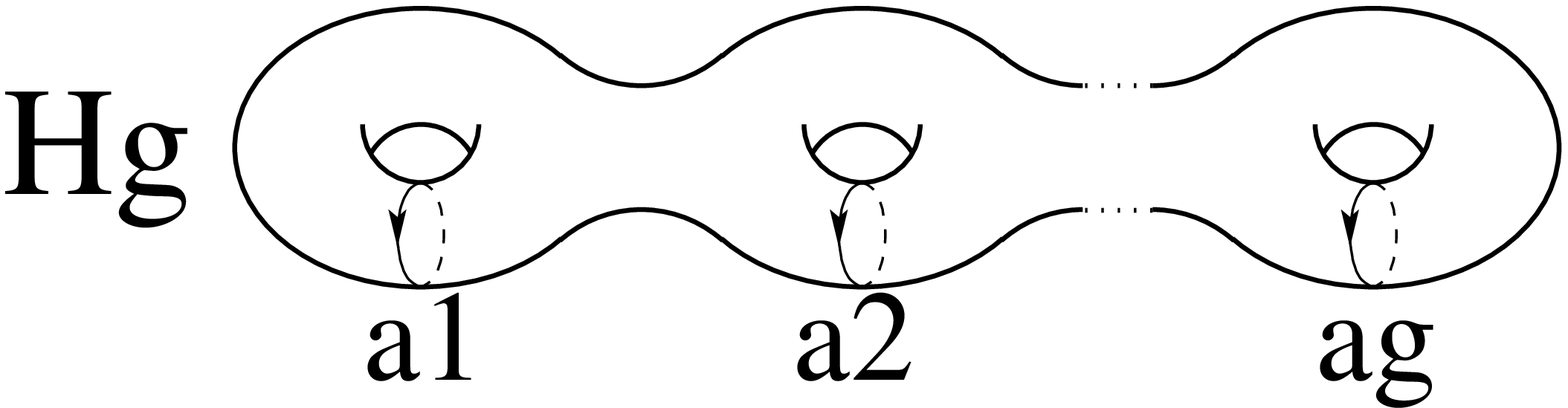}
\relabel {Hg}{$H_g$}
\relabel {a1}{$z_1$}
\relabel {a2}{$z_2$}
\relabel{ag}{$z_g$}
\endrelabelbox
\raisebox{0.8cm}[0pt][0pt]{\;\;\;\;$\subset S^3$}}
\refstepcounter{figure}  \label{fig}
\mbox{\small Figure \ref{fig}}
\end{center}
\end{figure}

\bp
Let us first prove the lemma when $B=H_g$ is the standard handlebody of genus $g$
with the boundary of $A$ identified with
the boundary $\Sigma_g$ of $H_g$ so that $\CL_A=\CL_{H_g}$. 
Embed $H_g$ trivially in $S^3$ so that  
$$\tilde{H}_g=S^3\setminus \IT(H_g)$$ 
is a standard $g$--handlebody.
Let $z_1,\dots ,z_g$ be the meridian curves of $H_g$ on 
$\Sigma_g$ presented in \mbox{Figure \ref{fig}.} 
$$M=\big( S^3\setminus {\rm Int}(H_g)\big) 
\cup_{\Sigma_g} A=\tilde{H}_g\cup
_{\Sigma _g} A.\leqno{\rm Let}$$ 
Then $M$ is a $\ZZ$--sphere.
Thus, by Theorem \ref{T1}, there exists a $Y$\!--link $G\subset M$ such that 
$M_G=S^3$. By isotopy, $G$ can avoid $\tilde{H}_g$. Then
$$S^3=A_G\cup_{\Sigma_g}\tilde{H}_g$$ 
Now $A_G$ is the complement in $S^3$
of a possibly knotted $g$--handlebody $\tilde{H}_g$.
Thanks to Lemma \ref{lll}, 
there exists a $Y$\!--link $G^{\prime}\subset \IT(A_G)$ such that 
$S^3_{G^{\prime}}= S^3$ and $\tilde{H}_g$ is embedded in $S^3_{G^{\prime}}$
so that the curves $z_i$ bound embedded discs in $A_{G\cup G^{\prime}}$. 
Thus $A_{G\cup G^{\prime}}=H_g$ with the expected boundary identification.
The general case follows easily with the help of Lemma~\ref{T2}. 
\eop

We have the following obvious lemma.

\begin{lemma} \label{ajout}
Let $D=\big( M;n;(A_i,B_i)\big) $ be an $n$--component
{\CL}P--surgery. 
Let $A^{\prime}_1$ be a $\ZZ$--handlebody such that
$\partial A^{\prime}_1$ and $\partial A_1$ are identified 
so that $\CL_{A^{\prime}_1}=\CL_{A_1}$. Let 
$$M_{A^{\prime}_1/A_1}=\big( M\setminus \IT(A_1)\big) \cup_{\partial A_1} A^{\prime}_1$$ 
denote the manifold obtained by surgery on $M$
along the pair $(A_1, A^{\prime}_1)$. Set 
$$\begin{array}{lcc}
D'&=&\big( M;n;(A_1,A^{\prime}_1),(A_2,B_2),\dots,(A_{n},B_{n})\big) \\ 
D''&=&\big( M_{A^{\prime}_1/A_1};n;(A^{\prime}_1,B_1),(A_2,B_2),\dots,(A_{n},B_{n})\big) .
\end{array}$$
Then $$[D]=[D']+[D''].$$
\end{lemma}

\bp[Proof of Proposition \ref{T5}]
Let $D=\big( M;n;(A_i,B_i)\big)$ be an $n$--component
{\CL}P--surgery. 
Thanks to \mbox{Lemma \ref{V1}}, for any $i\in \{1,\dots ,n\}$, 
there exists a $Y$\!--link 
$G^i\subset \mbox{\rm Int}(A_i)$ 
such that $(A_i)_{G^i}=B_i$. 
Let $k_i$ denote the minimal number of components for such a $G_i$.
Consider the sum $k=\sum_i k_i$.

If there exists $i\in \{ 1,\dots, n\}$ such that $k_i=0$, 
then $[D]=0\in \mathcal{F}_n$.
If, for all $i$, $k_i=1$, then $[D]\in \mathcal{F}_n$ by definition.  
Therefore $[D]\in \mathcal{F}_n$ if $k\leq n$.

If $k>n$, assume that
$k_1>1$, without loss of generality.
Then there exists a $\ZZ$--handlebody $A^{\prime}_1$ verifying the hypotheses of
\mbox{Lemma \ref{ajout}} such that $A^{\prime}_1$ can be obtained from $A_1$ by 
$Y$\!--surgery along a $Y$\!--graph in $\IT (A_1)$, and $B_1$ can be obtained
from $A^{\prime}_1$ by $Y$\!--surgery along a $Y$\!--link in $\IT (A^{\prime}_1)$ with 
$k_1-1$ components. Thus, with the notation of \mbox{Lemma \ref{ajout},}
$[D']\in \mathcal{F}_n$
and $[D'']\in \mathcal{F}_n$ by induction on $k$.
Then $[D] \in \mathcal{F}_n$ thanks to \mbox{Lemma \ref{ajout}.} 
The proposition follows. 
\eop 

\subsection{Review of the clover calculus} \label{clovercalculus}

In this section, we review the clover calculus following \cite{ggp}.
However we produce alternative proofs in the spirit of the present paper only 
based on \mbox{Lemma \ref{llll}}. Furthermore, we summarize what we shall use 
about the clover calculus in \mbox{Proposition \ref{Review}}.

\begin{figure}[ht!]
\centerline{\relabelbox \small
\epsfysize 2cm \epsfbox{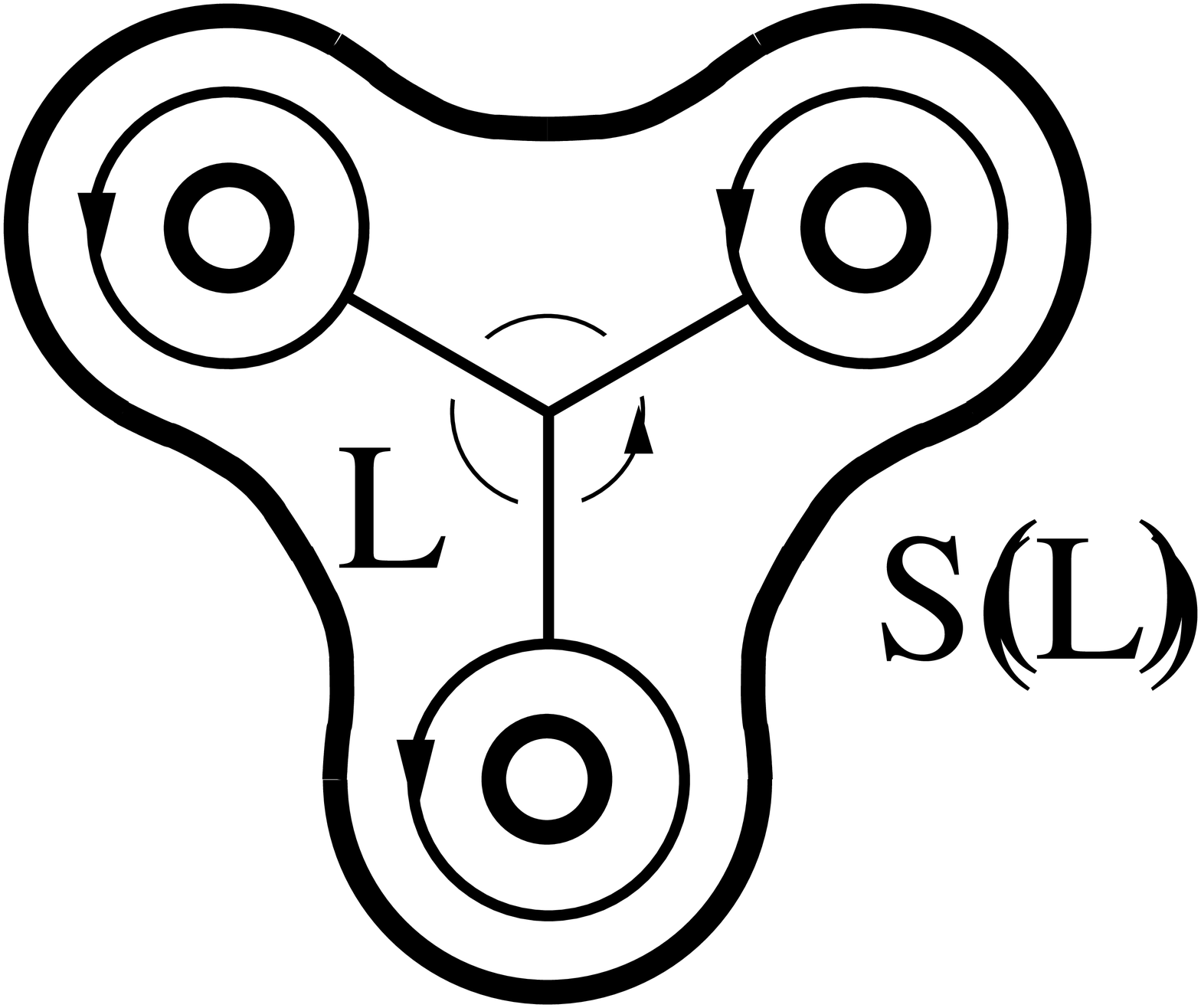}
\relabel {L}{$\Lambda$}
\relabel {S(L)}{$\Sigma (\Lambda)$}
\endrelabelbox}
\caption{oriented $Y$\!--graph} \label{frasurf}
\end{figure}

A $Y$\!--graph $\Lambda$ is {\em oriented\/} if its framing surface $\Sigma(\Lambda)$ is
equipped with an orientation.
Such an orientation provides an orientation for every leaf and (a cyclic order) for the set of leaves
of $\Lambda$. \mbox{Figure \ref{frasurf}} shows the induced orientations when 
$\Sigma (\Lambda)$ is given the standard orientation of $\RR ^2$. Reversing the
orientation of $\Sigma (\Lambda)$ reverses these four orientations.

Recall that a {\em framing\/} of a knot is a nonzero vector field normal to the knot, up to homotopy,
or a parallel
to the knot up to isotopy. In a homology sphere these two canonically equivalent notions are represented by the linking number of the knot and its parallel induced by the framing. This linking number is therefore also called the {\em framing\/} of the knot.

The goal of this subsection is to prove the following proposition.
  
\begin{proposition} \label{Review}
Let $G$ be an oriented $n$--component $Y$\!--link in a $\ZZ$--sphere $M$.
\begin{itemize}
\item[\rm(i)] The bracket $\overline{[(M;G)]}$ (in $\CG_n$) is a function independent of $M$ of
\begin{itemize}
\item  the linking numbers $\ell k(l,l^{\prime})$ where $l$ and $l^{\prime}$ are leaves in two 
distinct $Y$\!--components of $G$
\item  the products $\bar{f}(l_1)\bar{f}(l_2)\bar{f}(l_3)$ 
where $l_1$, $l_2$ and $l_3$ are leaves of a 
same $Y$\!--component, and where $\bar{f}(l)$ is the framing of $l$ in $\ZZ/2\ZZ$.
\end{itemize}
\item[\rm(ii)] Fix $G$ except for a leaf $l$ in the complement of the other parts of
$G$. Let $[l]$ denote the class of $l$ in $H_1(M\setminus \cup_{l^{\prime}\neq l}l^{\prime})$, 
where the union runs over all leaves $l^{\prime}$ distinct of $l$.
Then the bracket $(\overline{[(M;G)]} \in \CG_n)$ of $G$  is a linear map of 
$([l],\bar{f}(l))\in H_1(M\setminus \cup_{l^{\prime}\neq l}l^{\prime})\times \ZZ/2\ZZ$.
\end{itemize}
\end{proposition}

\begin{lemma} \label{firststep}
Let $G$ be an oriented $n$--component $Y$\!--link in a $\ZZ$--sphere $M$. 
The bracket $\overline{[(M;G)]}$ is a function independent of $M$ of
\begin{itemize}
\item  the linking numbers $\ell k(l,l^{\prime})$ where $l$ and $l^{\prime}$ are leaves of $G$
\item  the framings $f(l)$ where $l$ runs over the leaves of $G$.
\end{itemize}
\end{lemma}

\bp
Let $\Upsilon$ be the diagram made of $n$ copies of the diagram $\Lambda$ 
connected by an additional edge from the internal vertex of $\Lambda$ to a 
common $n$--valent vertex $p$. Embed $\Upsilon$ in $\RR^3$. Let $A$ be a regular
neighbourhood of $\Upsilon$ in $\RR^3$. Then $A$ is a union of a ball $B$ with 
$n$ copies of the genus $3$ handlebody $N$ that are glued on 
$\partial B$ along $n$ 
disjoint discs. Let $\phi _G \co A\longrightarrow M$ be an embedding of $A$ in $M$
that extends the embedding $G$.
Set $Z=M\setminus \IT(\phi _G(A))$. Then $Z$ is a genus $3n$ homology handlebody
whose Lagrangian $\CL_Z \subset H_1(\partial A)$ is fully determined by the framings
and by the linking numbers of the leaves of $G$. Therefore if
$G^{\prime}\subset M^{\prime}$ is another oriented $n$--component $Y$\!--link
with the same linking numbers and framing data, then $Z^{\prime}=M^{\prime}\setminus \IT(\phi_{G^{\prime}}(A))$
is a homology handlebody with the same lagrangian as $Z$ in $H_1(\partial A)$.

By \mbox{Lemma \ref{V1}}, there exists
a $Y$\!--link $G^{\prime \prime}\subset \IT(Z)$ such that $Z_{G^{\prime \prime}}=Z^{\prime}$. Then 
$[(M^{\prime};G^{\prime})]=[(M_{G^{\prime \prime}};G)]$. If $G^{\prime \prime}$ is a one-component $Y$\!--link, then 
$[(M;G\cup G^{\prime \prime})]=[(M;G)]-[(M_{G^{\prime \prime}};G)]$,
and $\overline{[(M;G)]}=\overline{[(M_{G^{\prime \prime}};G)]}$. By induction on the number
of components of $G^{\prime \prime}$, $\overline{[(M;G)]}=\overline{[(M_{G^{\prime \prime}};G)]}$.
Then $\overline{[(M;G)]}=\overline{[(M^{\prime};G^{\prime})]}$. \eop

\begin{figure}[ht!]
\centerline{\relabelbox \small
\epsfysize 2cm \epsfbox{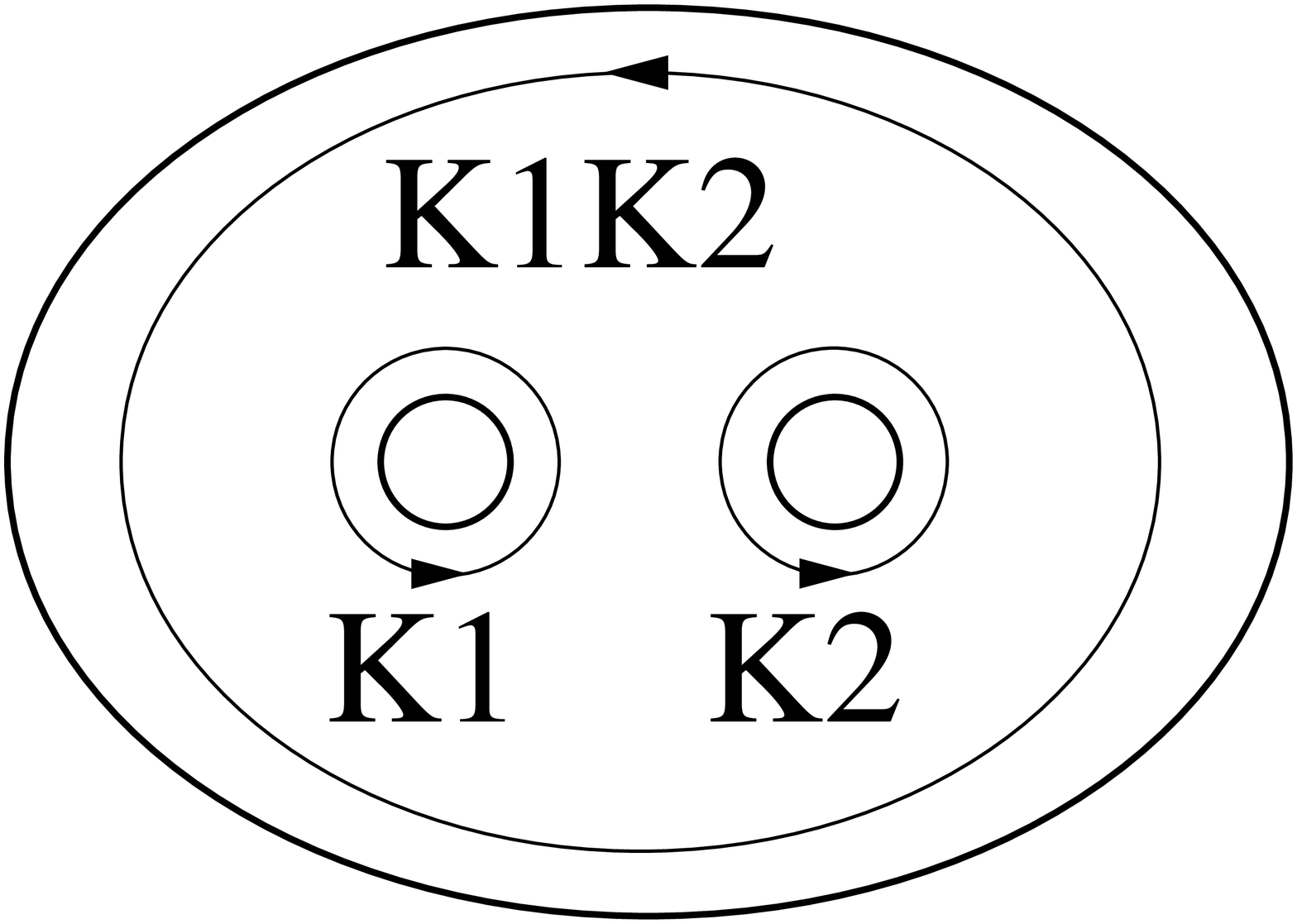}
\relabel {K1K2}{$K_1 \sharp _b K_2$}
\relabel {K1}{$K_1$}
\relabel {K2}{$K_2$}
\endrelabelbox}
\caption{Band sum of two knots} \label{2hole} 
\end{figure}

A framed knot $K_1\sharp _b K_2$ is a {\em band sum\/} of two framed oriented knots
$K_1$ and $K_2$ if there exists an embedding of a $2$--hole disk
\begin{itemize}
\item that factors the three knot embeddings by the embeddings of the three
curves pictured in \mbox{Figure \ref{2hole}} representing the disk, and
\item that induces the three framings.
\end{itemize} 
Note that
$$f(K_1\sharp _b K_2)=f(K_1)+f(K_2)+2 \ell k(K_1,K_2).$$

\begin{lemma}{\rm\cite[Theorem 3.1]{ggp}}\label{lemcutleaf}\qua
Let $G$ be an oriented framed $Y$\!--graph with leaves $K_1$, $K_2$, $K_3$ 
in a $\ZZ$--sphere $M$.
Assume that $K_3$ is a band sum of two framed knots $K^{2}_3$ and $K^{3}_3$. 
For $k=1$ and $2$, let $K_k^2$ and $K_k^3$ be two parallels of $K_k$ 
equipped with the framing $f(K_k)$ of $K_k$, and such that 
$\ell k(K_k^2,K_k^3)=f(K_k)$.
\begin{figure}[ht!]
\begin{center}
\centerline{\relabelbox \small
	\epsfysize 2cm \epsfbox{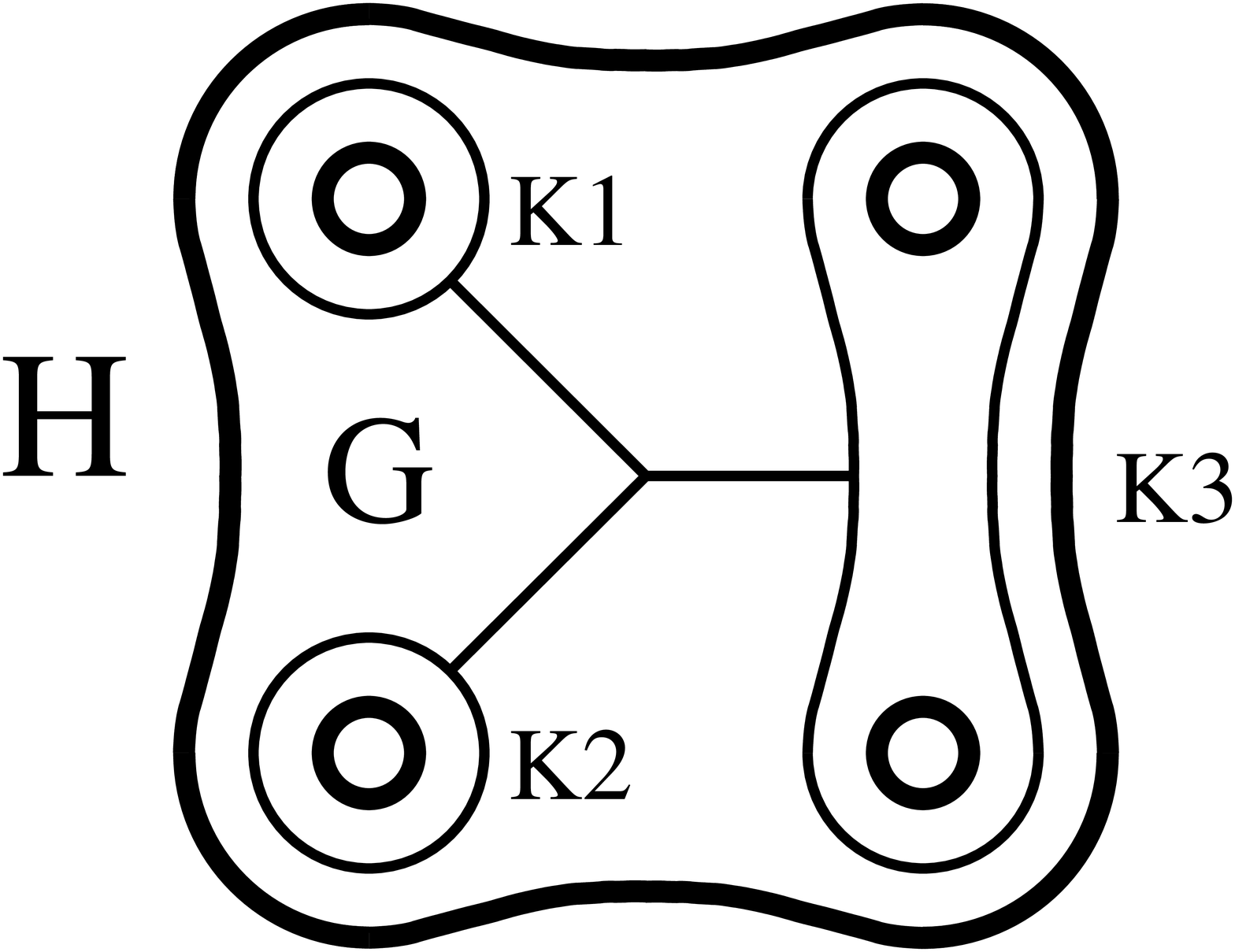}
\relabel {H}{$H$}
\relabel {G}{$G$}
\adjustrelabel <-1pt, 0pt> {K1}{$K_1$}
\adjustrelabel <-1pt, 0pt> {K2}{$K_2$}
\relabel {K3}{$K_3$}
\endrelabelbox
\hspace{40pt}
\relabelbox \small
\epsfysize 2cm \epsfbox{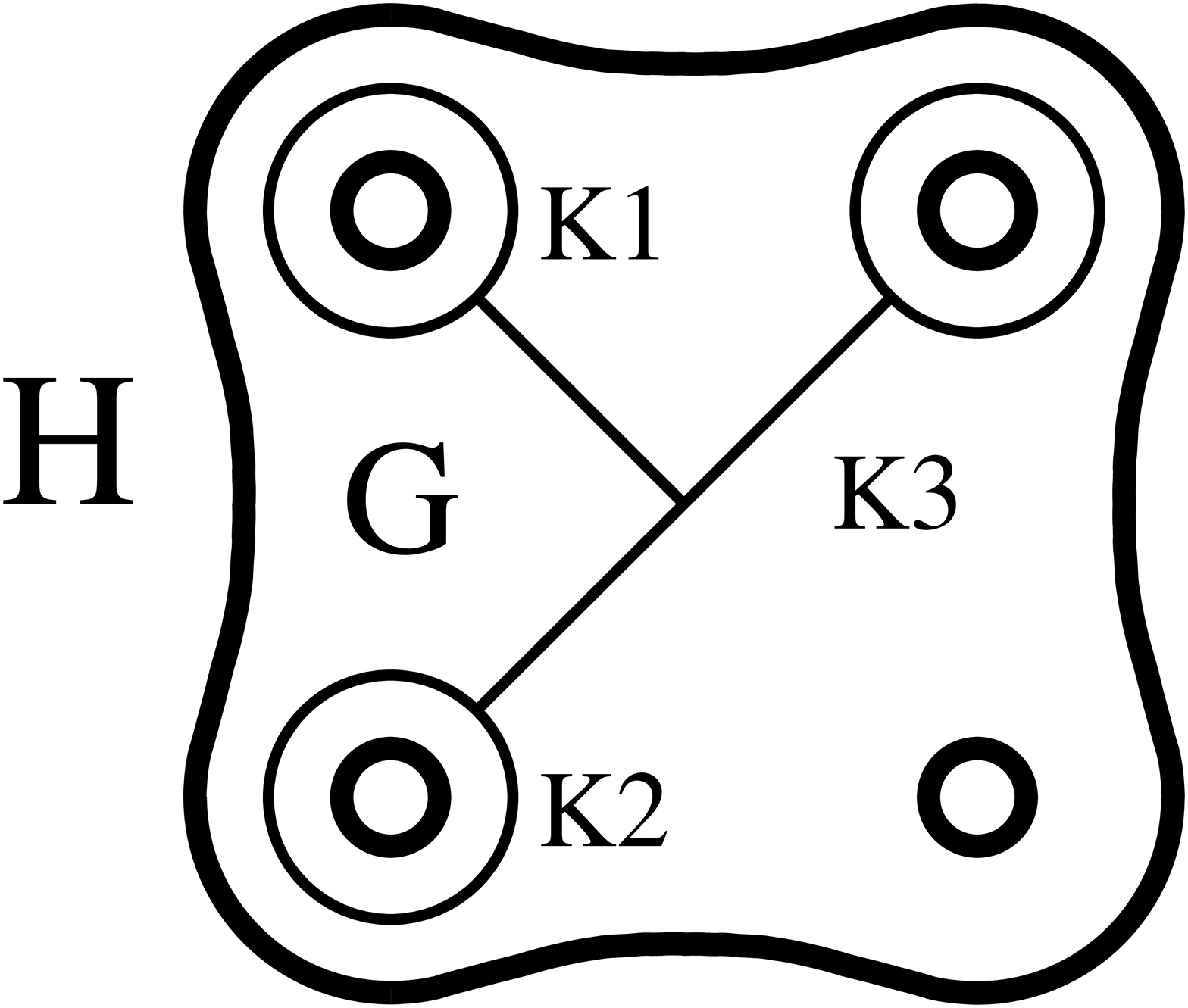}
\relabel {H}{$H$}
\relabel {G}{$G^2$}
\adjustrelabel <-1pt, 0pt> {K1}{$K_1^2$}
\adjustrelabel <-1pt, 0pt> {K2}{$K_2^2$}
\relabel {K3}{$K_3^2$}
\endrelabelbox
\hspace{10pt}
\relabelbox \small
\epsfysize 2cm \epsfbox{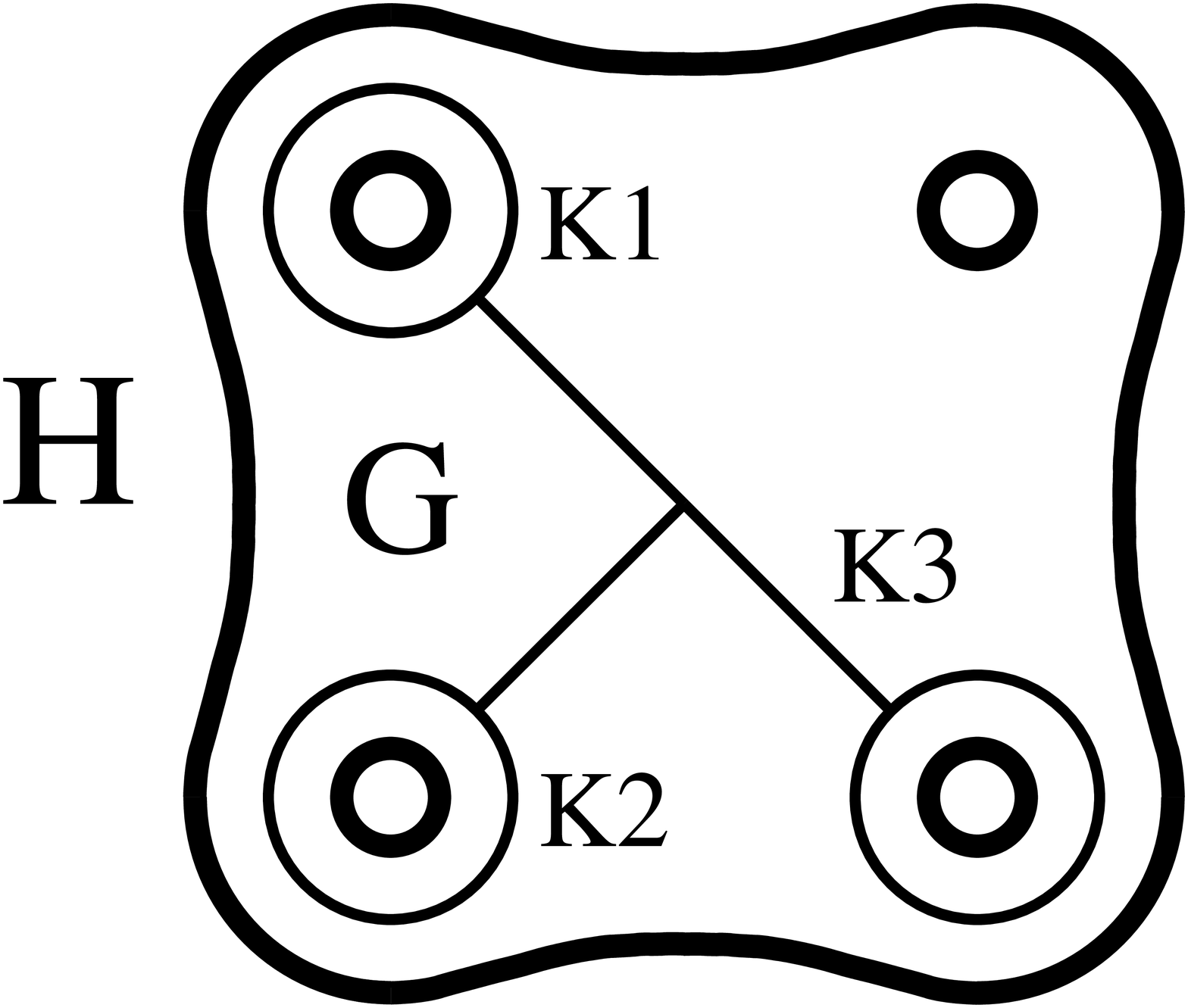}
\relabel {H}{$H$}
\relabel {G}{$G^3$}
\relabel {K1}{$K_1^3$}
\adjustrelabel <-1pt, 0pt> {K2}{$K_2^3$}
\adjustrelabel <-1pt, 0pt> {K3}{$K_3^3$}
\endrelabelbox }
\mbox{\small Figure \ref{fi6}: Splitting a leaf} 
\refstepcounter{figure} \label{fi6}
\end{center}
\end{figure}
Then
\begin{itemize}
\item[\rm(i)] There exist
two oriented disjoint framed $Y$\!--graphs $G^{2}$ and $G^{3}$ in $M$ whose 
framed leaves are
$K^{2}_1$, $K^{2}_2$, $K^{2}_3$ and $K^{3}_1$, $K^{3}_2$, $K^{3}_3$, 
respectively, such that the surgery along $G$ is equivalent to the surgery 
along $G^{2} \cup G^{3}$.
\item[\rm(ii)] For any $(n-1)$--component $Y$\!--link $L$ 
in the complement in $M$ of the embedded neighbourhood $H$ of $G$ represented
in \mbox{Figure \ref{fi6}},
$$\overline{[(M;L\cup G)]}=
\overline{[(M;L\cup G^2)]}+\overline{[(M;L \cup G^3)]}.$$
\end{itemize}
\end{lemma}

\proof
The surgery operation on $G$ is thought of as the move of two packs $P_2$ and $P_3$
of arcs of surgery components that go through the two holes on the right hand-side
of $H$ in \mbox{Figure \ref{fi6}} as in Remark~\ref{rksury}.
Then \mbox{Lemma \ref{llll}} says that the surgery along $G$ moves these
two framed packs of arcs by adding the boundary of a genus $1$ surface $\Sigma$. 
This 
operation can be made in two steps. Move $P_3$ first, 
which means do the surgery along a $Y$\!--graph $G^3$ whose leaves are
$K^3_1$, $K^3_2$ and $K^3_3$. Then move $P_2$ so that it 
is parallel to
$\partial \Sigma$ inside $\Sigma$. It can be done by a surgery along a $Y$\!--graph 
$G^2\subset H\setminus G^3$ whose leaves are $K^2_1$, $K^2_2$ and $K^2_3$. 
Then $M_{G^2\cup G^3}=M_G$. Therefore  
$$[(M;L\cup G^2\cup G^3)]=-[(M;L\cup G)]+[(M;L\cup G^2)]+[(M;L\cup G^3)]$$
and $$\overline{[(M;L\cup G)]}=\overline{[(M;L\cup G^2)]}+
\overline{[(M;L\cup G^3)]}.\eqno{\qed}$$   

\begin{lemma}{\rm\cite[Lemma 4.8]{ggp}}\label{framing}\qua
Let $G\subset M$ be an $n$--component $Y$\!--link. 
Suppose that a $Y$\!--component of $G$ contains
a $2$--framed leaf $l$ that bounds an embedded disc in $M\setminus G$.
Then $\overline{[M,G]}=0$.
\end{lemma}

\proof
If $l$ is a $2$--framed leaf that 
bounds an embedded disc in $M\setminus G$, then $l$ is a band sum of two knots
$K^2$ and $K^3$  that form a trivial Hopf link (see \mbox{Figure \ref{figHopf}}) in $M\setminus G$. 
Thanks to
\mbox{Lemma \ref{lemcutleaf}}, there exist two $n$--components $Y$\!--links 
$G^2$ and $G^3$ with a trivial leaf such that
$$\overline{[(M;G)]}=\overline{[(M;G^2)]}+\overline{[(M;G^3)]}=0.\eqno{\qed}$$

\bp[Proof of part (ii) of Proposition \ref{Review}]
Consider the bracket of $G$ as a\break function of a leaf $l$ of $G$ by fixing
$G\setminus l$. According to \mbox{Lemma \ref{firststep}}, the bracket of $G$ 
only depends on $[l]\in H_1(M\setminus \cup_{l^{\prime}\neq l}l^{\prime})$ and on $f(l)$.
Applying \mbox{Lemmas \ref{lemcutleaf}} and \ref{framing} 
when adding a disjoint $2$--framed
trivial knot to $l$ shows that this function of $l$ only depends on $[l]$ and on
$f(l)\  \mbox{\rm mod}\ 2$. Then \mbox{Lemma \ref{lemcutleaf}} implies (ii).
\eop

\begin{lemma}{\rm\cite[Lemma 2.3]{ggp}}\label{Hopf}\qua
Let $G\subset M$ be an $n$--component $Y$\!--link. Suppose that
$G$ contains a $Y$\!--component with two leaves $l$ and $l^{\prime}$ that form
the trivial Hopf link of \mbox{Figure \ref{figHopf}}. 
\begin{figure}[ht!]
\centerline{
\epsfysize 2cm \epsfbox{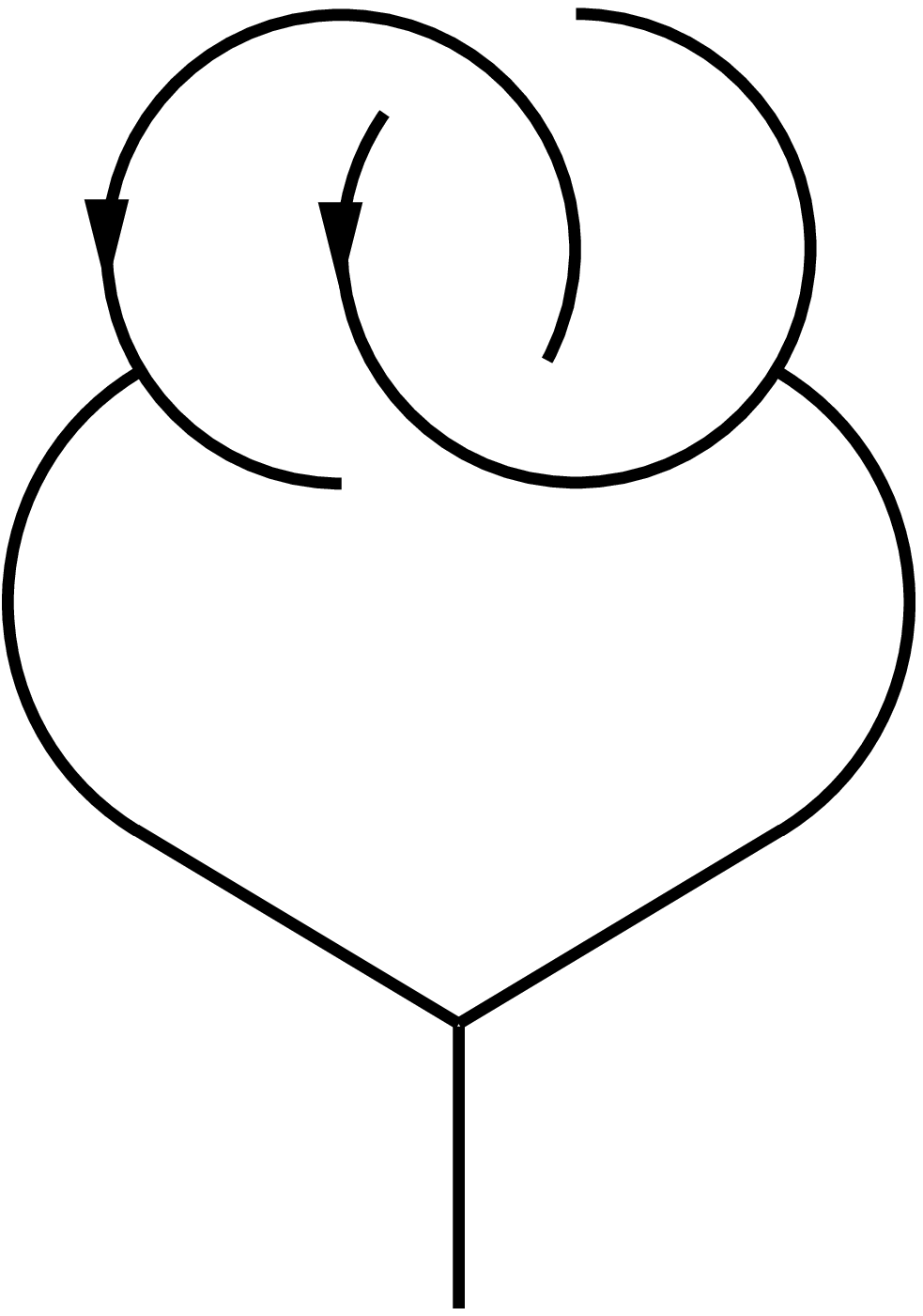}}
\caption{Two leaves that form the trivial Hopf link} \label{figHopf} 
\end{figure}
Then 
$\overline{[(M;G)]}=0$.
\end{lemma}

\bp
The surgery along a 
$Y$\!--graph with this trivial Hopf link is trivial: Think of this surgery as the 
move of surgery arcs along the boundary of the surface corresponding to these two
leaves as in \mbox{Lemma \ref{llll}}. It implies that the bracket of
$G$ vanishes.
\eop

\begin{definition}
Let $c$ be a curve in a surface $\Sigma$ that is the image of $S^1\times \{\pi\}$
under an orientation-preserving embedding 
$\phi \co S^1\times [0,2\pi]\longrightarrow \Sigma$. A {\em left-handed Dehn twist\/} of
$\Sigma$ along $c$ is the homeomorphism of $\Sigma$ that is the identity outside
$\phi(S^1\times ]0,2\pi[)$ and that maps $\phi(z,t)$ to $\phi(ze^{-it},t)$.
\end{definition}

\begin{lemma}{\cite[Theorem 3.1]{ggp}}\label{twist}\qua
Let $H$ be an oriented $Y$\!--graph in a $\ZZ$--sphere $M$.
Let $l^-$ and $l^{\prime}$ be two oriented leaves of $H$, and let $\Sigma$ be the genus 
one surface presented in \mbox{Figure \ref{lalere}}. Let $l$ be an oriented
parallel of $l^-$ in $\Sigma$ equipped with the framing induced by $\Sigma$. 
Let $l^{\prime \prime}$ be obtained from $l^{\prime}$ by a left-handed Dehn twist along $l$, and equipped
with the framing induced by the surface $\Sigma$, that is 
$f(l^{\prime \prime})=f(l)+f(l^{\prime})+2\ell k(l^-,l^{\prime})-1$.
Let $H^{\prime}$ be the $Y$\!--graph obtained from $H$ by changing $l^{\prime}$ into $l^{\prime \prime}$.
\begin{figure}[ht!]
\begin{center}
\refstepcounter{figure} \label{lalere}
\centerline{\relabelbox \small
\epsfysize 2cm \epsfbox{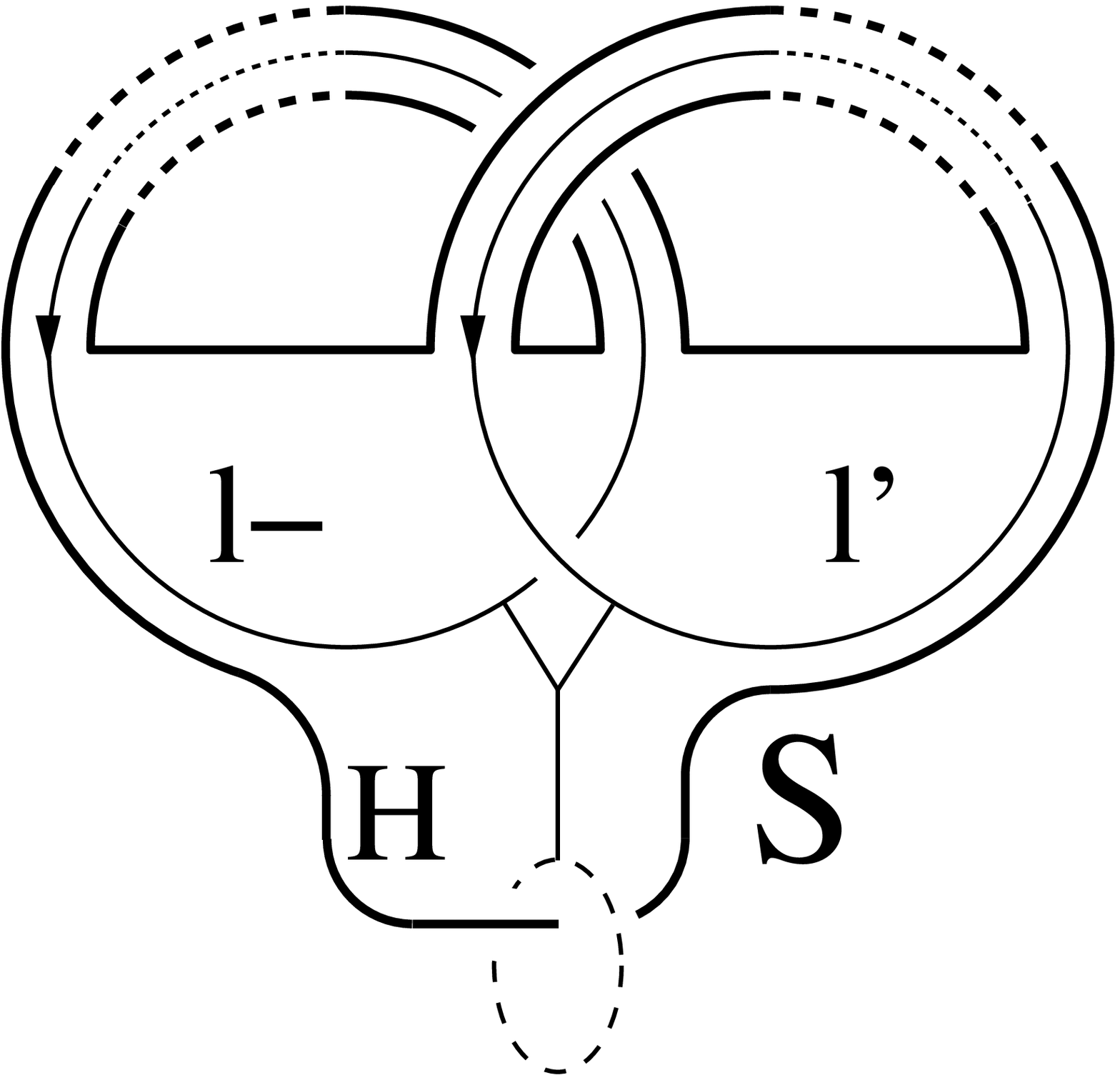}
\relabel {S}{$\Sigma$}
\relabel {H}{$H$}
\relabel {l'}{$l'$}
\relabel {l-}{$l^-$}
\endrelabelbox
\hspace{40pt}
\relabelbox \small
\epsfysize 2cm \epsfbox{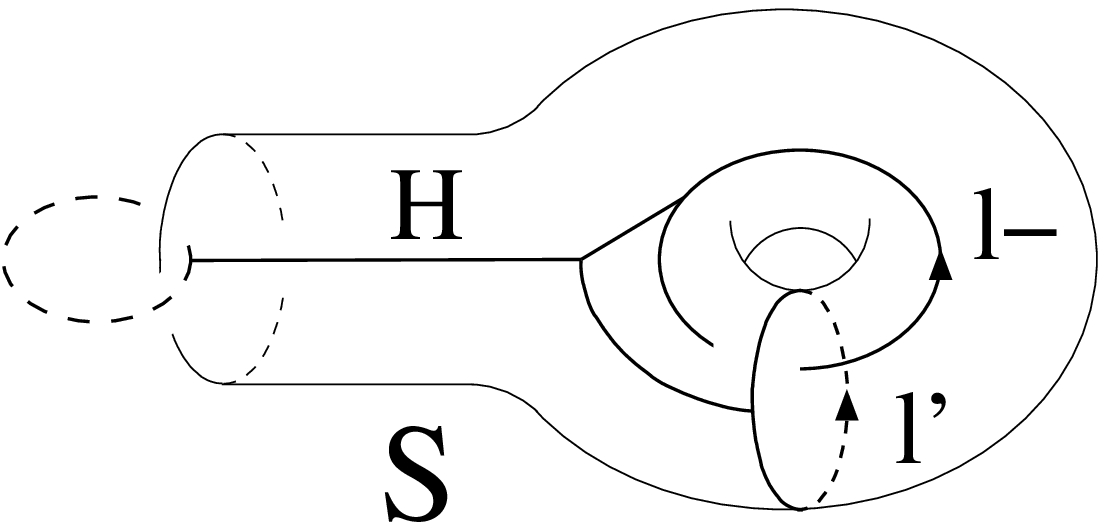}
\relabel {S}{$\Sigma$}
\relabel {H}{$H$}
\relabel {l'}{$l'$}
\relabel {l-}{$l^-$}
\endrelabelbox}
\mbox{\small Figure \ref{lalere}} 
\end{center}
\end{figure}
Then $M_H=M_{H^{\prime}}$ and, 
for any $Y$\!--link $L$ in the complement in $M$ of a 
neighbourhood of $H$, 
$${[(M;L\cup H^{\prime})]}={[(M;L\cup H)]}.$$
\end{lemma}

\bp
Thanks to \mbox{Lemma \ref{llll}}, the surgery on $H$ is uniquely determined by
$\Sigma$ that is unchanged by a Dehn twist  of $\Sigma$. It implies that $M_H=M_{H^{\prime}}$ and that, for any sublink $L(J)$ of $L$,
$M_{L(J)\cup H}=M_{L(J)\cup H^{\prime}}$. The equality of the brackets follows.
\eop

{\bf End of proof of Proposition \ref{Review}}
\begin{itemize}
\item Let $l$ be a leaf of an oriented $Y$\!--component $H$ of $G$. Let $l^{\prime}$ be the 
next leaf of $H$ (with respect to the cyclic order). We prove that increasing 
$\ell k(l,l^{\prime})$ does not change $\overline{[(M;G)]}$. 
By \mbox{Lemma \ref{lemcutleaf}}, adding a $0$--framed meridian $m_0(l^{\prime})$ of $l^{\prime}$ to $l$
adds $\overline{[(M;G(m_0(l^{\prime})/l))]}$ to $\overline{[(M;G)]}$, where 
$G(m_0(l^{\prime})/l)$ is 
obtained from $G$ by changing $l$ into $m_0(l^{\prime})$. Now, $l^{\prime}=m'+l^{\prime}_0$, where $l^{\prime}_0$
does not intersect a disk bounded by $m_0(l^{\prime})$, and $m'$ is a meridian of
$m_0(l^{\prime})$. Then $$\begin{array}{lll}\overline{[\big( M;G(m_0(l^{\prime})/l)\big)]}&=&
\overline{[\big( M;G( m_0(l^{\prime}),m'/l,l^{\prime})\big)]}\\
&&+\overline{[\big( M;G( m_0(l^{\prime}),l^{\prime}_0/l,l^{\prime}) \big)]}\\
&=&0\end{array}$$
since $G(m_0(l^{\prime}),m'/l,l^{\prime})$ is a $Y$\!--link with a trivial Hopf link
and since $G(m_0(l^{\prime}),l^{\prime}_0/l,l^{\prime})$ has a trivial leaf. 
\item To conclude, it is enough to show that if $l^-$ is a $0$--framed leaf of $G$,
if $l^{\prime}$ is the previous leaf in the component of $l^-$ in $G$ (w.r.t. the cyclic order), and if $l^{\prime}$ is $1$--framed, then changing
the framing of $l^{\prime}$ into $0$ does not change $\overline{[(M;G)]}$. By linearity, we 
may assume that $l^{\prime}$ is a trivial knot, and then it is enough to show that
$\overline{[(M;G)]}=0$.
By linearity on $l^-$, we can assume that $l^-$ is a 
$0$--framed meridian of some leaf
$l_0$ in another $Y$\!--component of $G$.
Let $G^{\prime}$ be the $Y$\!--link obtained from $G$ by changing $l^{\prime}$ into the twisted curve
$l^{\prime \prime}=l+l^{\prime}$ as in \mbox{Lemma \ref{twist}} 
so that 
$\bar{f}(l^{\prime \prime})=\bar{f}(l^{\prime})+1=0$. 
Then $l^-$ and $l^{\prime \prime}$ are $0$--framed meridians of $l_0$.
By 
linearity with respect to $l_0$, $\overline{[(M;G^{\prime})]}$ is the sum
of brackets of two $Y$\!--links with a trivial leaf. Then $\overline{[(M;G^{\prime})]}=\overline{[(M;G)]}=0$.\qed
\end{itemize}

\subsection{Proof of Theorem \ref{TT} for {\CL}P--surgeries induced
by $Y$\!--links} \label{induced}

Here we prove \mbox{Theorem \ref{TT}} when $D=(M;G)$ is an {\CL}P--surgery 
induced by a $Y$\!--link $G$.

\begin{lemma} \label{mute}
Let $A$ and $A^{\prime}$ be two $\ZZ$--handlebodies whose boundaries 
are identified so that $\CL_A=\CL_{A^{\prime}}$. Let $Z$ and $Z^{\prime}$ be two other 
$\ZZ$--handlebodies whose boundaries 
are identified so that $\CL_Z=\CL_{Z^{\prime}}$. 
Assume that $\partial A$ and $(-\partial Z)$ are identified so that
$A\cup_{\partial A} Z$ is a $\ZZ$--sphere.
Then
$$\mu(A\cup_{\partial A} Z)-\mu(A^{\prime}\cup_{\partial A} Z)=\mu(A\cup_{\partial A} Z^{\prime})-\mu(A^{\prime}\cup_{\partial A} Z^{\prime}).$$ 
\end{lemma}

\bp
For any $\ZZ$--sphere $M$, let $\lambda(M)$ be the {\em Casson invariant} of $M$.
Then $\mu(M)\equiv \lambda(M)\  \mbox{\rm mod} \;2$. Thanks to 
\mbox{\cite[Theorem 1.3]{les}}, 
$$\lambda (A\cup_{\partial A} Z)-\lambda (A^{\prime}\cup_{\partial A} Z)-(\lambda (A\cup_{\partial A} Z^{\prime})-\lambda (A^{\prime}\cup_{\partial A} Z^{\prime}))$$
is an even number. It implies the result.
\eop

Together with \mbox{Proposition \ref{Review}}, it
implies

\begin{corollary} \label{mu}
Let $H$ be a $Y$\!--graph in a $\ZZ$--sphere $M$. Let 
$$p=\bar{f}(l_1)\bar{f}(l_2)\bar{f}(l_3)$$
denote the product of the framings of 
the three leaves of $H$ in $\ZZ/2\ZZ$.
Then $$\mu(M_H)-\mu(M)=p.$$ 
\end{corollary}

\bp
First, we prove that $\mu$ vanishes on $\mathcal{F}_2$.
Let $G=G_1\cup G_2$ be a $2$--component $Y$\!--link in a $\ZZ$--sphere $M$. Let
$A$ be a regular neighbourhood of $G_1$. Let $Z$ be the complement of $\IT(A)$ in $M$.
Set $A^{\prime}=A_{G_1}$ and set $Z^{\prime}=Z_{G_2}$. Then $A$, $A^{\prime}$, $Z$ and $Z^{\prime}$ satisfy the 
assumptions of \mbox{Lemma \ref{mute}} and
$$\mu(M_G)-\mu(M_{G_1})-\mu(M_{G_2})+\mu(M)$$
$$= \mu(A^{\prime}\cup Z^{\prime})-\mu(A^{\prime}\cup Z)-
\mu(A\cup Z^{\prime})+\mu(A\cup Z)=0.$$
Then $\mu(M_H)-\mu(M)$ only depends  on 
$\overline{[(M;H)]} \in \CG_1$, thus it only depends on the product $p$ 
thanks to \mbox{Proposition \ref{Review}}, and it vanishes when $p=0$.

Then $\mu(M_H)-\mu(M)=p$
since $\mu$ is a non-trivial invariant on $\mathcal{F}$. 
\eop

\begin{lemma} \label{L2}
Let $\tilde{D}=(M;2k-1;(A_i,B_i)_{i=2,\dots,2k})$ be a $(2k-1)$--component
{\CL}P--surgery. Let $N$ be a 
$\ZZ$--handlebody in $M\setminus (\cup_{i=2}^{2k} A_i)$. Let $(A_1,B_1)$ be a pair of 
$\ZZ$--handlebodies such that
$A_1\subset \IT(N)$ and $\partial A_1$ and 
$\partial B_1$ are identified
so that $\CL _{A_1}=\CL _{B_1}$. Set $D(A_1,B_1)=\big( M;2k;(A_i,B_i)\big) $.
Let
$$i_{\ast} \co H_1(A_1)\longrightarrow H_1(N)$$ 
denote the homomorphism
induced by the inclusion map of $A_1$ into $N$.
Set 
$$\mathcal{J}_1(A_1,B_1)=\big( \otimes^3(i_{\ast}\circ 
\varphi_{A_1} ^{-1}) \big) 
\big( \CI(A_1,B_1)\big) \in \bigotimes_{i=1}^3 H_1(N)^{(i)}$$ 
where $\varphi_{A_1}$ has been defined in Notation~\ref{notation}.
Then, for any oriented degree $k$
Jacobi diagram $\Gamma$, there exists a linear form 
$\ell _N(\tilde{D};\Gamma)$ in $(\bigotimes ^3 H_1(N))^{\ast}$ such that,
for any pair $(A_1,B_1)$ as above,
$$\ell \big( D(A_1,B_1);\Gamma\big) =\langle\mathcal{J}_1(A_1,B_1)\ ,\  
\ell _N(\tilde{D};\Gamma)\rangle.$$
\end{lemma}

\proof
Let 
$\sigma$ be a coloration of 
$\Gamma$. Set
$$T'(\tilde{D};\Gamma;\sigma)=\mbox{\rm sign}(h)\ \bigotimes_{i=2}^{2k} 
 \CI (A_i, 
B_i)\in \otimes_{\{c\in H(\Gamma);\sigma\circ v(c)>1\}}X(c). $$ 
Apply all the contractions corresponding to the edges that do not contain 
any half-edge $c$ in $(\sigma\circ v)^{-1}(1)$. These are all the edges except the three edges
\{c,j(c)\} where $c \in (\sigma\circ v)^{-1}(1)$ and $j(c)$ is the other half-edge of $e(c)$, and the obtained tensor is
in $\otimes_{\{c\in H(\Gamma);\sigma\circ v(c)=1\}}X(j(c))$.
Now apply $\bigotimes_{\{c\in H(\Gamma);\sigma\circ v(c)=1\}} \varphi_{A_{\sigma\circ v(j(c))}}^{-1}$ in order to obtain the
tensor
$$\ell^{\prime}(\tilde{D};\Gamma;\sigma)\in \bigotimes_{\{c\in H(\Gamma);\sigma\circ v(c)=1\}}
H_1(A_{\sigma\circ v(j(c))}).$$ 
The linking number maps 
$H_1(A_{\sigma\circ v(j(c))})$ to $\big( H_1(N)\big) ^{\ast}$ and therefore maps 
$\ell^{\prime}(\tilde{D};\Gamma;\sigma)$ to an element 
$$\ell _N(\tilde{D};\Gamma;\sigma)\in  \bigotimes_{i=1}^3\big( H_1(N)^{\ast}
\big) ^{(i)}.$$
By definition, $\ell(\tilde{D};\Gamma;\sigma)$ is the contraction of 
$\ell _N(\tilde{D};\Gamma;\sigma)\otimes \mathcal{J}_1(A_1,B_1)$.
Then 
$$ \ell (D;\Gamma)= \langle\mathcal{J}_1(A_1,B_1)\ ,\  
\sum_{\sigma \in \mbox{\rm Bij}(\Gamma)/\mbox{\rm Aut}(\Gamma)} 
\frac{\ell _N (\tilde{D};\Gamma; \sigma)}{\sharp \mbox{Aut}_V(\Gamma)}
\rangle.\eqno{\qed}$$

\begin{proposition}\label{Ylink}
Let $n\in \mathbb{N}$. Let $M$ be a $\ZZ$--sphere.
Let $\Gamma$ be an oriented degree $k$ Jacobi diagram with $2k\leq n$. 
Consider $\ell((M;G);\Gamma)$ as a function of oriented $n$--component $Y$\!--links
$G$ in $M$.
Then 
\begin{itemize}
\item The linking number 
$\ell((M;G);\Gamma)$ only depends on
\begin{itemize}
\item  the linking numbers $\ell k(l,l^{\prime})$, where $l$ and $l^{\prime}$ are leaves in two 
distinct
$Y$\!--components of $G$
\item  the products $\bar{f}(l_1)\bar{f}(l_2)\bar{f}(l_3)$, where  
$l_1$, $l_2$ and $l_3$ are leaves in a 
same $Y$\!--component.
\end{itemize}
\item Considered as a map of a leaf $l$ of $G$,
$\ell((M;G);\Gamma)$ is a linear map in 
$([l],\bar{f}(l))\in H_1(M\setminus \cup_{l^{\prime}\neq l}l^{\prime})\times \ZZ/2\ZZ$.
\end{itemize}
\end{proposition}

\bp
Assume $n=2k$. \mbox{Lemma \ref{L2}} and the expression of 
$\CI(A_1,B_1)$
given in Subsection~\ref{S26} 
show that
$\ell((M;G);\Gamma)$ does not depend on the framing of 
the leaves of $G$ and that 
$\ell((M;G);\Gamma)$, seen as a map of a leaf $l$ of a component $H$ of $G$, 
linearly depends on the
homology class of $l$ in $M\setminus(G\setminus H)$. 
It implies the result. If $2k<n$, the result follows from
\mbox{Corollary \ref{mu}} and from the previous case.
\eop

\begin{proposition} \label{moitie}
Let $G$ be an $n$--component $Y$\!--link in a $\ZZ$--sphere $M$.
Then Theorem \ref{TT} is true when $D=(M;G)$.
\end{proposition}

\bp
The simultaneous multilinearities 
of the bracket in \mbox{Proposition \ref{Review}} 
and of the linking number of {\CL}P--surgeries induced by $Y$\!--links 
in \mbox{Proposition \ref{Ylink}} allow us to cut
the leaves of $G$ and to reduce the proof in the case where
\begin{itemize}
\item The non-zero-framed leaves are $\pm1$--framed and bound discs disjoint from 
$G$.
\item Any $0$--framed leaf is a meridian of another leaf.
\end{itemize}
If a $\pm 1$--framed leaf is in a component with a $0$--framed leaf, its framing can be 
changed without changing either side of the equality. Then we can assume than 
the only $\pm 1$--framed leaves are parts of components like $Y_{I\!I\!I}$.
Similarly, we can assume that any $0$--framed leaf is a meridian of one leaf in 
another component of $G$. 
Then, up to orientation changes of leaves, we can assume that $G$ is a $Y$\!--link
induced by a Jacobi diagram.
Since \mbox{Theorem \ref{TT}} is true for {\CL}P--surgeries induced by Jacobi 
diagrams, the result follows.
\eop
 

\subsection{Proof of Theorem \ref{TT} in the general case}

Let $(A,B)$ be a pair of $\ZZ$--handlebodies whose boundaries
$\partial A$ and $\partial B$  are identified so that 
$\CL_A=\CL_B$. In what follows, 
$\CI(A,B)$ denotes the linear form on
 $\bigotimes^3\CL_A$ 
induced by the intersection
form on $\bigotimes^3H_2\big( A\cup(-B)\big) $, and
$$\varphi_A \co  H_1(A)\longrightarrow \CL_A^{\ast}$$
denotes the isomorphism presented in \mbox{Notation \ref{notation}}. 

\begin{lemma} \label{L1}
Let $A$, $B$ and $C$ be three $\ZZ$--handlebodies with the same genus.
Assume that $\partial A$, $\partial B$ and $\partial C$ are identified so that
$\CL_A=\CL_B=\CL_C$. Then 
$$\CI(A,B)=\CI(A,C)+\CI(C,B).$$
\end{lemma}

\bp
Let $a_1$, $a_2$ and $a_3$ be three oriented curves in $\partial A$ 
that represent 
elements 
of $\CL_A$ still denoted by $a_1$, $a_2$ and $a_3$ such that the curves 
$a_i$ do not intersect each other. 
Let $M=A\cup(-B)$.
For any $i\in \{1,2,3\}$, 
let $S^i_A$ and $S^i_B$ be oriented surfaces in $A$ and in $B$, respectively,
such that
$a_i$ bounds $S^i_A$ in $A$ and $a_i$ bounds $S^i_B$ in $B$.
Assume that all the surfaces are transverse to each other and to $\partial A$.
Set 
$$\Sigma_M^i= S^i_A\cup_{a_i} (-S^i_B) \subset M.$$
The orientation of $\Sigma_M^i$
and the orientation of $M$ induce a positive normal vector field $n_i$ on $\Sigma_M^i$.
The algebraic intersection $\langle\Sigma_M^1,\Sigma_M^2, \Sigma_M^3\rangle_M$
is the sum of the signs of the intersection
points, where the sign is 
defined as follows. For any intersection point $x$, the sign is $+1$ 
if $\big( n_1(x),n_2(x),n_3(x)\big) $ is a direct basis of $T_xM$ according to 
the orientation of $M$, and $-1$ otherwise. By definition,  
$$\big( \CI(A,B)\big) \big( a_1 \otimes a_2 \otimes a_3\big)\  =\ 
\langle\Sigma_M^1,\Sigma_M^2, \Sigma_M^3\rangle_M.$$
Then
$$\begin{array}{rl}
\big( \CI(A,B)\big) \big( a_1 \otimes a_2 \otimes a_3\big) &=\ \ \ 
\langle\Sigma_M^1,\Sigma_M^2, \Sigma_M^3\rangle_M\\
=& \langle S_A^1,S_A^2, S_A^3\rangle_M + \langle(-S_B^1),(-S_B^2), (-S_B^3)\rangle_M\\
=& \langle S_A^1,S_A^2, S_A^3\rangle_A + \langle(-S_B^1),(-S_B^2), (-S_B^3)\rangle_{(-B)}.
\end{array}$$
Note that the normal vector field $n_B^i$ to $S_B^i$ induced by the 
orientation of $S_B^i$ and by the orientation of $B$ is equal to
the normal vector field to $(-S_B^i)$ induced by 
the orientation
of $(-S_B^i)$ and by the orientation of $(-B)$.
Now, for each point of $S_B^1\cap S_B^2\cap S_B^3$,
$(n_B^1,n_B^2,n_B^3)$ is direct according to the orientation of $B$ if
and only if it is not direct according to the orientation of $(-B)$. 
It implies that
$$ \big( \CI(A,B)\big) \big( a_1 \otimes a_2 \otimes a_3\big)\  =\ 
\langle S_A^1,S_A^2, S_A^3\rangle_A - \langle S_B^1,S_B^2,S_B^3 \rangle_B$$
and the lemma follows. \eop

\begin{lemma} \label{l11}
Under the hypotheses of \mbox{Lemma \ref{ajout}}, 
for any Jacobi diagram $\Gamma$,
$$\ell(D;\Gamma)=\ell(D';\Gamma)+ \ell(D'';\Gamma).$$
\end{lemma}
 
\proof
The result follows from the equality $\CI(A_1,B_1)=\CI(A_1,A^{\prime}_1)+\CI(A^{\prime}_1,B_1)$
given by \mbox{Lemma \ref{L1}} and from the equality
$$\CL (D(\{1\}))=\CL (D'(\{1\}))+\CL (D''(\{1\})).\eqno{\qed}$$

\begin{lemma} \label{l10}
Consider an $n$-component {\CL}P--surgery 
$$D=(M;n;(A_1,B_1),(A_2,B_2),\dots (A_n,B_n)).$$
Let $A^{\prime}_1$ and $B^{\prime}_1$ be 
two $\ZZ$--handlebodies such that
\begin{itemize}
\item $A^{\prime}_1\subset \IT(A_1)$
\item $\partial A^{\prime}_1$ and $\partial B^{\prime}_1$  are identified so that
$\CL_{A^{\prime}_1}=\CL_{B^{\prime}_1}$
\item $B_1=(A_1)_{B^{\prime}_1/A^{\prime}_1}$ is the $\ZZ$--handlebody obtained
from  $A_1$ by replacing $A^{\prime}_1$ by $B^{\prime}_1$. 
\end{itemize}
Set
$D'=\big( M;n; (A^{\prime}_1,B^{\prime}_1),(A_2,B_2),\dots, (A_{n},B_{n})\big) $. 
Then 
$[D']=[D]$ while, for any Jacobi diagram $\Gamma$,
$$\ell (D';\Gamma)=\ell(D;\Gamma).$$
\end{lemma}

\begin{sublemma} \label{l30}
Under the hypotheses of \mbox{Lemma \ref{l10},} let 
$$\begin{array}{cccc}
\partial_{A_1} \co  & H_2(A_1,\partial A_1) &\to  &\CL_{A_1}\\ 
\partial_{A^{\prime}_1} \co  & H_2(A^{\prime}_1,\partial A^{\prime}_1) &\to  &\CL_{A^{\prime}_1}
\end{array}$$ 
denote the isomorphisms
induced by the long exact homology sequences. 
Let 
$$\begin{array}{cccc}
i_{A_1} \co     &H_2(A_1,\partial A_1)      &\to    &H_2\big( A_1, 
A_1\setminus \IT (A^{\prime}_1)
\big)\\
i_{A^{\prime}_1} \co  &H_2(A^{\prime}_1,\partial A^{\prime}_1)    &\to    &H_2\big( A_1, A_1\setminus 
\IT (A^{\prime}_1)
\big)
\end{array}$$
be the homomorphisms induced by the inclusion maps. 
Then $i_{A^{\prime}_1}$ is
an isomorphism by the excision axiom.  
Set $$\Phi=\partial_{A^{\prime}_1}\circ i_{A^{\prime}_1}^{-1}\circ 
i_{A_1} \circ \partial _{A_1}^{-1} \co \CL_{A_1} \to \CL_{A^{\prime}_1}.$$
Then
$$\big( \CI(A^{\prime}_1,B^{\prime}_1) \big) \circ (\otimes ^3 \Phi )=\CI(A_1,B_1).$$  
\end{sublemma}

\bp{\bf of \mbox{Lemma \ref{l10}} assuming \mbox{Sublemma \ref{l30}}} \;
The assertion $[D']=[D]$ is obvious. 
Since $\mu(M_{B_1/A_1})=\mu(M_{B^{\prime}_1/A^{\prime}_1})$, it suffices to prove that 
$$\ell (D;\Gamma)=\ell(D';\Gamma)$$
when $n=2k$ is even and when $\Gamma$ is a degree $k$ Jacobi diagram.
Set 
$\tilde{D}=\big( M;2k-1;(A_i,B_i)_{i=2,\dots , 2k}\big)$.
Let $$i_{\ast} \co H_1(A^{\prime}_1) \longrightarrow H_1(A_1)$$ 
be the map induced by the inclusion map
of $A^{\prime}_1$ into $A_1$.
By Lemma \ref{L2}, there exists a linear form 
$$\ell _{A_1}\big( \tilde{D};\Gamma\big) \in \big( \bigotimes^3 
H_1(A_1)\big) ^{\ast}$$ 
such that
$$\begin{array}{ccc}
\ell (D;\Gamma) &=      &\langle\big( \otimes ^3 \varphi _{A_1}^{-1}\big)
\big( \CI (A_1,B_1)\big) \ ,\  \ell _{A_1}(\tilde{D};\Gamma)\rangle\\ 
\ell (D';\Gamma)        &=&
\langle\big( \otimes ^3 (i_{\ast}\circ \varphi _{A^{\prime}_1}^{-1})\big) \big( \CI (A^{\prime}_1,B^{\prime}_1)
\big) ,
\ell _{A_1}(\tilde{D};\Gamma)\rangle.
\end{array}$$
The following diagram is commutative
$$\begin{array}{ccccc}
H_1(A^{\prime}_1) & \ &\stackrel{\varphi _{A^{\prime}_1}}{\longrightarrow} & \ &
\CL _{A^{\prime}_1}^{\ast}\\
&&&& \\
i_{\ast} \downarrow &&      & &\downarrow \Phi^{\ast}\\
&&&& \\
H_1(A_1)  && \stackrel{\varphi _{A_1}}{\longrightarrow} & &\CL _{A_1}^{\ast}.
\end{array}$$
Indeed both compositions, seen as elements of 
$\big( H_1(A^{\prime}_1)\otimes \CL_{A_1}\big) ^{\ast},$
map $([b],[a]) \in H_1(A^{\prime}_1)\times \CL_{A_1}$ to the 
algebraic intersection in $A_1$ of a surface bounded
by $a$ and of the curve $b$. 
Thus, 
$$\big( \otimes^3 (i_{\ast}\circ \varphi _{A^{\prime}_1}^{-1})\big)
\big(\CI (A^{\prime}_1,B^{\prime}_1)\big)=\big( \otimes^3 \varphi _{A_1}^{-1}\big) \big( \CI (A^{\prime}_1,B^{\prime}_1) \circ 
(\otimes^3 \Phi)\big).$$ 
By Sublemma \ref{l30}, 
the lemma follows.
\eop

\bp{\bf of Sublemma \ref{l30}} \;
Let $a_1$, $a_2$ and $a_3$ be three  oriented curves in $\partial A_1$
that represent elements of $\CL_{A_1}$ still denoted by $a_1$, $a_2$ and $a_3$ 
such 
that the 
curves $a_i$ do not intersect each other. Let $a'_1,a'_2,a'_3$ be oriented 
curves in $\partial A^{\prime}_1$
that represent the elements $\Phi(a_i)$ of $\CL_{A^{\prime}_1}$ such that the 
curves $a'_i$ do not intersect each other. 
For any $i\in\{1,2,3\}$, the curves $a_i$ and $a'_i$ cobound an oriented 
surface
$\Sigma^i$ in $A_1\setminus \IT(A^{\prime}_1)$. The curve $a'_i$ bounds an oriented 
surface
$\sigma_{A^{\prime}_1}^i$ in $A^{\prime}_1$ and bounds an oriented surface $\sigma_{B^{\prime}_1}^i$ in 
$B^{\prime}_1$.
Set
$$\begin{array}{cccc}
S_i=    & \big( \sigma_{A^{\prime}_1}^i\cup \Sigma^i\big) \cup 
\big( -(\Sigma^i\cup \sigma_{B^{\prime}_1}^i)\big) &
\subset & A_1\cup (-B_1)\\
S'_i=   &\sigma_{A^{\prime}_1}^i \cup (-\sigma_{B^{\prime}_1}^i)    & \subset & A^{\prime}_1\cup (-B^{\prime}_1).
\end{array}$$
Set
$$\begin{array}{ccc}
\mathcal{J}_{A^{\prime}_1B^{\prime}_1}&=&\big( \CI(A^{\prime}_1,B^{\prime}_1)\big) \big( 
\Phi(a_1)\otimes \Phi(a_2) \otimes
\Phi(a_3)\big) \\
\mathcal{J}_{A_1B_1}&=&\big( \CI(A_1,B_1)\big) 
\big( a_1\otimes a_2 \otimes a_3\big) .
\end{array}$$
By definition, 
$\mathcal{J}_{A_1B_1}$ is the intersection in $A_1\cup (-B_1)$ of the oriented
surfaces $S_i$ and $\mathcal{J}_{A^{\prime}_1B^{\prime}_1}$ is the intersection in 
$A^{\prime}_1\cup (-B^{\prime}_1)$ of the oriented surfaces $S'_i$.
Then $\left(\mathcal{J}_{A_1B_1}-\mathcal{J}_{A^{\prime}_1B^{\prime}_1}\right)$ is the contribution of the 
intersection of the surfaces $\Sigma^i\cup (-\Sigma^i)$. 
This contribution vanishes when $A^{\prime}_1=B^{\prime}_1$ because 
$\CI(A_1,A_1)= \CI(A^{\prime}_1,A^{\prime}_1)=0$ by  \mbox{Lemma \ref{L1}.}
Hence it always vanishes. \eop

\bp[Proof of \mbox{Theorem \ref{TT}}]
Lemmas \mbox{\ref{V1},} \mbox{\ref{ajout}} and  
\mbox{\ref{l11}} allow us to reduce the proof to the case of an {\CL}P--surgery
$D=\big( M;n;(A_i,B_i) \big)$ where $B_i$ is obtained from $A_i$ by a surgery on 
a $Y$\!--graph embedded in $A_i$.
By Lemma \ref{l10}, $D$ can next be considered as an {\CL}P--surgery
induced by an $n$--component $Y$\!--link
in $M$. Then 
\mbox{Theorem \ref{TT}} follows from \mbox{Proposition \ref{moitie}}.
\eop

\appendix
\section{The  IHX relation}
\setcounter{equation}{0}

For self-containedness, we finish the proof of Theorem \ref{base} by proving that 
$\Psi_n$ factors through the IHX relation. This is a consequence of the following 
proposition.
\begin{proposition} \label{ihx}
Let $G_1$ be an oriented $Y$\!--link in a $\ZZ$--sphere $M$ that admits the following 
two-component 
sublink $I_1$.
For $i=2,3$, let $G_i$ be obtained from $G_1$ by changing
$I_1$ into $I_i$. The four leaves $a$, $b$, $c$ and $d$ are identical.
\begin{center}
\begin{pspicture}[.4](-1.6,-.2)(2,3.3) 
\rput(1.3,.1){\small $I_1$}
\psline{-*}(0,.6)(0,.9)
\rput[l](-.4,2.8){\small $c$}
\psarc{->}(-.8,2.8){.3}{-180}{0}
\psarc[linestyle=dashed,dash=2pt 2pt](-.8,2.8){.3}{0}{180}
\rput[l](1.2,2.8){\small $b$}
\psarc{->}(.8,2.8){.3}{-180}{0}
\psarc[linestyle=dashed,dash=2pt 2pt](.8,2.8){.3}{0}{180}
\rput[l](.4,.3){\small $a$}
\psarc[linestyle=dashed,dash=2pt 2pt]{->}(0,.3){.3}{-180}{0}
\psarc(0,.3){.3}{0}{180}
\rput[b](-.9,1.3){\small $d$}
\psarc{->}(-.9,.9){.3}{-90}{90}
\psarc[linestyle=dashed,dash=2pt 2pt](-.9,.9){.3}{90}{-90}
\psline{-*}(-.8,2.5)(0,2.3)
\psline{-}(.8,2.5)(0,2.3)
\psline{-}(0,1.9)(0,2.3)
\psline{-}(0,.9)(0,1.3)
\psline(-.6,.9)(0,.9)
\psarc{->}(0,1.5){.2}{60}{0}
\psarc{->}(0,1.7){.2}{-120}{180}
\end{pspicture}
\begin{pspicture}[.4](-2,-.2)(1.6,3.3) 
\rput(1.3,.1){\small $I_2$}
\psline{-*}(0,.6)(0,.9)
\rput[l](-.4,2.8){\small $c$}
\psarc{->}(-.8,2.8){.3}{-180}{0}
\psarc[linestyle=dashed,dash=2pt 2pt](-.8,2.8){.3}{0}{180}
\rput[l](1.2,2.8){\small $b$}
\psarc{->}(.8,2.8){.3}{-180}{0}
\psarc[linestyle=dashed,dash=2pt 2pt](.8,2.8){.3}{0}{180}
\rput[l](.4,.3){\small $a$}
\psarc[linestyle=dashed,dash=2pt 2pt]{->}(0,.3){.3}{-180}{0}
\psarc(0,.3){.3}{0}{180}
\psline{-*}(-.8,2.5)(0,.9)
\psline{-}(.15,1.2)(0,.9)
\psline{-}(.45,1.8)(.6,2.1)
\pscurve[border=1pt]{-}(-.6,.9)(0,.75)(1.3,1.6)(.6,2.1)
\psline{-*}(.8,2.5)(.6,2.1)\rput[b](-.9,1.3){\small $d$}
\psarc{->}(-.9,.9){.3}{-90}{90}
\psarc[linestyle=dashed,dash=2pt 2pt](-.9,.9){.3}{90}{-90}
\psarc{-}(.35,1.6){.17}{-140}{165}
\psarc{-}(.25,1.4){.17}{45}{-10}
\end{pspicture}
\begin{pspicture}[.4](-2,-.2)(2,3.3) 
\rput(1.3,.1){\small $I_3$}
\psline{-*}(0,.6)(0,.9)
\rput[l](-.4,2.8){\small $c$}
\psarc{->}(-.8,2.8){.3}{-180}{0}
\psarc[linestyle=dashed,dash=2pt 2pt](-.8,2.8){.3}{0}{180}
\rput[l](1.2,2.8){\small $b$}
\psarc{->}(.8,2.8){.3}{-180}{0}
\psarc[linestyle=dashed,dash=2pt 2pt](.8,2.8){.3}{0}{180}
\rput[l](.4,.3){\small $a$}
\psarc[linestyle=dashed,dash=2pt 2pt]{->}(0,.3){.3}{-180}{0}
\psarc(0,.3){.3}{0}{180}
\psline{-*}(.8,2.5)(0,.9)
\psline{-}(-.15,1.2)(0,.9)
\psline{-}(-.45,1.8)(-.6,2.1)
\pscurve[border=1pt]{-}(-.6,.9)(0,.75)(1.3,1.6)(-.6,2.1)
\psline{-*}(-.8,2.5)(-.6,2.1)
\rput[b](-.9,1.3){\small $d$}
\psarc{->}(-.9,.9){.3}{-90}{90}
\psarc[linestyle=dashed,dash=2pt 2pt](-.9,.9){.3}{90}{-90}
\psarc{-}(-.35,1.6){.17}{-80}{220}
\psarc{-}(-.25,1.4){.17}{100}{45}
\end{pspicture}
\end{center}

Then $$\overline{[(M;G_1)]} + \overline{[(M;G_2)]}+ \overline{[(M;G_3)]}=0.$$
\end{proposition}
\bp First recall that Proposition~\ref{Review} implies that the actual embeddings of the $Y$\!--graph edges
do not affect the $Y$\!--link brackets. Therefore these embeddings will not be specified in the proof below.
Consider an embedded product $D_{abc} \times [0,1]$ of the three-hole disk $D_{abc}$ whose three
inner boundary components are $a$, $b$ and $c$, and the three generators of $\pi_1(D_{abc})$,
$\alpha$, $\beta$ and $\gamma$. We first prove that there exists some two-component $Y$\!--link $J_1$ in 
$\left(M \setminus (G_1 \setminus I_1)\right)$ that is obtained from $I_1$ by changing (the edge adjacent to $d$ and) the leaves $a$, $b$ and $c$ into leaves 
$a(1)$, $b(1)$ and $c(1)$ that are homologous to $a$, $b$ and $(-c)$ 
 in $\left(M \setminus (G_1 \setminus I_1)\right)$, respectively, 
such that surgery on $J_1$ makes
a pack $T$ of surgery arcs in a surgered disk bounded by $d$
describe the element $[\beta\alpha\beta^{-1},[\gamma,\beta]]$ of $\pi_1(D_{abc})$ in $D_{abc} \times [0,1]$ in an ascending way with respect to the height of $[0,1]$.

Indeed the second part of the following picture shows such a path that is ascending (after sliding two tongues), and that cobounds a genus one surface 
with the initial shown portion of $T$. 
\begin{center}
\begin{pspicture}[.4](-2,-1.4)(2,1.4) 
\psecurve[linestyle=dashed,dash=3pt 2pt]{->}(-1,-1)(0,0)(.5,.9)(.9,.5)
\psecurve[linestyle=dashed,dash=3pt 2pt](0,0)(.5,.9)(.9,.5)(.5,.1)(0,0)(-1,-1)
\psecurve[linestyle=dashed,dash=3pt 2pt]{->}(-1,1)(0,0)(.5,-.1)(.9,-.5)(.5,-.9)(0,0)
\psecurve[linestyle=dashed,dash=3pt 2pt](.9,-.5)(.5,-.9)(0,0)(-1,1)
\psecurve[linestyle=dashed,dash=3pt 2pt]{->}(1,-1)(0,0)(-.5,.1)(-.9,.5)(-.5,.9)(0,0)
\psecurve[linestyle=dashed,dash=3pt 2pt](-.9,.5)(-.5,.9)(0,0)(1,-1)
\psecurve{->}(0,0)(-.2,-.2)(-.9,-.5)(-.9,-.9)(-.5,-.9)
\psccurve[border=.1](-1.3,1.3)(0,1.3)(1.3,1.3)(1.3,0)(1.3,-1.3)(-.5,-.5)
\psline[border=.1]{->}(-.7,-.7)(-1.1,-.3)
\psline[border=.1]{>-}(-.3,-1.1)(-.4,-1)
\psecurve[border=.1](-.9,-.5)(-.9,-.9)(-.5,-.9)(-.2,-.2)(0,0)
\psline{*-}(-.2,-.2)(0,0)
\psarc{<-}(-.5,.5){.2}{-180}{180}
\psarc{<-}(.5,.5){.2}{-180}{180}
\psarc{<-}(.5,-.5){.2}{-180}{180}
\rput[l](-.6,.5){\small $c$}
\rput[l](.4,.5){\small $b$}
\rput[l](.4,-.5){\small $a$}
\rput[rt](-1.2,-.4){\small $T$}
\rput[rt](-1,-1){\small $d$}
\rput[b](.5,1){\small $\beta$}
\rput[b](-.5,1){\small $\gamma$}
\rput[t](.6,-1){\small $\alpha$}
\rput[r](1.3,0){\small $D_{abc}$}
\end{pspicture}
\begin{pspicture}[.4](-4,-1.5)(4,4) 
\psline[doubleline=true, linearc=.3,doublesep=.2](0,.6)(0,3.8)(3.8,3.8)(3.8,-1.2)(1.5,-1.2)(1.5,-.2)(3.4,.7)(3.4,3.4)(1,3.4)(1,.6)
\pspolygon[linearc=.3,linestyle=dashed,dash=2pt 5pt](.4,.3)(0,.6)(0,3.8)(3.8,3.8)(3.8,-1.2)(1.5,-1.2)(1.5,-.2)(3.4,.7)(3.4,3.4)(1,3.4)(1,.6)
\psline[border=.1,doubleline=true, linearc=.3,doublesep=.2](3,2)(3,2.8)(0.5,3.3)(0.5,.6)
\psecurve(.1,.9)(.1,.6)(.4,.6)(.4,.9)
\psecurve(.6,.9)(.6,.6)(.9,.6)(.9,.9)
\psline[linearc=.3]{>-}(1,-1.2)(-.1,-.3)(-.1,-.2)
\psline(-.1,0)(-.1,.6)
\psline[border=.2,linearc=.3](-.9,-.1)(1.1,-.1)(1.1,.6)
\psline(-1.5,-.1)(-1.15,-.1)
\psline[border=.2,linearc=.3]{->}(-1.5,.1)(-1,.1)(-1,-.4)(-2.5,-.4)
\psline[linearc=.3,linestyle=dashed,dash=2pt 5pt](3,2)(3,2.8)(0.5,3.3)(0.5,.6)
\psline[border=.1,doubleline=true, linearc=.3,doublesep=.2](-1.5,0)(-2.5,0)(-2.5,3)(2.5,2.5)(2.5,1.5)(1.5,1.5)(-1,2.5)(-2,2.5)(-2,1.5)(1.5,1)(3,1)(3,2)
\psline[linearc=.3,linestyle=dashed,dash=2pt 5pt](0.5,.6)(.4,.3)(-.1,-.1)(-.5,-.3)(-1,-.1)(-1.5,0)(-2.5,0)(-2.5,3)(2.5,2.5)(2.5,1.5)(1.5,1.5)(-1,2.5)(-2,2.5)(-2,1.5)(1.5,1)(3,1)(3,2)
\psdot(.4,.3)
\psarc{<-}(2,2){.3}{-180}{180}
\psarc{<-}(-1.5,2){.3}{-180}{180}
\psarc{<-}(2,-.7){.3}{-180}{180}
\rput[l](-1.7,2){\small $c$}
\rput[l](1.8,2){\small $b$}
\rput[l](1.8,-.7){\small $a$}
\rput[lt](-2.3,-.5){\small $T$ after surgery}
\end{pspicture}
\end{center}
Therefore by Lemma~\ref{llll}, this path is obtained from $T$ by surgery
on a $Y$\!--graph with leaves $d$, $a(1)$ and $f$, where $a(1)$ and $f$ are the two dashed handle cores of the genus one surface, $a(1)$ is homotopic
to $\beta\alpha\beta^{-1}$, and $f$ is homotopic to $[\beta,\gamma]$. Thus, $f$ is obtained from a trivial leaf $f_0$
by surgery on a $Y$\!--graph, with one trivial leaf that makes a Hopf link together with $f_0$, and two other leaves $b(1)$
and $c(1)$
that do not link $a(1)$ and that are homotopic to $\beta$ and $\gamma^{-1}$, respectively, as Lemma~\ref{llll} and the next picture show.
\begin{center}
\begin{pspicture}[.4](-1.6,-.1)(2,3.1) 
\rput(1.3,.1){\small $J_1$}
\psline{-*}(0,.6)(0,.9)
\rput[l](-.4,2.8){\small $c(1)$}
\psarc{->}(-.8,2.8){.3}{-180}{0}
\psarc[linestyle=dashed,dash=2pt 2pt](-.8,2.8){.3}{0}{180}
\rput[l](1.2,2.8){\small $b(1)$}
\psarc{->}(.8,2.8){.3}{-180}{0}
\psarc[linestyle=dashed,dash=2pt 2pt](.8,2.8){.3}{0}{180}
\rput[l](.4,.3){\small $a(1)$}
\psarc[linestyle=dashed,dash=2pt 2pt]{->}(0,.3){.3}{-180}{0}
\psarc(0,.3){.3}{0}{180}
\rput[b](-.9,1.3){\small $d$}
\psarc{->}(-.9,.9){.3}{-90}{90}
\psarc[linestyle=dashed,dash=2pt 2pt](-.9,.9){.3}{90}{-90}
\psline{-*}(-.8,2.5)(0,2.3)
\psline{-}(.8,2.5)(0,2.3)
\psline{-}(0,1.9)(0,2.3)
\psline{-}(0,.9)(0,1.3)
\psline(-.6,.9)(0,.9)
\psarc{->}(0,1.5){.2}{60}{0}
\psarc{->}(0,1.7){.2}{-120}{180}
\rput[l](.3,1.5){\small $f_0$}
\end{pspicture}
\begin{pspicture}[.4](-4,0)(4,3.3) 
\psline[linearc=.3](3,2)(3,2.8)(0.5,3.3)(0.5,2.85)
\psline(0.5,2.65)(0.5,2.26)
\psline(0.5,2.06)(0.5,1.3)
\psline[linearc=.3](-.5,.8)(.5,.8)(.5,1.1)
\psline[linearc=.3]{<-}(-.5,.8)(-1.5,.8)(-1.5,3)(2.5,2.5)(2.5,1.5)(1.5,1.5)(0,2.5)(-1,2.5)(-1,1.5)(1.5,1)(3,1)(3,2)
\psline[linearc=.3,linestyle=dotted]{->}(.5,2.5)(.5,3.1)(2.8,2.6)(2.8,1.9)
\psline[linearc=.3,linestyle=dotted](2.8,1.9)(2.8,1.2)(.5,1.4)(.5,2)
\psline[linearc=.3,linestyle=dashed,dash=2pt 5pt]{->}(-1.3,1.9)(-1.3,1.56)(1,1.1)(1.2,1.7)(-.3,2.7)(-1.3,2.7)(-1.3,1.9)
\psarc{<-}(2,2){.3}{-180}{180}
\psarc{<-}(-.5,2){.3}{-180}{180}
\rput[l](-.7,2){\small $c$}
\rput[l](1.8,2){\small $b$}
\rput[l](3.1,1.9){\small $b(1)$}
\rput[r](-1.6,1.9){\small $c^+(1)$}
\rput[t](-.5,.7){\small $f$}
\end{pspicture}
\end{center}
Similarly, for $i=2$ or $3$, there exists a two-component $Y$\!--link $J_i$, that is obtained from $I_i$ by changing (the edge adjacent to $d$ and) the leaves $a$, $b$ and $c$ into leaves 
$a(i)$, $b(i)$ and $c(i)$ that are homologous to $\alpha(i)a$, $\beta(i)b$ and $c$ 
 in $\left(M \setminus (G_1 \setminus I_1)\right)$, respectively, 
where $\alpha(i),\beta(i)\in \{-1,1\}$, and $\alpha(i)\beta(i)=-1$, such that surgery on $J_i$ makes
a pack $T$ of surgery arcs in a surgered disk bounded by $d$
describe the element of $\pi_1(D_{abc})$, $[\gamma\beta\gamma^{-1},[\alpha,\gamma]]$ for $i=2$, or  $[\alpha\gamma\alpha^{-1},[\beta,\alpha]]$ for $i=3$, in $D_{abc} \times [0,1]$ in an ascending way with respect to the height of $[0,1]$.

In particular, the following identity in the free group generated by $\alpha$, $\beta$ and $\gamma$ -whose verification is straightforward-
$$[\beta\alpha\beta^{-1},[\gamma,\beta]][\gamma\beta\gamma^{-1},[\alpha,\gamma]][\alpha\gamma\alpha^{-1},[\beta,\alpha]]=1$$
ensures that the surgery on $J_1 \cup J_2 \cup J_3$ is trivial.

Therefore, if $H_i$ is obtained from $G_i$ by changing
$I_i$ into $J_i$,
$$\begin{array}{l}{[M_{J_1 \cup J_2},H_3]} + {[M_{J_1},H_2]}+ {[M,H_1]}\\
=
\sum_{J \subset (G_1 \setminus I_1)}(-1)^{\sharp J}
\left( \begin{array}{l}M_{J \cup J_1 \cup J_2 \cup J_3}-M_{J \cup J_1 \cup J_2}\\
+M_{J \cup J_1 \cup J_2}-M_{J \cup J_1}\\
+M_{J \cup J_1} -M_{J}\end{array}\right)
=0.\end{array}$$

Furthermore, Proposition~\ref{Review} ensures that
$$\overline{[(M;G_i)]}=-\overline{[(M;H_i)]}=
-\overline{[(M_{\cup_{j<i}J_j};H_i)]},$$
and allows us to conclude the proof.
\eop

This proof shows how the Jacobi IHX relation comes from the Lie algebra structure on the graded space associated to the lower central series of a free group. See \cite{MKS}.

\Addresses\recd


\begin{thebibliography}{AAA}


\bibitem[GGP]{ggp} {\bf S Garoufalidis}, {\bf M Goussarov}, {\bf M Polyak},
{\em Calculus of clovers and finite type invariants of 3--manifolds},
\gtref5{2001}{3}{75}{108}
  \MR{1812435}


\bibitem[GM]{gm} {\bf L Guillou}, {\bf A Marin}, {\em Notes sur l'invariant
de Casson des sph\`eres d'homologie de dimension trois},
Enseign. Math. {38} (1992)  233--290
  \MR{1189008}


\bibitem[Hbg]{hbg} {\bf N Habegger}, {\em Milnor, Johnson, and Tree Level 
Perturbative Invariants}, preprint (2000)
\url{http://www.math.sciences.univ-nantes.fr/~habegger}

\bibitem[Hbo]{hab} {\bf K Habiro}, {\em Claspers and finite type invariants 
of links}, \gtref4{2000}{1}{1}{83}
  \MR{1735632}

\bibitem[KT]{KT}
{\bf  G Kuperberg}, {\bf D\,P Thurston},
{\em Perturbative 3--manifold invariants by cut-and-paste topology}, 
\arxiv{math.GT/9912167}


\bibitem[Le]{le} {\bf T\,T\,Q Le}, {\em An invariant of integral
homology $3$--spheres which is universal for all finite type
invariants}, from: ``Solitons, geometry, and topology: on the crossroad'',
Amer. Math. Soc. Transl. 179 (1997) 75--100
\MR{1437158}

\bibitem[Les]{les}
{\bf C Lescop}, {\em A sum formula for the Casson--Walker invariant},
Invent. Math. 133 (1998) 613--681
  \MR{1645066}

\bibitem[L2]{L2} {\bf C Lescop},
{\em Splitting formulae for the Kontsevich--Kuperberg--Thurston invariant of rational homology $3$--spheres}, \arxiv{math.GT/0411431}

\bibitem[Lic]{Lic} {\bf W\,B\,R Lickorish}, {\em An Introduction to Knot 
Theory}, Graduate Texts in Mathematics 175,
Springer--Verlag, New York (1997)
  \MR{1472978}

\bibitem[MKS]{MKS} {\bf W Magnus}, {\bf A Karrass}, {\bf D Solitar},
   {\em Combinatorial group theory. Presentations of groups in terms
   of generators and relations}, Second revised edition, Dover
   Publications, Inc. New York (1976) \MR{0422434}


\bibitem[Mat]{mat} {\bf S\,V Matveev}, 
{\it Generalized surgery of three-dimensional manifolds and representations of homology spheres}, 
Mat. Zametki  {42} (1987) 268--278, 345 (Russian);
English translation: Math. Notes {42} (1987) 651--656  
  \MR{0915115}

\bibitem[MN]{mn} {\bf H Murakami}, {\bf Y Nakanishi}, {\em On a certain move 
generating link-homology}, 
Math. Ann. {284} (1989) 75--89 
  \MR{0995383}


\bibitem[Oht]{oht} {\bf T Ohtsuki}, {\em Finite type invariants of integral homology $3$--spheres}, 
J. Knot Theory 
Ramifications {5} (1996) 101--115
  \MR{1373813}


\bibitem[Rol]{Rol} {\bf D Rolfsen}, {\em Knots and Links,\/}
Mathematics Lecture Series, No.7. 
Publish or Perish, Inc. Berkeley, Calif. (1976)
\MR{0515288}

\bibitem[Thu]{Thu} {\bf W\,P Thurston}, 
{\em The Geometry and Topology of $3$--manifolds, Chapter $13$},
Princeton University (1978)\nl
\url{http://www.msri.org/publications/books/gt3m/}

\end{thebibliography}
\end{document}